\documentclass[12pt]{amsart}
\usepackage{amssymb,amscd,amsthm,amsmath,color}
\usepackage{fullpage}
\usepackage{epsfig}
\usepackage{graphicx}
\usepackage{xypic}
\usepackage{url}
\usepackage{color}
\usepackage{hyperref}
\usepackage{enumerate}

\numberwithin{equation}{section}

\newcommand{\PP}{\mathbb{P}}
\newcommand{\OO}{\mathcal{O}}

\newcommand{\Spec}{\operatorname{Spec}}
\newcommand{\codim}{\operatorname{codim}}

\newcommand{\rk}{\operatorname{rk}}

\newcommand{\Eff}{\overline{\operatorname{Eff}}}
\newcommand{\Nef}{\operatorname{Nef}}
\newcommand{\Supp}{\operatorname{Supp}}
\newcommand{\Mor}{\operatorname{Mor}}

\newcommand{\SP}{\operatorname{SP}}

\newcommand{\Sec}{\operatorname{Sec}}

\newcommand{\lct}{\operatorname{lct}}
\newcommand{\ord}{\operatorname{ord}}
\newcommand{\Fit}{\operatorname{Fit}}

\newtheorem{theorem}{Theorem}[section]
\newtheorem{lemma}[theorem]{Lemma}
\newtheorem{proposition}[theorem]{Proposition}
\newtheorem{corollary}[theorem]{Corollary}

\theoremstyle{definition}
\newtheorem{definition}[theorem]{Definition}
\newtheorem{notation}[theorem]{Notation}

\newtheorem{remark}[theorem]{Remark}

\newtheorem{question}[theorem]{Question}
\newtheorem{example}[theorem]{Example}
\newtheorem{construction}[theorem]{Construction}
\newtheorem{motivatingquestion}[theorem]{Motivating Question}

\input xy
\xyoption{all}

\begin{document}

\title{Codimension of jumping loci}

\author{Brian Lehmann}
\address{Department of Mathematics \\
Boston College  \\
Chestnut Hill, MA \, \, 02467}
\email{lehmannb@bc.edu}

\author{Eric Riedl}
\address{Department of Mathematics \\
University of Notre Dame  \\
255 Hurley Hall \\
Notre Dame, IN 46556}
\email{eriedl@nd.edu}

\author{Sho Tanimoto}
\address{Graduate School of Mathematics, Nagoya University, Furocho Chikusa-ku, Nagoya, 464-8602, Japan}
\email{sho.tanimoto@math.nagoya-u.ac.jp}

\begin{abstract}
Suppose that $\mathcal{E}$ is a vector bundle on a smooth projective variety $X$.  Given a family of curves $C$ on $X$, we study how the Harder-Narasimhan filtration of $\mathcal{E}|_{C}$ changes as we vary $C$ in our family.  Heuristically we expect that the locus where the slopes in the Harder-Narasimhan filtration jump by $\mu$ should have codimension which depends linearly on $\mu$.  We identify the geometric properties which determine whether or not this expected behavior holds. We then apply our results to study rank $2$ bundles on $\mathbb{P}^{2}$ and to study singular loci of moduli spaces of curves.
\end{abstract}

\maketitle

\section{Introduction}

Let $X$ be a smooth projective complex variety and let $\mathcal{E}$ be a locally free coherent sheaf on $X$.  Given a dominant family of curves $C$ on $X$, we analyze how the Harder-Narasimhan filtrations of the restrictions $\mathcal{E}|_{C}$ change as we vary $C$.  Results which bound the behavior of $\mathcal{E}|_{C}$ are known as ``restriction theorems''; such theorems have played a key role in the theory of slope stability since the work of \cite{Barth77} and \cite{Hulek79} on jumping lines in $\mathbb{P}^{n}$.  (In this paper ``stability'' will always refer to ``slope stability with respect to a nef curve class'' as first developed by \cite{CP11}; see Notation \ref{nota:hnfilt}.)

Most restriction theorems -- such as the Mehta-Ramanathan theorem (\cite{MR82,MR84}) or the Grauert-M\" ulich theorem  (\cite{GM75, Spindler79,FHS80,Maruyama81,PatelRiedlTseng,LRT23}) -- address the case when $C$ is general in its parameter space.  Our focus will be on non-general curves $C$:

\begin{motivatingquestion} \label{motques}
Suppose $W \subset \mathcal{M}_{g,0}(X,\beta)$ parametrizes a dominant family of curves $s: C \to X$.  Can we bound the codimension of $W$ in terms of the difference in slopes between the Harder-Narasimhan filtration of $s^* \mathcal{E}$ and the ``expected value'' given by the Harder-Narasimhan filtration of $\mathcal{E}$ with respect to the numerical class $s_{*}[C]$?
\end{motivatingquestion}

In this paper we will analyze Motivating Question \ref{motques} using techniques from stability theory.  Throughout we aim for the maximal possible generality: our main theorem will describe the qualitative behavior of restrictions of vector bundles to dominant families of curves. Loosely speaking, we show that when the evaluation map over $W$ has connected fibers then the codimension of $W$ is bounded by a linear function in the difference in slopes for the Harder-Narasimhan filtrations.  Conversely, when the evaluation map does not have connected fibers one cannot always obtain a linear bound.

Motivating Question \ref{motques} has found many applications in the literature; we present two in this paper.  First, we use our work to bound the codimension of the singular locus of the moduli spaces of curves.  The same argument establishes geometric analogues of arithmetic conjectures of \cite{Peyre17}.  Second, we analyze jumping loci for rank 2 bundles on $\mathbb{P}^2$, an important classical subject that continues to play a role in modern research (for example in Terao's Conjecture or the Weak Lefschetz Property for Artinian algebras). Our results explain when the ideal jumping behavior does or does not happen based on the geometry of the families of jumping curves.

\subsection{Main result}

We will need to precisely quantify the difference in slopes between two filtrations of a vector bundle.  We will use the following definition:

\begin{definition} \label{defi:slopepanel}
Let $X$ be a smooth projective variety and let $\alpha \in \Nef_{1}(X)$ be a nef curve class.  Suppose that $\mathcal{E}$ is a non-zero torsion-free sheaf of rank $r$.  Write
\begin{equation*}
0 = \mathcal{F}_{0} \subset \mathcal{F}_{1} \subset \mathcal{F}_{2} \subset \ldots \subset \mathcal{F}_{s} = \mathcal{E}
\end{equation*}
for the $\alpha$-Harder-Narasimhan filtration of $\mathcal{E}$. 
The slope panel $\SP_{X,\alpha}(\mathcal{E})$ is the $r$-tuple of rational numbers consisting of the union over all indices $i$ of $\rk(\mathcal{F}_{i}/\mathcal{F}_{i-1})$ copies of $\mu_{\alpha}(\mathcal{F}_{i}/\mathcal{F}_{i-1})$ arranged in non-increasing order:
\begin{equation*}
\SP_{X,\alpha}(\mathcal{E}) = ( \underbrace{\mu_{\alpha}(\mathcal{F}_{1}/\mathcal{F}_{0}), \ldots}_{\rk(\mathcal{F}_{1}/\mathcal{F}_{0}) \textrm{ copies}}, \underbrace{\mu_{\alpha}(\mathcal{F}_{2}/\mathcal{F}_{1}), \ldots}_{\rk(\mathcal{F}_{2}/\mathcal{F}_{1}) \textrm{ copies}},\ldots,\underbrace{\mu_{\alpha}(\mathcal{F}_{s}/\mathcal{F}_{s-1}), \ldots}_{\rk(\mathcal{F}_{s}/\mathcal{F}_{s-1}) \textrm{ copies}})
\end{equation*}
For a curve $C$, we write $\SP_{C}(\mathcal{E})$ for the slope panel with respect to the class of a point on $C$.
\end{definition}

Returning to our motivating question, suppose that $U \to W$ is a family of curves equipped with a dominant evaluation morphism $ev: U \to X$.  \cite{LRT23} shows that the global slope panel $\SP_{X,s_{*}[C]}(\mathcal{E})$ is the ``expected value'' of the slope panel for the restriction of $\mathcal{E}$ to a curve.  Suppose that the slope panel for the general curve $s: C \to X$ parametrized by $W$ differs from the expected value by $\mu$, that is,
\begin{equation*}
\Vert \SP_{X,s_{*}[C]}(\mathcal{E}) - \SP_{C}(s^{*}\mathcal{E}) \Vert_{sup} = \mu.
\end{equation*}
If $W$ parametrizes every curve whose difference in slope panels is $\mu$, then basic deformation theory provides an upper bound on the codimension of $W$ that is linear in $\mu$.  In some situations we can also obtain lower bounds on the codimension of $W$ that are linear in $\mu$ (e.g.~when the map from $W$ to the versal deformation space of $s^{*}\mathcal{E}$ is smooth near $s$). 

Our main theorem identifies a set of mild assumptions which guarantee that there is indeed a lower bound on the codimension of $W$ that is linear in $\mu$.  It applies to all curves whose classes are contained in any slightly smaller subcone $\mathcal{C} \subset \Nef_{1}(X)$.

\begin{theorem} \label{theo:maintheorem1}
Let $X$ be a smooth projective variety and let $\mathcal{E}$ be a non-zero locally free coherent sheaf on $X$.  Fix a genus $g$, an ample divisor $H$ on $X$, and a closed cone $\mathcal{C} \subset N_{1}(X)_{\mathbb{R}}$ such that $\mathcal{C} \backslash \{0\}$ is contained in the interior of $\Nef_{1}(X)$.  Then there are affine linear functions $S: \mathbb{R} \to \mathbb{R}$ and $L: \mathbb{R} \to \mathbb{R}$ with positive leading coefficients which have the following property.

Let $W \to \mathcal{M}_{g,0}(X)$ be a generically finite morphism from an irreducible variety and let $ev^{\nu}: U^{\nu} \to X$ denote the evaluation map on the normalization $U^{\nu}$ of the universal family over $W$.  Define
\begin{equation*}
\mu = \Vert \SP_{X,s_{*}[C]}(\mathcal{E}) - \SP_{C}(s^{*}\mathcal{E}) \Vert_{sup}
\end{equation*}
where $s: C \to X$ is a general curve parametrized by $W$.
Suppose that:
\begin{itemize}
\item $ev^{\nu}$ is dominant with connected fibers,
\item the general map parametrized by $W$ is birational onto its image, and
\item the class of the curves parametrized by $W$ is contained in the cone $\mathcal{C}$. 
\end{itemize}
Then the codimension of the image of $W$ satisfies either
\begin{enumerate}
\item $\codim(W) \geq S(\mu)$, or
\item $\codim(W) \geq L(H \cdot s_{*}C)$.
\end{enumerate}
\end{theorem}

\begin{remark} \label{rema:functiondependence}
The functions $S, L$ depend on $X$, $g(C)$, the sheaf $\mathcal{E}$, the cone $\mathcal{C}$, and the ample divisor $H$, as well as two choices: the choice of a fixed curve class (called $\beta$ in the proof) contained in $\mathcal{C}$ and the choice of a direct sum of line bundles that surjects onto $\mathcal{E}$.
\end{remark}

The second bulleted assumption in Theorem \ref{theo:maintheorem1} is natural: if every morphism $s$ parametrized by $W$ factors through an intermediate curve $s': C' \to X$ then we should analyze $s^{*}\mathcal{E}$ by first understanding $s'^{*}\mathcal{E}$.  The third bulleted assumption is used in our proof but we are not sure whether it is necessary.  The first bulleted assumption is essential as explained in the following remark.

\begin{remark} \label{rema:connfibers}
In Theorem \ref{theo:maintheorem1} it is necessary to assume that $ev^{\nu}$ has connected fibers.  If $ev^{\nu}$ fails to satisfy this property then by using Stein factorizations one can construct a generically finite morphism $f: Y \to X$ such that $ev^{\nu}$ factors rationally through $f$.  We can only expect -- and indeed Theorem \ref{theo:maintheorem1} proves -- that $s^{*}\mathcal{E}$ is controlled by the stability of $f^{*}\mathcal{E}$ on $Y$ instead of the stability of $\mathcal{E}$ on $X$.  (Note that there are typically many different curve classes $\alpha'$ on $Y$ which pushforward to a given curve class $\alpha$ on $X$, so the Harder-Narasimhan filtration of $f^* \mathcal{E}$ can be very different from that of $\mathcal{E}$ depending on $\alpha'$.)  See Example \ref{exam:lines} for a demonstration of this phenomenon.

When $X$ is a Fano variety it is possible to understand dominant families of curves whose evaluation map has disconnected fibers using the theory of accumulating maps developed by \cite{LRT23}.  In Theorem \ref{theo:mainfanocodim} we combine both approaches to prove a more comprehensive result for Fano varieties.
\end{remark}

\subsection{Applications}

\subsubsection{Rank 2 bundles on $\mathbb{P}^{2}$}
Since Theorem \ref{theo:maintheorem1} is aiming for maximal generality, it is not reasonable to expect our bounds to be sharp.  However, in specific examples one can emulate our approach with more careful estimates to obtain precise information.  In Section \ref{sect:r2p2} we carry out this program for rank $2$ bundles on $\mathbb{P}^{2}$.

Mathematicians have studied the jumping behavior for many special types of Fano varieties, including Fano varieties ruled by lines (see \cite{MOS12} and the references therein).  However the most well-studied examples are vector bundles on $\PP^n$, especially rank 2 vector bundles on $\PP^2$. There is the classical theory of jumping lines (which has too many applications to list here) and also work on the jumping locus of higher degree curves (\cite{Stromme84}, \cite{Manaresi90}, \cite{Ran01}, \cite{Vitter04}, \cite{Marangone24}). Much of the previous work has been example-based, studying the jumping loci for particular vector bundles and observing whether they did or did not have the expected dimension. 

Our work identifies precisely the geometric obstructions to obtaining the expected behavior for special families of curves on $\mathbb{P}^{2}$. 
The cleanest result is for rational curves: we show that families of rational curves $W$ have the expected codimension unless there is a clear geometric reason why they cannot (see cases (a) and (b) in Theorem \ref{theo:introp2} below). Furthermore, we show these exceptional geometries only occur when the codimension of $W$ is large relative to the degree.

\begin{theorem} \label{theo:introp2}
Let $\mathcal{E}$ be a stable rank $2$ bundle on $\mathbb{P}^{2}$.  There is an explicit constant $\zeta$ depending on the Chern classes of $\mathcal{E}$ such that the following holds.

Fix a constant $\mu$ and suppose that $V \subset \mathcal{M}_{0,0}(\mathbb{P}^{2},d)$ is an irreducible component of the locus of degree $d$ maps $s: C \to \mathbb{P}^{2}$ satisfying
\begin{equation*}
\Vert \SP_{\mathbb{P}^{2},s_{*}[C]}(\mathcal{E}) - \SP_{C}(s^{*}\mathcal{E}) \Vert_{sup} = \mu.
\end{equation*}
Let $W \to V$ be a generically finite dominant morphism from a variety.

Assume that the general map $s: C \to \mathbb{P}^{2}$ parametrized by $W$ is a birational immersion.  Then either
\begin{enumerate}
\item $\codim(W)$ has the expected value $\sup\{ 2\mu - 1, 0\}$, or
 \item $\codim(W) \geq \zeta d$.
\end{enumerate}
Furthermore, the only way (1) can fail is when either (a) the evaluation map for the normalization of the universal family over $W$ fails to have connected fibers or (b) there is a birational model $\phi: X' \to \mathbb{P}^{2}$ flattening the family of curves $V$ such that $\phi^{*}\mathcal{E}$ fails to be semistable with respect to the strict transforms of the curves.
\end{theorem}

This result systematizes and generalizes various examples in the literature which demonstrate that Theorem \ref{theo:introp2}.(1) can fail (see Example \ref{exam:lines} and Example \ref{exam:conics}).

\begin{remark}
In a different direction, \cite{Bogomolov95} and \cite{Kopper20} show that the restriction of a stable bundle to an embedded smooth curve of sufficiently high degree in a surface will always be stable.  In contrast, the pullback of a bundle to the normalization of an embedded singular curve is more complicated and cannot admit such a simple description.  Note that if we fix a genus $g$ then most plane curves of geometric genus $g$ will be singular.  
\end{remark}

\subsubsection{Locus of non-free curves}
When studying the moduli space of curves $\Mor(C,X)$ on a Fano variety $X$, it is crucial to develop a good understanding of the singular locus of the moduli space.  For example, singularities can affect the computation of the Kodaira dimension of $\Mor(C,X)$ (as in \cite{Starr03, dJS17}).  Since a curve $s: C \to X$ can only be contained in the singular locus if $H^{1}(C,s^{*}T_{X}) \neq 0$, we can bound the codimension of the singular locus of $X$ by studying the Harder-Narasimhan filtration of the restricted tangent bundle.

More generally, recall that a curve $s: C \to X$ is said to be $m$-free if $\mu^{min}(s^{*}T_{X}) \geq m+2g$.  (The singular locus of $\Mor(C,X)$ is a subset of the curves which fail to be $0$-free.)  For certain types of variety the locus of $m$-free curves has been analyzed in detail, including projective space (\cite{Ascenzi88}, \cite{Ramella}, \cite{HK93}, \cite{BR97}, \cite{Hein00}, \cite{Ran01}, \cite{GHI13}, \cite{AR15}, \cite{Larson16}, \cite{CR18}, \cite{Ascenzi22}), Grassmannians (\cite{BR00}, \cite{Mandal}), and hypersurfaces (\cite{BrowningSawin}).  Our work allows us to bound the dimension of the locus of curves which fail to be $m$-free for arbitrary Fano varieties.

Our results are motivated by analogous questions in arithmetic geometry.  In \cite{Peyre17} Peyre makes a number of compelling conjectures about the behavior of arithmetic slopes of rational points on Fano varieties.  More precisely, he gives a variant of the counting problem in Manin's Conjecture based on slopes and conjectures this modified counting problem still has the expected asymptotic behavior.  Loosely speaking, Peyre's conjectures suggest that ``most'' rational points should have large slopes.

Analogously, we can expect that for ``most'' curves on a Fano variety the restricted tangent bundle has large minimal slope.  As usual in Manin's Conjecture, we must remove an ``exceptional set'' of curves which have pathological behavior.  For this purpose we will use the theory of accumulating maps as developed by \cite{LRT23} (see Definition \ref{defi:accumulatingmap}).  By combining Theorem \ref{theo:maintheorem1} with the results of \cite{LRT23}, we show that the codimension of the locus of curves which fails to be $m$-free grows linearly in the degree \emph{unless} the curves come from an accumulating map.  This establishes a geometric analogue of the arithmetic conjectures of \cite{Peyre17}.

\begin{theorem} \label{theo:intrononfree}
Let $X$ be a smooth Fano variety.  Fix a genus $g$, a positive constant $m$, and a closed cone $\mathcal{C}$ such that $\mathcal{C} \backslash \{0\}$ is contained in the interior of $\Nef_{1}(X)$.  There is an affine linear function $T: \mathbb{R} \to \mathbb{R}$ with a positive leading coefficient which satisfies the following property.

Let $W$ be a variety equipped with a generically finite morphism $W \to {\mathcal M}_{g,0}(X)$ and let $ev^{\nu}: U^{\nu} \to X$ denote the evaluation map on the normalization $U^{\nu}$ of the universal family over $W$.
Assume that the class of the curves parametrized by $W$ is contained in the cone $\mathcal{C}$.

If the general morphism $s: C \to X$ parametrized by $W$ fails to be $m$-free (as in Definition \ref{defi:mfree}), then one of the following properties holds:
\begin{enumerate}
\item $\codim(W) \geq T(-K_{X} \cdot s_{*}C)$,
\item the evaluation map $ev^{\nu}$ factors rationally through an accumulating morphism $f: Y \to X$, or
\item the image of the general map $s: C \to X$ parametrized by $W$ is a rational curve of anticanonical degree $\leq 2$.
\end{enumerate}
\end{theorem}

\begin{remark}
The function $T$ depends on all the choices listed in Remark \ref{rema:functiondependence} as well as the finite set of slope panels of $T_{X}$ with respect to nef curve classes on $X$ and the (inexplicit but finite) list of possible Fujita invariants of value $>\frac{1}{2}$ for pairs of dimension $\leq \dim(X)$ arising from Birkar's solution to the Borisov-Alexeev-Borisov conjecture as in \cite{HL20}.
\end{remark}

\subsection{Strategy}

Retaining the notation of Theorem \ref{theo:maintheorem1}, suppose that we have a generically finite morphism $W \to \mathcal{M}_{g,0}(X)$ such that $ev^{\nu}$ is dominant with connected fibers.  The best situation is when the evaluation morphism satisfies an additional property: the general curve $s: C \to X$ in our family is contained in the locus where $ev^{\nu}$ is flat.  In this case, the Grauert-M\" ulich theorem of \cite{LRT23} gives excellent control on the codimension of $W$.  When this condition does not apply then we choose a birational model $\phi: X' \to X$ such that the strict transforms of the curves parametrized by $W$ define a flat family of curves on $X'$.  Letting $s': C' \to X'$ denote a strict transform curve, the Grauert-M\" ulich theorem relates $s^{*}\mathcal{E}$ to the $s'_{*}[C']$-Harder-Narasimhan filtration of $\phi^{*}\mathcal{E}$ on $X'$.

The key challenge is to understand how the stability of $\mathcal{E}$ changes when we pass from the class $s_{*}[C]$ on $X$ to the class $s'_{*}[C']$ on $X'$.  Although stability is preserved by \emph{pullback} of curve classes, it usually changes under the \emph{strict transform} of curve classes.  To obtain Theorem \ref{theo:maintheorem1} in full generality, we must develop a systematic theory of how stability changes under birational maps.  This study makes up the technical core of the paper.

We leverage the following tension.  Suppose that $\mathcal{Q}$ is a rank $r$ quotient of $\mathcal{E}$ such that the induced map $\phi^{*}\mathcal{E} \to (\phi^{*}\mathcal{Q})_{tf}$  is $\alpha'$-destabilizing for some class $\alpha' \in \Nef_{1}(X')$ which pushes forward to $\alpha \in \Nef_{1}(X)$.  The difference between $\mu_{\alpha}(\mathcal{Q})$ and $\mu_{\alpha'}((\phi^{*}\mathcal{Q})_{tf})$ is controlled by the $r$th Fitting ideal of $\mathcal{Q}$.  If $\mathcal{Q}$ has small slope compared to $\mathcal{E}$, then the Fitting ideal of $\mathcal{Q}$ must be comparatively simple so that the difference in slopes will be small.  If $\mathcal{Q}$ has large slope compared to $\mathcal{E}$, then the Fitting ideal can be more complicated but also a larger change in slope is required for $\mathcal{Q}$ to become destabilizing.

We prove Theorem \ref{theo:maintheorem1} by bounding the log canonical threshold of the $r$th Fitting ideal of $\mathcal{Q}$.  Using the log-canonical threshold, we can relate the relative anticanonical degree of $\alpha'$ with the change in $\alpha'$-slope induced by the Fitting ideal.  In turn, the relative anticanonical degree allows us to control the codimension of the corresponding family of curves.

\

\noindent
{\bf Acknowledgements:}
The authors would like to thank the referee for detailed suggestions which significantly improved the exposition of the paper.
Part of this project was conducted at the SQuaRE workshop ``Geometric Manin's Conjecture in characteristic $p$'' at the American Institute of Mathematics. Furthermore this project was supported by the National Science Foundation under Grant No. DMS-1929284 while the authors were in residence at the Institute for Computational and Experimental Research in Mathematics in Providence, RI, during the Counting Sections of Fano Fibrations Collaborate@ICERM program. The authors would like to thank AIM and ICERM for the generous support.

Brian Lehmann was supported by Simons Foundation grant Award Number 851129.
Eric Riedl was supported by NSF CAREER grant DMS-1945944.  Sho Tanimoto was partially supported by JST FOREST program Grant number JPMJFR212Z, by JSPS Bilateral Joint Research Projects Grant number JPJSBP120219935, and by JSPS KAKENHI Grand-in-Aid (B) 23K25764.

\section{Background}

We work over the ground field $\mathbb{C}$.  Throughout all our schemes will be assumed to be separated schemes and every connected component will have finite type over the ground field.  A variety is a scheme that is reduced and irreducible.  Given a coherent sheaf $\mathcal{F}$ on a variety $V$, we denote by $\mathcal{F}_{tors}$ the torsion subsheaf of $\mathcal{F}$ and by $\mathcal{F}_{tf}$ the quotient of $\mathcal{F}$ by its torsion subsheaf.  We say that a morphism is generically finite if it is generically finite onto its image.

When $X$ is a projective variety, we let $N^{1}(X)_{\mathbb{R}}$ denote the space of $\mathbb{R}$-Cartier divisors up to numerical equivalence and $N_{1}(X)_{\mathbb{R}}$ denote the dual space of $\mathbb{R}$-curves up to numerical equivalence.  We let $N^{1}(X)_{\mathbb{Z}}$ and $N_{1}(X)_{\mathbb{Z}}$ denote the lattices of $\mathbb{Z}$-classes in these two vector spaces.  We define $\Nef_{1}(X)$ to be the nef cone of curves, i.e.~the set of $\mathbb{R}$-curve classes $\alpha$ which satisfy $E \cdot \alpha \geq 0$ for every effective divisor $E$.  Given a curve $C$, we denote its numerical class by $[C]$.  We also define $\Eff^{1}(X)$ to be the pseudo-effective cone of $\mathbb{R}$-Cartier divisors.

We let $\overline{\mathcal{M}}_{g,n}(X)$ denote the Kontsevich moduli stack of stable maps and let $\mathcal{M}_{g,n}(X) \subset \overline{\mathcal{M}}_{g,n}(X)$ denote the open substack parametrizing maps with smooth irreducible domain.

\subsection{Slope stability}

Let $X$ be a smooth projective variety and let $\mathcal{E}$ be a torsion-free coherent sheaf on $X$.  We briefly review the notion of slope stability of $\mathcal{E}$ with respect to nef curve classes $\alpha$ as developed in \cite{CP11}.

\begin{definition}
Let $X$ be a smooth projective variety and let $\mathcal{E}$ be a torsion-free coherent sheaf on $X$ of rank $r$.  The first Chern class of $\mathcal{E}$, denoted by $c_{1}(\mathcal{E})$, is any divisor representing the line bundle
\begin{equation*}
\bigwedge^{[r]} \mathcal{E} := \left( \bigwedge^{r} \mathcal{E} \right)^{\vee \vee}.
\end{equation*}
Note that $c_{1}(\mathcal{E})$ is only well-defined up to linear equivalence; in practice this mild abuse of notation will be harmless.

For any coherent sheaf $\mathcal{E}'$, we can define the first Chern class by taking a finite resolution of $\mathcal{E}'$ by locally free sheaves and defining $c_{1}(\mathcal{E}')$ as an alternating sum of the first Chern classes of these sheaves. Note that these two definitions of the first Chern class are compatible.
\end{definition}

\begin{definition}
Let $X$ be a smooth projective variety and let $\alpha \in \Nef_{1}(X)$.  For any non-zero torsion-free sheaf $\mathcal{E}$ on $X$, we define the $\alpha$-slope
\begin{equation*}
\mu_{\alpha}(\mathcal{E}) = \frac{c_{1}(\mathcal{E}) \cdot \alpha}{\rk(\mathcal{E})}.
\end{equation*}
We say that $\mathcal{E}$ is $\alpha$-semistable if for every non-zero torsion-free subsheaf $\mathcal{F} \subset \mathcal{E}$ we have $\mu_{\alpha}(\mathcal{F}) \leq \mu_{\alpha}(\mathcal{E})$.
\end{definition}

As explained by \cite{CP11,GKP14,GKP16}, this notion of $\alpha$-semistability naturally leads to $\alpha$-Harder-Narasimhan filtrations.

\begin{notation} \label{nota:hnfilt}
Let $X$ be a smooth projective variety and let $\mathcal{E}$ be a non-zero torsion-free coherent sheaf of rank $r$ on $X$.  Fix an $\alpha \in \Nef_{1}(X)$.  The $\alpha$-Harder-Narasimhan filtration of $\mathcal{E}$ is the unique sequence of subsheaves
\begin{equation*}
0 = \mathcal{F}_{0} \subset \mathcal{F}_{1} \subset \mathcal{F}_{2} \subset \ldots \subset \mathcal{F}_{s} = \mathcal{E}.
\end{equation*}
such that each $\mathcal{F}_{i}/\mathcal{F}_{i-1}$ is $\alpha$-semistable and the $\alpha$-slopes of the quotients $\mathcal{F}_{i}/\mathcal{F}_{i-1}$ are strictly decreasing in $i$.

We denote by $\mu^{max}_{\alpha}(\mathcal{E})$ the maximal slope of any torsion-free subsheaf, i.e.,~$\mu^{max}_{\alpha}(\mathcal{E}) = \mu(\mathcal{F}_{1})$.  We denote by $\mu^{min}_{\alpha}(\mathcal{E})$ the minimal slope of any torsion-free quotient, i.e.,~$\mu^{min}_{\alpha}(\mathcal{E}) = \mu(\mathcal{E}/\mathcal{F}_{s-1})$.
\end{notation}

Given a filtration of $\mathcal{E}$ (possibly different from the Harder-Narasimhan filtration), we will need a notation which describes the slopes of the graded pieces of the filtration.  We define the slope panel in the analogous way:

\begin{definition}
Let $X$ be a smooth projective variety, let $\mathcal{E}$ be a non-zero torsion-free coherent sheaf on $X$, and let $\alpha \in \Nef_{1}(X)$.  Suppose that
\begin{equation*}
0 = \mathcal{G}_{0} \subset \mathcal{G}_{1} \subset \ldots \subset \mathcal{G}_{t} = \mathcal{E}
\end{equation*}
is a sequence of distinct saturated subsheaves of $\mathcal{E}$.  We define the slope panel $\SP_{X,\alpha}(\mathcal{E};\mathcal{G}_{\bullet})$ of $\mathcal{E}$ with respect to the sequence $\mathcal{G}_{\bullet}$ as follows.  
\begin{equation*}
\SP_{X,\alpha}(\mathcal{E};\mathcal{G}_{\bullet}) = ( \underbrace{\mu_{\alpha}(\mathcal{G}_{1}/\mathcal{G}_{0}), \ldots}_{\rk(\mathcal{G}_{1}/\mathcal{G}_{0}) \textrm{ copies}}, \underbrace{\mu_{\alpha}(\mathcal{G}_{2}/\mathcal{G}_{1}), \ldots}_{\rk(\mathcal{G}_{2}/\mathcal{G}_{1}) \textrm{ copies}},\ldots,\underbrace{\mu_{\alpha}(\mathcal{G}_{t}/\mathcal{G}_{t-1}), \ldots}_{\rk(\mathcal{G}_{t}/\mathcal{G}_{t-1}) \textrm{ copies}})
\end{equation*}
\end{definition}

Recall from the introduction that the slope panel $\SP_{X,\alpha}(\mathcal{E})$ (with no mention of a filtration) simply means the slope panel of $\mathcal{E}$ with respect to the $\alpha$-Harder-Narasimhan filtration.  Our next goal is Lemma \ref{lemm:hnmaininequality} which gives a combinatorial characterization of the $\alpha$-Harder-Narasimhan filtration.

\begin{lemma} \label{lemm:saturated}
Let $X$ be a smooth projective variety and let $\mathcal{E}$ be a torsion-free coherent sheaf on $X$.  Suppose $\mathcal{F},\mathcal{G} \subset \mathcal{E}$ are saturated subsheaves.  Then $\mathcal{F} \cap \mathcal{G}$ is a saturated subsheaf of $\mathcal{G}$.
\end{lemma}

\begin{proof}
We have an injection $\mathcal{E}/\mathcal{F} \cap \mathcal{G} \to (\mathcal{E}/\mathcal{F}) \oplus (\mathcal{E}/\mathcal{G})$, showing that $\mathcal{E}/\mathcal{F} \cap \mathcal{G}$ is torsion-free.  Thus the subsheaf $\mathcal{G}/\mathcal{F} \cap \mathcal{G}$ is also torsion-free.
\end{proof}

\begin{lemma} \label{lemm:hnmaininequality}
Let $X$ be a smooth projective variety, let $\mathcal{E}$ be a non-zero torsion-free coherent sheaf on $X$ of rank $r$, and let $\alpha \in \Nef_{1}(X)$.  Fix a sequence of distinct saturated subsheaves
\begin{equation*}
0 = \mathcal{G}_{0} \subset \mathcal{G}_{1} \subset \ldots \subset \mathcal{G}_{t} = \mathcal{E}.
\end{equation*}
Set
\begin{align*}
(a_{1},\ldots,a_{r})  & = \SP_{X,\alpha}(\mathcal{E}) \\
(b_{1},\ldots,b_{r}) & = \SP_{X,\alpha}(\mathcal{E};\mathcal{G}_{\bullet})
\end{align*}
For any $1 \leq \ell \leq r$, we have
\begin{equation*}
\sum_{i=1}^{\ell} a_{i} \geq \sum_{i=1}^{\ell} b_{i} \qquad \textrm{and} \qquad \sum_{i=\ell+1}^{r} a_{i} \leq \sum_{i=\ell+1}^{r} b_{i}.
\end{equation*}
\end{lemma}

\begin{proof}
Since $\sum_{i=1}^{r} b_{i} = c_{1}(\mathcal{E}) = \sum_{i=1}^{r} a_{i}$, the second statement is a consequence of the first and so we may focus on the first statement.  Let
\begin{equation*}
0 = \mathcal{F}_{0} \subset \mathcal{F}_{1} \subset \mathcal{F}_{2} \subset \ldots \subset \mathcal{F}_{s} = \mathcal{E}
\end{equation*}
denote the $\alpha$-Harder-Narasimhan filtration of $\mathcal{E}$.  For every index $1 \leq j \leq t$, consider the filtration
\begin{equation*}
0 = \mathcal{G}_{j} \cap \mathcal{F}_{0} \subset \mathcal{G}_{j} \cap \mathcal{F}_{1} \subset \ldots \subset \mathcal{G}_{j} \cap \mathcal{F}_{s} = \mathcal{G}_{j}
\end{equation*}
where some of the entries are possibly equal.  Lemma \ref{lemm:saturated} shows that each $\mathcal{G}_{j} \cap \mathcal{F}_{k}$ is saturated in $\mathcal{G}_{j}$; thus two successive terms will have the same rank if and only if they are equal.  Let $(v^{j}_{\bullet})$ denote the vector of length $s$ whose $k$th entry is $\rk(\mathcal{G}_{j} \cap \mathcal{F}_{k}) - \rk(\mathcal{G}_{j} \cap \mathcal{F}_{k-1})$.  Note that
\begin{equation*} 
\sum_{k=1}^{s} v^{j}_{k} = \rk(\mathcal{G}_{j}) \qquad \textrm{and} \qquad
\sum_{k=1}^{s} v^{j}_{k} \mu_{\alpha}(\mathcal{G}_{j} \cap \mathcal{F}_{k}/\mathcal{G}_{j} \cap \mathcal{F}_{k-1}) = c_{1}(\mathcal{G}_{j}).
\end{equation*}
Next we inductively define the vectors $(w^{j}_{\bullet})$ of length $s$ by setting $(w^{1}_{\bullet}) = (v^{1}_{\bullet})$ and defining $(w^{j}_{\bullet}) = (v^{j}_{\bullet}) - (v^{j-1}_{\bullet})$ for $2 \leq j \leq t$.   We have
\begin{equation} \label{eq:sumformula2}
\sum_{j=1}^{t} w^{j}_{k} = v^{1}_{k} + \sum_{j=2}^{t} (v^{j}_{k} - v^{j-1}_{k}) = v^{t}_{k} = \rk(\mathcal{F}_{k}/\mathcal{F}_{k-1}).
\end{equation}

Finally we inductively construct a vector $(c_{\bullet})$.  Suppose we have constructed the first $\rk(\mathcal{G}_{j-1})$ entries of $(c_{\bullet})$.  We then append the next $\rk(\mathcal{G}_{j}/\mathcal{G}_{j-1})$ entries as follows: for every $1 \leq k \leq s$ we insert $w^{j}_{k}$ copies of $\mu_{\alpha}(\mathcal{F}_{k}/\mathcal{F}_{k-1})$ to $(c_{\bullet})$ and then arrange these $\rk(\mathcal{G}_{j}/\mathcal{G}_{j-1})$ entries in non-increasing order. 

By Equation \eqref{eq:sumformula2} we see that $(c_{\bullet})$ is a rearrangement of $(a_{\bullet})$.  Since $(a_{\bullet})$ is arranged in non-increasing order, for every $1 \leq \ell \leq r$ we have
\begin{equation*}
\sum_{i=1}^{\ell} a_{i} \geq \sum_{i=1}^{\ell} c_{i}.
\end{equation*}
Thus it suffices to prove
\begin{equation*}
\sum_{i=1}^{\ell} c_{i} \geq \sum_{i=1}^{\ell} b_{i}.
\end{equation*}
We first prove this when $\ell = \rk(\mathcal{G}_{j})$ for some $j$.  In this case we have
\begin{align*}
\sum_{i=1}^{\ell} c_{i} & = \sum_{m=1}^{j} \sum_{k=1}^{s} w^{m}_{k} \mu_{\alpha}(\mathcal{F}_{k}/\mathcal{F}_{k-1}) \\
& = \sum_{k=1}^{s} v^{j}_{k} \mu_{\alpha}(\mathcal{F}_{k}/\mathcal{F}_{k-1}) \\
& \geq \sum_{k=1}^{s} v^{j}_{k} \mu_{\alpha}(\mathcal{G}_{j} \cap \mathcal{F}_{k}/\mathcal{G}_{j} \cap \mathcal{F}_{k-1}) \\
& = c_{1}(\mathcal{G}_{j}) \\
& = \sum_{i=1}^{\ell} b_{i}
\end{align*}
where we use the semistability of $\mathcal{F}_{k}/\mathcal{F}_{k-1}$ at the inequality step.  This proves the claim for the special case $\ell = \rk(\mathcal{G}_{j})$.

For other values of $\ell$, we note that for the indices satisfying $\rk(\mathcal{G}_{j-1}) < i \leq \rk(\mathcal{G}_{j})$ the entries of $(b_{\bullet})$ are constant while the entries of $(c_{\bullet})$ are non-increasing.  Since we have the desired inequality when $\ell = \rk(\mathcal{G}_{j-1}), \rk(\mathcal{G}_{j})$ we must also have it for the intermediate values of $\ell$.
\end{proof}

\subsection{Bounded families of sheaves}

\begin{definition} \label{defi:boundedfamily}
Let $X$ be a smooth projective variety.  A finite-type family of torsion-free sheaves on $X$ consists of a finite-type scheme $S$ over our ground field and a coherent sheaf $\mathcal{G}$ on $X \times S$ that is flat over $S$ such that for every point $s \in S$ the restriction $\mathcal{G}|_{X_{s}}$ is torsion-free.

We say that a set of torsion-free sheaves $\{ \mathcal{E}_{i} \}$ on $X$ is a bounded family if there is a finite-type family of torsion-free sheaves such that every $\mathcal{E}_{i}$ is parametrized by the family.
\end{definition}

Often our bounded families of sheaves will all be quotients of a fixed coherent sheaf.

\begin{definition}
Suppose $X$ is a smooth projective variety and $\mathcal{E}$ is a torsion-free sheaf on $X$.  A finite-type family of torsion-free quotients of $\mathcal{E}$ consists of a finite-type family of torsion-free sheaves $\mathcal{G}'$ on $X \times S$ equipped with a surjection $\pi_{1}^{*}\mathcal{E} \to \mathcal{G}'$ where $\pi_1 : X \times S \to X$ is the projection.

We say that a set of torsion-free quotients $\{ \mathcal{E} \to \mathcal{Q}_{i} \}$ on $X$ is a bounded family if there is a finite-type family of torsion-free quotients such that every $\mathcal{E} \to \mathcal{Q}_{i}$ is parametrized by the family.
\end{definition}

Using standard properties of the Quot scheme, we have:

\begin{lemma} \label{lemm:familytoquotient}
Let $X$ be a smooth projective variety and let $\mathcal{E}$ be a torsion-free coherent sheaf on $X$.  Suppose $\mathcal{G}$ is a bounded family of torsion-free sheaves on $X$.  Then the set of all quotients $\mathcal{E} \to \mathcal{Q}$ such that $\mathcal{Q}$ is parametrized by $\mathcal{G}$ is a bounded family of torsion-free quotients of $\mathcal{E}$.
\end{lemma}

We also use the following result which is a variant of a theorem of Grothendieck.

\begin{theorem}[{\cite[Lemma 1.7.9]{HL97}}, {\cite[Theorem 2.30]{GKP16}}] \label{theo:grothendieckbounded}
Let $X$ be a smooth projective variety and fix a compact set $\mathcal{T}$ contained in the interior of $\Nef_{1}(X)$.  Let $\mathcal{E}$ be a torsion-free sheaf on $X$.  Fix a constant $R$ and consider the set of torsion-free quotients $\mathcal{E} \to \mathcal{Q}$ such that $\mu_{\alpha}(\mathcal{Q}) < R$ for some $\alpha \in \mathcal{T}$.  This is a bounded family of torsion-free quotients of $\mathcal{E}$.
\end{theorem}

\begin{proof}
Since there is a surjection $\mathcal{O}_{X}(L)^{\oplus d} \to \mathcal{E}$ for some divisor $L$ and some positive integer $d$, it suffices to prove the statement for $\mathcal{O}_{X}(L)^{\oplus d}$.  In fact, after twisting we may assume that our initial sheaf is $\mathcal{O}_{X}^{\oplus d}$.  One can then repeat the argument of \cite[Theorem 2.30]{GKP16} with essentially no change to show that $c_{1}(\mathcal{Q})$ lies in a bounded subset of $N^1(X)_{\mathbb{Z}}$.  The boundedness of the set of quotients then follows from \cite[Lemma 1.7.9]{HL97} and Lemma \ref{lemm:familytoquotient}.
\end{proof}

\subsection{Log canonical thresholds}

In this subsection we give a brief reminder of the definition and basic properties of the log canonical threshold of an ideal.

\begin{definition}
Let $X$ be a smooth projective variety and let $\mathcal{I}$ be a coherent ideal sheaf on $X$.  Let $\phi: X' \to X$ denote a log resolution of $\mathcal{I}$ and set $\mathcal O(-F) = \phi^{-1}\mathcal{I} \cdot \mathcal{O}_{X'}$.  The log canonical threshold $\lct(\mathcal{I})$ is defined to be the infimum over all real numbers $c$ such that
\begin{equation*}
\phi_{*}\mathcal{O}_{X'}(K_{X'/X} - \lfloor cF \rfloor) \subsetneq \mathcal{O}_{X}.
\end{equation*}
If the subscheme defined by $\mathcal{I}$ has codimension at least $2$ in $X$, then $\lct(\mathcal{I})$ is equivalently the maximum over all real numbers $d$ such that $K_{X'/X} + \mathrm{Exc}_{\phi} \geq dF$ where $\mathrm{Exc}_{\phi}$ denotes the sum of all reduced $\phi$-exceptional divisors (see \cite[Example 9.3.16]{Lazarsfeld04b}).
\end{definition}

Note that if we have an inclusion of ideal sheaves $\mathcal{I} \subset \mathcal{J}$ then we also have an inequality $\lct(\mathcal{I}) \leq \lct(\mathcal{J})$.  We will use this fact repeatedly in what follows.

The next theorem gives a useful bound on the value of the log canonical threshold.  This result is certainly well-known to experts; for example, some version of this statement appears in \cite[Proof of Theorem 9.14]{BHJ17}.

\begin{theorem} \label{theo:lctestimate}
Let $X$ be a smooth projective variety of dimension $n$.  Fix a curve class $\beta \in N_{1}(X)_{\mathbb{Z}}$ which lies in the interior of $\Nef_{1}(X)$.  There is a positive constant $\upsilon$ (which depends on $X$ and $\beta$) with the following property.

Let $\mathcal{I}$ be an ideal sheaf on $X$.  Let $D$ be an effective Cartier divisor such that $\mathcal{O}_{X}(-D) \subset \mathcal{O}_{X}$ is contained in $\mathcal{I}$.  Then
\begin{equation*}
\lct(\mathcal{I}) \geq \frac{1}{\upsilon (D \cdot \beta) }.
\end{equation*}
\end{theorem}

Since we could not find a version of this statement that could be clearly cited we will give a proof here, following the strategy of \cite{BHJ17}.

\begin{proof}
We first define the constant $\upsilon$.  By Lemma \ref{lemm:nakayama} there is an ample Cartier divisor $H$ with the following property: for every positive integer $d$, we have $H^{0}(X,\mathcal{O}_{X}(dH - B)) > 0$ for every effective Cartier divisor $B$ satisfying $B \cdot \beta \leq d$. We claim that there is a constant $\upsilon$ such that for every positive integer $d$ and every effective divisor $L \sim dH$ we have
\begin{equation*}
\sup_{x \in X} \ord_{x}(L) \leq \upsilon d.
\end{equation*}
Indeed, suppose that $\phi: X' \to X$ is the blow-up of a point $x \in X$ and let $E$ denote the exceptional divisor.  By Seshadri's criterion (see \cite[Theorem 1.4.13]{Lazarsfeld04a}) 
there is some $\epsilon > 0$ (depending on $X,H$ but not on $x$) such that $\phi^{*}H - \epsilon E$ is ample.  Then
\begin{align*}
L^{n} = \phi^{*}L \cdot (\phi^{*}L - \epsilon d E)^{n-1} & \geq  \ord_{x}(L) E \cdot (\phi^{*}L - \epsilon d E)^{n-1} \\
& = d^{n-1} \epsilon^{n-1} \ord_{x}(L) 
\end{align*}
Thus $\upsilon = \frac{H^{n}}{\epsilon^{n-1}}$ satisfies the claim above.

We next verify that $\upsilon$ has the desired property.  First recall that if we have an inclusion of ideals $\mathcal{J} \subset \mathcal{I}$ then $\lct(\mathcal{J}) \leq \lct(\mathcal{I})$.  In particular, $\lct(D) \leq \lct(\mathcal{I})$.  

Set $r = D \cdot \beta$.  The defining property of $H$ implies that there is some effective Cartier divisor $F$ such that $F + D$ is linearly equivalent to $rH$.  Note we have $\lct(D) \geq \lct(F + D)$.  In turn, \cite[Lemma 8.10]{Kollar97} shows that
\begin{align*}
\lct(F+D) & \geq \inf_{x \in X} \frac{1}{\ord_{x}(F+D)} \\
& \geq \frac{1}{\upsilon r}
\end{align*}
We conclude that $\lct(\mathcal{I}) \geq \frac{1}{\upsilon r}$ as desired.
\end{proof}

The following well-known lemma was used in the previous proof.

\begin{lemma} \label{lemm:nakayama}
Let $X$ be a smooth projective variety and let $\mathcal{C} \subset \Eff^{1}(X)$ denote a compact subset.  There is a ample divisor $H$ with the following property: for any positive integer $d$ and any Cartier divisor $B$ whose numerical class lies in $d\mathcal{C} \cap N^{1}(X)_{\mathbb{Z}}$ we have $H^{0}(X,\mathcal{O}_{X}(dH-B)) > 0$.
\end{lemma}

\begin{proof}
Since $\mathcal{C}$ is bounded, there is an effective ample divisor $A$ such that $A - \mathcal{C}$ is contained in the ample cone of $X$.  By \cite[Corollary V.1.4]{Nakayama04} there is an ample divisor $L$ such that $H^{0}(X,\mathcal{O}_{X}(L+D)) > 0$ for every pseudo-effective Cartier divisor $D$.

Set $H = A+L$.  Then for any divisor $B$ as in the statement, we have
\begin{align*}
dH - B = L + (d-1)L + (dA - B).
\end{align*}
Since $(d-1)L + (dA - B)$ is pseudo-effective, $dH - B$ has a non-vanishing section.
\end{proof}

We will use the log canonical threshold to control the anticanonical divisor.

\begin{lemma} \label{lemm:lctvskx}
Suppose that $\phi: X' \to X$ is a birational morphism of smooth projective varieties that resolves an an ideal sheaf $\mathcal{I}$ on $X$ which defines a closed subscheme of codimension $\geq 2$.  Write $\phi^{-1}\mathcal{I} \cdot \mathcal{O}_{X'} = \mathcal{O}_{X'}(-D)$.

Let $\alpha'$ be any nef curve class on $X'$.  Then
\begin{equation*}
K_{X'/X} \cdot \alpha' \geq \frac{1}{2} \lct(\mathcal{I}) (D \cdot \alpha').
\end{equation*}
\end{lemma}

\begin{proof}
By definition of the log canonical threshold, we have
\begin{equation*}
K_{X'/X} + \mathrm{Exc}_{\phi} \geq \lct(\mathcal{I}) \cdot D
\end{equation*}
where $\mathrm{Exc}_{\phi}$ denotes the sum of the $\phi$-exceptional divisors with coefficients $1$.
Since $X'$ and $X$ are smooth, every exceptional divisor appears in $K_{X'/X}$ with coefficient at least $1$.  Thus 
\begin{equation*}
2K_{X'/X} \geq K_{X'/X} + \mathrm{Exc}_{\phi}.
\end{equation*}
Combining we see that $2K_{X'/X} - \lct(\mathcal{I}) \cdot D$ is effective.  Thus the claim follows from the fact that $\alpha'$ is nef.
\end{proof}

\section{Birational geometry of Fitting ideals}

This section is devoted to a number of results and constructions involving the Fitting ideal of a coherent sheaf on a smooth projective variety.  

\subsection{Fitting ideals}

Let $X$ be a smooth projective variety and let $\mathcal{Q}$ denote a torsion-free sheaf on $X$. We would like to find a birational morphism $\phi: X' \to X$ such that $(\phi^{*}\mathcal{Q})_{tf}$ is locally free.  The construction of such a morphism $\phi$ is given by the theory of Fitting ideals, which we now recall.

First suppose we take an open affine $\Spec(A) \subset X$ and that $S$ is a coherent $A$-module with a presentation
\begin{equation*}
A^{d} \xrightarrow{M} A^{n} \to S \to 0
\end{equation*}
Then the $j$th Fitting ideal $\Fit_{j}(S)$ on $\Spec(A)$ is generated by the $(n-j)$-minors of the matrix $M$.  It turns out that the resulting ideal sheaf $\Fit_{j}(S)$ is independent of the choice of presentation of $S$. This construction globalizes to give a Fitting ideal for a coherent sheaf $\mathcal{Q}$ on $X$.  Fitting ideals have the following important properties:
\begin{itemize}
\item The formation of the $j$th Fitting ideal commutes with base change.
\item $\Supp \Fit_{j}(\mathcal{Q}) = \{ x \in X | \dim_{\kappa(x)} ( \mathcal{Q} \otimes \kappa(x) ) > j\}$.
\item Since a torsion-free sheaf $\mathcal{Q}$ on a smooth variety is locally free in codimension $1$, for any $j \geq \rk(\mathcal{Q})$ the support of $\Fit_{j}(\mathcal{Q})$ defines a closed subscheme of codimension $\geq 2$.
\end{itemize}
The following result shows that the birational morphism $\phi$ which makes a torsion-free sheaf $\mathcal{Q}$ locally free is the same as the resolution of a Fitting ideal.

\begin{theorem}[{\cite[Chapter 4, \S 3, Lemma 1]{Raynaud72}}]
\label{theo:resolvingfittingideal}
Let $X$ be a smooth projective variety and let $\mathcal{Q}$ be a coherent sheaf on $X$ of rank $r$.  
Let $\phi: X' \to X$ denote a smooth birational model resolving $\Fit_{r}(\mathcal{Q})$.  Then $(\phi^{*}\mathcal{Q})_{tf}$ is locally free of rank $r$.
\end{theorem}

The following lemma identifies the geometric meaning of the divisor resolving the Fitting ideal.

\begin{lemma} \label{lemm:fittidealcomputation}
Let $X$ be a smooth projective variety and let $\mathcal{Q}$ be a coherent torsion-free sheaf on $X$ of rank $r$.  Suppose that $\phi: X' \to X$ is a birational morphism from a smooth projective variety that resolves the $r$th Fitting ideal $\Fit_{r}(\mathcal{Q})$.  Write $\phi^{-1}\Fit_{r}(\mathcal{Q}) \cdot \mathcal{O}_{X'} = \mathcal{O}_{X'}(-D)$ for some effective $\phi$-exceptional divisor $D$.
Then
\begin{equation*}
D = \phi^{*}c_{1}(\mathcal{Q}) - c_{1}((\phi^{*}\mathcal{Q})_{tf}).
\end{equation*}
\end{lemma}

\begin{proof}
Consider the exact sequence
\begin{equation*}
0 \to (\phi^{*}\mathcal{Q})_{tors} \to \phi^{*}\mathcal{Q} \to (\phi^{*}\mathcal{Q})_{tf} \to 0.
\end{equation*}
Note that the rightmost term is locally free of rank $r$ by Theorem \ref{theo:resolvingfittingideal}.  In particular $\Fit_{j}((\phi^{*}\mathcal{Q})_{tf})$ is $0$ for $j < r$ and is $\mathcal{O}_{X'}$ for $j \geq r$.

Choose an open cover of $X'$ by open affines $\Spec(A) \subset X'$ so that the restriction of $(\phi^{*}\mathcal{Q})_{tf}$ to $\Spec(A)$ is a projective $A$-module.  In particular the restriction of the exact sequence above to $\Spec(A)$ splits.  Applying \cite[Tag 07ZA]{stacks}, we see that on $\Spec(A)$ we have
\begin{align*}
\Fit_{r}(\phi^{*}\mathcal{Q}) & = \sum_{i+j = r} \Fit_{i}((\phi^{*}\mathcal{Q})_{tors}) \Fit_{j}((\phi^{*}\mathcal{Q})_{tf}) \\
& = \Fit_{0}((\phi^{*}\mathcal{Q})_{tors}).
\end{align*}
Of course this local equality for every $\Spec(A)$ implies that we also have a global equality.  
Since Fitting ideals are compatible with base-change, we conclude that $\mathcal{O}_{X'}(-D) \cong \Fit_{0}((\phi^{*}\mathcal{Q})_{tors})$. 

It only remains to show that the first Chern class of $\Fit_{0}((\phi^{*}\mathcal{Q})_{tors})$ can be identified with $c_{1}((\phi^{*}\mathcal{Q})_{tf}) - \phi^{*}c_{1}(\mathcal{Q})$.  Equivalently, we must show that the first Chern class of $\Fit_{0}((\phi^{*}\mathcal{Q})_{tors})$ is the negative of the first Chern class of $(\phi^{*}\mathcal{Q})_{tors}$.
By choosing a surjection to $(\phi^{*}\mathcal{Q})_{tors}$ from a locally free sheaf we obtain an exact sequence
\begin{equation*}
0 \to \mathcal{E}_{2} \xrightarrow{\psi} \mathcal{E}_{1} \to (\phi^{*}\mathcal{Q})_{tors} \to 0
\end{equation*}
where $\mathcal{E}_{1}$ is locally free and $\mathcal{E}_{2}$ is torsion-free of the same rank.  We let $t$ denote the common rank of $\mathcal{E}_{1}$ and $\mathcal{E}_{2}$.  This exact sequence shows that $c_{1}((\phi^{*}\mathcal{Q})_{tors})$ is equal to $c_{1}(\mathcal{E}_{1}) - c_{1}(\mathcal{E}_{2})$.

On the other hand, let $U \subset X'$ denote the open locus where $\mathcal{E}_{2}$ is locally free.  Note that the complement of $U$ has codimension $\geq 2$.  By taking top exterior powers in the exact sequence above we see that $\Fit_{0}((\phi^{*}\mathcal{Q})_{tors})|_{U}$ is the ideal sheaf corresponding to the effective divisor $F_{U}$ defining the map of line bundles
\begin{equation*}
\bigwedge^{t} \psi|_{U}: \bigwedge^{t}\mathcal{E}_{2}|_{U} \xrightarrow{\cdot F_{U}} \bigwedge^{t} \mathcal{E}_{1}|_{U}
\end{equation*}
We conclude that $F_{U}$ is equal to $c_{1}(\mathcal{E}_{1}|_{U}) - c_{1}(\mathcal{E}_{2}|_{U})$.  Since the complement of $U$ has codimension $\geq 2$ we can take closures in $X'$ to obtain the equality $\Fit_{0}((\phi^{*}\mathcal{Q})_{tors}) = -(c_{1}(\mathcal{E}_{1}) - c_{1}(\mathcal{E}_{2}))$, finishing the proof.
\end{proof}

Our next result will rely on the following lemma:

\begin{lemma}
\label{lemma:quotient_effective}
Let $X$ be a smooth projective variety and let $\mathcal{Q}$ denote a non-zero torsion-free coherent sheaf on $X$ of rank $r$.  Suppose that $L$ is a Cartier divisor and $d$ is a positive integer such that there is a surjection $\mathcal{O}_{X}(L)^{\oplus d} \to \mathcal{Q}$.  Then $c_{1}(\mathcal{Q}) - rL$ is (linearly equivalent to) an effective divisor.
\end{lemma}

\begin{proof}
Taking $r$th exterior powers and a double dual, we obtain a non-zero morphism $\mathcal{O}_{X}(rL) \to \mathcal{O}_{X}(c_{1}(\mathcal{Q}))$.  Thus the difference between these divisors is effective.
\end{proof}

Our main result in this subsection constructs functions in the Fitting ideal of bounded degree.  We will use it to estimate the log canonical threshold of Fitting ideals.

\begin{theorem} \label{theo:fittingidealgenerators}
Let $X$ be a smooth projective variety and let $\mathcal{Q}$ be a torsion-free sheaf on $X$ of positive rank $r$.  Suppose that $L$ is a Cartier divisor and $d$ is a positive integer such that we have a surjection $\mathcal{O}_{X}(L)^{\oplus d} \to \mathcal{Q}$.  Then $\Fit_{r}(\mathcal{Q})$ contains $\mathcal{O}_{X}(-D) \subset \mathcal{O}_{X}$ for an effective divisor $D$ satisfying
\begin{equation*}
D \in \left| (d-r) r \cdot (c_{1}(\mathcal{Q}) - rL) \right|.
\end{equation*}
\end{theorem}

\begin{proof}
For notational convenience we write $\mathcal{F} := \mathcal{O}_{X}(L)^{\oplus d}$.  We denote the kernel of the surjection  $\mathcal{F} \to \mathcal{Q}$ by $\mathcal{K}$.
Consider the composition
\begin{equation*}
\vartheta: \mathcal{F} \otimes \bigwedge^{r-1} \mathcal{F} \to \bigwedge^{r} \mathcal{F} \to \left( \bigwedge^{r} \mathcal{Q} \right)^{\vee \vee}
\end{equation*}
First suppose we look at the open locus $U \subset X$ where $\mathcal{Q}$ is locally free.  Note that a local section $s$ of $\mathcal{F}$ on $U$ will be contained in the kernel $\mathcal{K}$ if and only for every element $t \in \bigwedge^{r-1} \mathcal{F}(U)$ we have that $\vartheta(s \otimes t) = 0$.  In other words, if we define the induced map
\begin{equation*}
\psi: \mathcal{F} \to \left(\bigwedge^{r-1} \mathcal{F}^{\vee} \right) \otimes \mathcal{O}_{X}(c_{1}(\mathcal{Q})).
\end{equation*}
then $\mathcal{K}|_{U}$ is the same as $(\ker \psi)|_{U}$.  Since both $\mathcal{Q}$ and the image of $\psi$ are torsion-free, both $\mathcal{K}$ and the kernel of $\psi$ are saturated subsheaves of $\mathcal{F}$.  Since these subsheaves agree on $U$, they are equal on all of $X$.
In other words, $\mathcal{Q}$ is isomorphic to the image of a map
\begin{equation*}
\mathcal{O}_{X}(L)^{\oplus d} \xrightarrow{N} \mathcal{O}_{X}(c_{1}(\mathcal{Q}) - (r-1)L)^{\oplus {d \choose r-1}}
\end{equation*}
defined by a matrix $N$ of forms in $H^{0}(X,\mathcal{O}_{X}(c_{1}(\mathcal{Q}) - rL))$.

We know that the restriction of the matrix $N$ to the generic point of $X$ has rank $r$.  Let $A$ be an $r \times r$ submatrix of $N$ such that the restriction to the generic point has full rank and let $G$ be the $r \times d$ submatrix of $N$ whose rows have the same indices as $A$.  
Letting $\mathrm{adj}(A)$ denote the adjugate matrix, if we take the product $\mathrm{adj}(A) \cdot G$ 
then there is an $r \times r$ submatrix equal to $\det(A) \cdot \mathrm{Id}_{r}$.  For each of the $d-r$ columns not in this submatrix, we obtain an element in the kernel of $\mathrm{adj}(A) \cdot G$ which has entries in $H^{0}(X, \mathcal{O}_{X}(\det (A)))$ by comparing it against the $r$ columns in the submatrix $\mathrm{adj}(A) \cdot G$.  Furthermore, these $d-r$ elements span the kernel of the restriction of $\mathrm{adj}(A) \cdot G$ to the generic point, which by construction is the same as the kernel of the restriction of $N$ to the generic point.

Since the line bundle $\mathcal{O}_{X}(\det(A))$ is the same as $\mathcal{O}_{X}(r(c_{1}(\mathcal{Q}) - rL))$, altogether the construction above yields an injective morphism 
\begin{equation*}
\phi: \mathcal{O}_{X}\left(L - r \cdot (c_{1}(\mathcal{Q}) - rL) \right)^{\oplus d-r} \to \mathcal{O}(L)^{\oplus d}
\end{equation*}
whose image is a full-rank subsheaf of the kernel $\mathcal{K}$.  Note that $\phi$ is defined by a matrix $B$ whose entries lie in $H^{0}\left(X,\mathcal{O}_{X}\left( r \cdot (c_{1}(\mathcal{Q}) - rL) \right)\right)$.

Next, choose some other locally free sheaf $\mathcal{T}$ equipped with a morphism $\rho: \mathcal{T} \to \mathcal{O}(L)^{\oplus d}$ such that the image of
\begin{equation*}
\mathcal{T} \oplus \mathcal{O}_{X}\left(L - r \cdot (c_{1}(\mathcal{Q}) - rL) \right)^{\oplus d-r} \xrightarrow{\rho \oplus \phi} \mathcal{O}(L)^{\oplus d}
\end{equation*}
is exactly equal to $\mathcal{K}$.  We can locally compute $\Fit_{r}(\mathcal{Q})$ by taking $(d-r)$-minors of the matrix defining $\rho \oplus \phi$.
In particular, the $(d-r)$-minors of the matrix $B$ defining $\phi$ will be contained in $\Fit_{r}(\mathcal{Q})$.  In this way we see that $\Fit_{r}(\mathcal{Q})$ contains $\mathcal{O}_{X}(-D)$ where $D \in H^{0}\left(X,\mathcal{O}_{X}\left( (d-r) r \cdot (c_{1}(\mathcal{Q}) - rL) \right)\right)$.
\end{proof}

\begin{corollary} \label{coro:lctforfitting}
Let $X$ be a smooth projective variety.  Fix a curve class $\beta \in N_{1}(X)_{\mathbb{Z}}$ which lies in the interior of $\Nef_{1}(X)$.  There is a positive constant $\upsilon$ (which depends on $X$ and $\beta$) with the following property.

Let $\mathcal{Q}$ be a torsion-free sheaf on $X$ of positive rank $r$.  Suppose that $L$ is a Cartier divisor and $d$ is a positive integer such that we have a surjection $\mathcal{O}_{X}(L)^{\oplus d} \to \mathcal{Q}$.  Then
\begin{equation*}
\lct(\Fit_{r}(\mathcal{Q})) \geq \frac{1}{\upsilon (d-r) r  ((c_{1}(\mathcal{Q}) - rL) \cdot \beta) }.
\end{equation*}
\end{corollary}

\begin{proof}
Define $\upsilon$ as in Theorem \ref{theo:lctestimate}.  
Theorem \ref{theo:fittingidealgenerators} shows that $\Fit_{r}(\mathcal{Q})$ contains $\mathcal{O}_{X}(-D)$ where $D$ is a divisor satisfying
\begin{equation*}
D \in \left| (d-r) r \cdot (c_{1}(\mathcal{Q}) - rL) \right|.
\end{equation*}
Theorem \ref{theo:lctestimate} provides the desired bound on the log canonical threshold.
\end{proof}

\subsection{Resolving Fitting ideals in families}

Finally, we will need a construction of a bounded family of birational maps which simultaneously resolves a family of Fitting ideals.

\begin{construction} \label{cons:familyofbirationalmodels}
Let $X$ be a smooth projective variety.  Suppose that $\mathcal{G}$ on $X \times S$ is a finite-type family of torsion-free sheaves on $X$ (in the sense of Definition \ref{defi:boundedfamily}). 

Recall that Fitting ideals are compatible with base change: for every positive integer $r$ the restriction of $\Fit_{r}(\mathcal{G})$ to a fiber $X_{s}$ is the same as $\Fit_{r}(\mathcal{G}|_{X_{s}})$.  In particular, since the rank of $\mathcal{G}|_{X_{s}}$ is constant in $S$, the restriction of the top Fitting ideal for $\mathcal{G}$ will be the top Fitting ideal for $\mathcal{G}|_{X_{s}}$.

We claim that there is a finite stratification of $S$ into locally closed subvarieties $S = \sqcup_{k} S_{k}$, a collection of birational morphisms $\psi_{k}: \mathfrak{Y}_{k} \to X \times S_{k}$, and a locally free sheaf $\mathfrak{G}_{k}$ on $\mathfrak{Y}_{k}$ satisfying the following properties for every $k$: 
\begin{enumerate}
\item the induced morphism $\mathfrak{Y}_{k} \to S_{k}$ is smooth projective with connected fibers and the restriction of $\psi_{k}$ to the fiber over any closed point of $S_{k}$ is birational,
\item $\mathfrak{G}_{k}$ is isomorphic to $(\psi_{k}^{*}\mathcal{G}|_{X \times S_{k}})_{tf}$.
\end{enumerate}
We construct this stratification of $S$ by Noetherian induction.  We first construct a smooth birational model $Y \to X \times S$ that resolves the top Fitting ideal of $\mathcal{G}$.  The morphism $Y \to S$ is generically smooth; we let $S_{0} \subset S$ denote the open locus over which the map is smooth and let $\mathfrak{Y}_{0}$ be the preimage of $S_{0}$.  It is clear that $\psi_{0}$ resolves the top Fitting ideal of $\mathcal{G}|_{X \times S_{0}}$.  It is also clear that for each point $s \in S_{0}$ the map $\mathfrak{Y}_{0,s} \to X_{s}$ is birational and resolves the top Fitting ideal of $\mathcal{G}|_{X_{s}}$.  Thus $\psi_{0}: \mathfrak{Y}_{0} \to X \times S_{0}$ satisfies both desired properties.  We then repeat the construction on the irreducible components of $X \times (S \backslash S_{0})$.

\end{construction}

\section{Sequences of quotients}

Our main theorems involve sequences of quotients of sheaves (often arising from Harder-Narasimhan filtrations).  In this section we describe some basic constructions involving sequences of quotients.

Suppose we have a sequence $(\mathcal{Q}_{\bullet})$ of torsion-free quotients of a torsion-free coherent sheaf $\mathcal{E}$:
\begin{equation*}
\mathcal{E} = \mathcal{Q}_{s} \to \mathcal{Q}_{s-1} \to \ldots \to \mathcal{Q}_{1} \to \mathcal{Q}_{0} = 0.
\end{equation*}
This contains the same information as the sequence of saturated subsheaves $(\mathcal{G}_{\bullet})$ of $\mathcal{E}$ defined by $\mathcal{G}_{s-i} = \ker(\mathcal{E} \to \mathcal{Q}_{i})$.  In this setting we define $\SP_{X,\alpha}(\mathcal{E};\mathcal{Q}_{\bullet}) = \SP_{X,\alpha}(\mathcal{E};\mathcal{G}_{\bullet})$.  Henceforth we will work exclusively with sequences of quotients and not filtrations.

\subsection{Birational transforms}

We first discuss how to push sequences of quotients back and forth over birational morphisms of smooth projective varieties.

\begin{definition}
Let $\phi: X' \to X$ be a birational morphism of smooth projective varieties.  Suppose that $\mathcal{E}$ is a torsion-free coherent sheaf on $X$ and we have a sequence of torsion-free quotients
\begin{equation*}
\mathcal{E} = \mathcal{Q}_{s} \to \mathcal{Q}_{s-1} \to \ldots \to \mathcal{Q}_{1} \to \mathcal{Q}_{0} = 0.
\end{equation*}
The birational transform of the sequence $\mathcal{Q}_{\bullet}$ is the sequence of torsion-free quotients of $(\phi^{*}\mathcal{E})_{tf}$ defined by
\begin{equation*}
\mathcal{Q}'_{i} := (\phi^{*}\mathcal{Q}_{i})_{tf}.
\end{equation*}
\end{definition}

\begin{remark}
When we discuss the birational transform of a Harder-Narasimhan filtration of $\mathcal{E}$, this should be interpreted as applying the construction to the corresponding sequence of torsion-free \emph{quotients} of $\mathcal{E}$.  Note that the torsion-free parts of the pullbacks of subsheaves defining the filtration may not filter the torsion-free part of the pullback of $\mathcal{E}$.
\end{remark}

\subsection{Birational pushforwards}

Suppose we have a birational morphism $\phi: X' \to X$ and a sequence of torsion-free quotients of a sheaf on $X'$.  We would like to turn this into a sequence of sheaves on $X$.  The basic construction is the following.

\begin{lemma} \label{lemm:birationalconstruction}
Let $X$ be a smooth projective variety and let $\mathcal{E}$ be a torsion-free sheaf on $X$.  Suppose that $\phi: X' \to X$ is a birational morphism from a smooth projective variety $X'$.  For every torsion-free quotient $(\phi^{*}\mathcal{E})_{tf} \to \mathcal{Q}'$, there is a unique torsion-free quotient $\mathcal{E} \to \mathcal{Q}$ such that there is an isomorphism of quotients
\begin{equation*}
\xymatrix{
(\phi^{*}\mathcal{E})_{tf} \ar[r] \ar[dr] & (\phi^{*}\mathcal{Q})_{tf} \ar[d]^{\cong} \\
& \mathcal{Q}'
}
\end{equation*}
\end{lemma}

\begin{proof}
We define $\mathcal{Q}$ to be the image of the composed morphism
\begin{equation*}
\mathcal{E} \to \phi_{*}\phi^{*}\mathcal{E} \to \phi_{*}\mathcal{Q}'.
\end{equation*}
Since torsion-freeness is preserved by pushforward, $\phi_{*}\mathcal{Q}'$ and thus also $\mathcal{Q}$ are torsion-free.  We define $\mathcal{F}$ to be the kernel of $\mathcal{E} \to \mathcal{Q}$ and define $\mathcal{F}'$ to be the kernel of $(\phi^{*}\mathcal{E})_{tf} \to \mathcal{Q}'$.

Note that the image of $\phi^{*}\mathcal{F} \to (\phi^{*}\mathcal{E})_{tf}$ agrees with $\mathcal{F}'$ when restricted to the open locus where $\phi$ is an isomorphism.  Since $\mathcal{F}'$ is saturated, we obtain an induced morphism $\phi^{*}\mathcal{F} \to \mathcal{F}'$ whose cokernel is torsion.  Applying the Snake Lemma to the diagram
\begin{equation*}
\xymatrix{
& \phi^{*}\mathcal{F} \ar[r] \ar[d] & \phi^{*}\mathcal{E} \ar[r] \ar[d] & \phi^{*}\mathcal{Q} \ar[r] \ar@{.>}[d] & 0 \\
0 \ar[r] & \mathcal{F}' \ar[r] & (\phi^{*}\mathcal{E})_{tf} \ar[r] & \mathcal{Q}' \ar[r] & 0
}
\end{equation*}
shows that the induced morphism $\phi^{*}\mathcal{Q} \to \mathcal{Q}'$ has torsion kernel and trivial cokernel.  Thus $(\phi^{*}\mathcal{Q})_{tf} \to \mathcal{Q}'$ is an isomorphism.

It only remains to show the uniqueness of $\mathcal{Q}$.  Suppose that $\mathcal{Q}_{1}, \mathcal{Q}_{2}$ are two torsion-free quotients with this property.  For $i=1,2$ let $\mathcal{F}_{i}$ denote the kernel of $\mathcal{E} \to \mathcal{Q}_{i}$.  On the open set $U$ where $\phi$ is an isomorphism, we have $\mathcal{F}_{1}|_{U} = \mathcal{F}_{2}|_{U}$ as subsheaves of $\mathcal{E}|_{U}$.  There is a unique saturated subsheaf of $\mathcal{E}$ whose restriction to $U$ agrees with this subsheaf.  Since $\mathcal{Q}_{1},\mathcal{Q}_{2}$ are torsion-free we conclude that $\mathcal{F}_{1} = \mathcal{F}_{2}$ as subsheaves of $\mathcal{E}$.
\end{proof}

\begin{definition}
Let $\phi: X' \to X$ be a birational morphism of smooth projective varieties.  Suppose $\mathcal{E}$ is a torsion-free sheaf on $X$ and we are given a sequence
\begin{equation*}
(\phi^{*}\mathcal{E})_{tf} = \mathcal{Q}'_{s} \to \mathcal{Q}'_{s-1} \to \ldots \to \mathcal{Q}'_{1} \to \mathcal{Q}'_{0} = 0
\end{equation*}
of torsion-free quotients of $(\phi^{*}\mathcal{E})_{tf}$.  The birational pushforward of $(\mathcal{Q}'_{\bullet})$ to $X$ is given by repeatedly applying Lemma \ref{lemm:birationalconstruction} to construct successive quotients of $\mathcal{E}$.  Precisely, we apply Lemma \ref{lemm:birationalconstruction} to $\mathcal{Q}'_{s} \to \mathcal{Q}'_{s-1}$ to construct $\mathcal{E} \to \mathcal{Q}_{s-1}$, then apply Lemma \ref{lemm:birationalconstruction} to $\mathcal{Q}'_{s-1} \to \mathcal{Q}'_{s-2}$ to construct $\mathcal{Q}_{s-1} \to \mathcal{Q}_{s-2}$, and so on.
\end{definition}

\begin{lemma}
Let $\phi: X' \to X$ be a birational morphism of smooth projective varieties and let $\mathcal{E}$ be a torsion-free sheaf on $X$.
\begin{enumerate}
\item Suppose $(\mathcal{Q}_{\bullet})$ is a sequence of torsion-free quotients of $\mathcal{E}$.  Then the birational pushforward of the birational transform of this sequence $(\mathcal{Q}_{\bullet})$ is the same as the original sequence.
\item Suppose $(\mathcal{Q}'_{\bullet})$ is a sequence of torsion-free quotients of $(\phi^{*}\mathcal{E})_{tf}$.  Then the birational transform of the birational pushforward of this sequence is the same as the original sequence.
\end{enumerate}
\end{lemma}

\begin{proof}
(1) Follows from the uniqueness of the quotient $\mathcal{Q}$ constructed in Lemma \ref{lemm:birationalconstruction}.

(2) Follows immediately from the definitions.
\end{proof}

\subsection{Graded pieces}

When working with sequences of quotients $(\mathcal{Q}_{\bullet})$, often we care not just about the quotients but also about the graded pieces $\mathcal{T}_{i} := \ker(\mathcal{Q}_{i} \to \mathcal{Q}_{i-1})$.  The following lemma indicates how to ``improve'' the graded pieces of a sequence of quotient sheaves by taking birational transforms.

\begin{lemma} \label{lemm:resolvingsingsofsequence}
Let $X$ be a smooth projective variety and let $\mathcal{E}$ be a torsion-free sheaf on $X$.  Let $(\mathcal{Q}_{\bullet})$ be a sequence of torsion-free quotients of $\mathcal{E}$.  There is a birational morphism $\phi: X' \to X$ from a smooth projective variety $X'$ such that each of the birational transforms $\mathcal{Q}_{i}' := (\phi^{*}\mathcal{Q}_{i})_{tf}$
and each of the graded pieces $\mathcal{T}_{i}' := \ker(\mathcal{Q}_{i}' \to \mathcal{Q}_{i-1}')$ is locally free.
\end{lemma}

\begin{proof}
By inductively resolving Fitting ideals, we can first find a birational morphism $\psi: \widetilde{X} \to X$ such that each birational transform $\widetilde{\mathcal{Q}}_{i} := (\psi^{*}\mathcal{Q}_{i})_{tf}$ is locally free.  Let $\widetilde{\mathcal{T}}_{i}$ denote the kernel of $\widetilde{\mathcal{Q}}_{i} \to \widetilde{\mathcal{Q}}_{i-1}$.  Now suppose that $\rho: X' \to \widetilde{X}$ is any smooth birational model and consider the exact sequence
\begin{equation*}
\xymatrix{
 \rho^{*}\widetilde{\mathcal{T}}_{i} \ar[r]  & \rho^{*}\widetilde{\mathcal{Q}}_{i} \ar[r] & \rho^{*}\widetilde{\mathcal{Q}}_{i-1} \ar[r] & 0
}
\end{equation*}
The leftmost arrow factors through $( \rho^{*}\widetilde{\mathcal{T}}_{i} )_{tf}$.  Since the induced morphism $( \rho^{*}\widetilde{\mathcal{T}}_{i} )_{tf} \to \rho^{*}\widetilde{\mathcal{Q}}_{i}$ is a generically injective morphism of torsion-free sheaves, it is injective.  In other words, the kernel of $\rho^{*}\widetilde{\mathcal{Q}}_{i} \to \rho^{*}\widetilde{\mathcal{Q}}_{i-1}$ is exactly the birational transform of $\widetilde{\mathcal{T}}_{i}$.  Thus we can finish the construction by choosing a birational morphism $\rho: X' \to \widetilde{X}$ from a smooth projective variety $X'$ that resolves the appropriate Fitting ideals of the sheaves $\widetilde{\mathcal{T}}_{i}$.
\end{proof}

\subsection{Families of quotient sequences}
Our next topic is families of quotient sequences, extending our earlier discussion of families of torsion-free sheaves.

\begin{definition} \label{defi:boundedfamilyquotientsequences}
Let $X$ be a smooth projective variety and let $\mathcal{E}$ be a torsion-free coherent sheaf on $X$.  A finite-type family of sequences of torsion-free quotients of $\mathcal{E}$ consists of:
\begin{itemize}
\item a finite-type separated scheme $S$ over our ground field, and
\item for every connected component $S^{j}$ of $S$ a sequence of quotients $(\mathcal{Q}_{i}^{j})$ of the sheaf $\pi_{1}^{*}\mathcal{E}$ on the product $X \times S^{j}$, where $\pi_{1}: X \times S^{j} \to X$ is the projection map
\end{itemize}
such that
\begin{enumerate}
\item for every connected component $S^{j}$ of $S$ and every index $i$ the sheaf $\mathcal{Q}_{i}^{j}$ is flat over $S^{j}$, and
\item for every connected component $S^{j}$, for every point $s \in S^{j}$, and for every index $i$ the restriction $\mathcal{Q}_{i}^{j}|_{X_{s}}$ is torsion-free.
\end{enumerate}

We say that a set of sequences of torsion-free quotients is bounded if there is a finite-type family of sequences of torsion-free quotients that parametrizes every sequence in our set.
\end{definition}

We will need a boundedness criterion for families of sequences of quotients.

\begin{lemma} \label{lemm:boundedfrombounded}
Let $X$ be a smooth projective variety and let $\mathcal{E}$ be a torsion-free coherent sheaf on $X$.  Suppose $\mathcal{G}$ on $X \times S$ is a finite-type family of torsion-free quotients of $\mathcal{E}$. 
Then the set of sequences of torsion-free quotients $(\mathcal{Q}_{\bullet})$ of $\mathcal{E}$ such that each $\mathcal{Q}_{i}$ is parametrized by $\mathcal{G}$ is also a bounded family.
\end{lemma}

\begin{proof}
We have an inclusion of relative Quot schemes $\mathrm{Quot}(\mathcal{G}/S) \to \mathrm{Quot}(\pi_{1}^{*}\mathcal{E}/S) \cong \mathrm{Quot}(\mathcal{E}) \times S$ where $\pi_{1}: X \times S \to X$ is the projection.  For every positive integer $1 \leq t \leq \rk(\mathcal{E})$, consider the product $X \times S^{\times t}$ equipped with the projection maps $p_{j}: X \times S^{\times t} \to X \times S$ onto the product of $X$ with the $j$th factor of $S^{\times t}$ and $r: X \times S^{\times t} \to X$.  We can construct the set of all sequences of quotients of total length $t$ as the closed subscheme $R$ of $S^{\times t}$ such that for every $1 \leq i < j \leq t$ the map $\mathrm{Quot}(p_{j}^{*}\mathcal{G}/S^{\times t})|_{R} \to \mathrm{Quot}(r^{*}\mathcal{E}/S^{\times t})|_{R}$ factors through the map $\mathrm{Quot}(p_{i}^{*}\mathcal{G}/S^{\times t})|_{R} \to \mathrm{Quot}(r^{*}\mathcal{E}/S^{\times t})|_{R}$.
\end{proof}

Finally, we discuss how to improve families of quotient sequences using birational morphisms.

\begin{construction} \label{cons:resolvingfamiliesofsequences}
Let $X$ be a smooth projective variety and let $\mathcal{E}$ be a torsion-free sheaf on $X$.  Suppose we have a finite-type family of quotient sequences of $\mathcal{E}$ consisting of
\begin{itemize}
\item a finite-type separated scheme $S$ and
\item for every connected component $S^{j}$ of $S$ a sequence of torsion-free quotients $(\mathcal{Q}_{i}^{j})$ of the sheaf $\pi_{1}^{*}\mathcal{E}$ on the product $X \times S^{j}$.
\end{itemize}
Then there exists a stratification of $S$ into locally closed subschemes $S = \sqcup_{k} S_{k}$ and a collection of smooth schemes $\mathfrak{Y}_{k}$ equipped with:
\begin{itemize}
\item surjective morphisms $\pi_{k}: \mathfrak{Y}_{k} \to S_{k}$ which have smooth projective connected fibers and
\item birational morphisms $\psi_{k}: \mathfrak{Y}_{k} \to X \times S_{k}$ whose restriction to every fiber over a point in $S_{k}$ is birational
\end{itemize}
that satisfies the following properties.
\begin{enumerate}
\item For every $S_{k} \subset S^{j}$ and for every term $\mathcal{Q}_{i}^{k} := \mathcal{Q}_{i}^{j}|_{X \times S_{k}}$ in the restriction of the quotient sequence, the birational transform sheaf $\mathfrak{Q}_{i}^{k, \prime}$ on $\mathfrak{Y}_{k}$ is locally free and
\item for all indices $i,k$ the restriction of $\ker(\mathfrak{Q}_{i}^{k,\prime} \to \mathfrak{Q}_{i-1}^{k,\prime})$ to the fibers of $\mathfrak{Y}_{k} \to S_{k}$ is a locally free sheaf.
\end{enumerate}
Indeed, first suppose we are given a single quotient sequence $(\mathcal{Q}_{\bullet})$ of $\mathcal{E}$.  Then Lemma \ref{lemm:resolvingsingsofsequence} shows that by repeatedly resolving Fitting ideals we can construct a single birational map $\phi: Y \to X$ such that both the birational transforms of the quotients and their successive kernels are locally free.  Since one can resolve Fitting ideals in families by Construction \ref{cons:familyofbirationalmodels}, a repeated application of Construction \ref{cons:familyofbirationalmodels} in the analogous way allows us to construct the desired $\mathfrak{Y}_{k}$.
\end{construction}

\section{Uniform bounds on birational stability} \label{sect:uniformbounds}

Suppose we have a birational morphism of smooth projective varieties $\phi: X' \to X$.  Let $\alpha'$ be a nef curve class on $X'$ and set $\alpha = \phi_{*}\alpha'$.  Given a locally free sheaf $\mathcal{E}$ on $X$, our goal in this section is to analyze the difference between the $\alpha'$-slope panel of $\phi^{*}\mathcal{E}$ and the $\alpha$-slope panel of $\mathcal{E}$.

\subsection{Birational behavior of stability}

There is one situation in which Harder-Narasimhan filtrations are compatible under birational transforms.

\begin{lemma}[{\cite[Proposition 2.8]{GKP16}}]
Let $\phi: X' \to X$ be a birational morphism of smooth projective varieties.  Let $\alpha \in \Nef_{1}(X)$ be a non-zero nef curve class. Suppose that $\mathcal{E}$ is a torsion-free sheaf on $X$ and let $(\mathcal{Q}_{\bullet})$ be the sequence of quotients defined by the $\alpha$-Harder-Narasimhan filtration of $\mathcal{E}$.

Then the $\phi^{*}\alpha$-Harder-Narasimhan filtration of $(\phi^{*}\mathcal{E})_{tf}$ is given by the birational transform of $(\mathcal{Q}_{\bullet})$.
\end{lemma}

Of course if $\alpha' \in \Nef_{1}(X')$ satisfies $\phi_{*}\alpha' = \alpha$ but $\alpha' \neq \phi^{*}\alpha$ then in general there is no reason to expect the $\alpha'$-Harder-Narasimhan filtration of $(\phi^{*}\mathcal{E})_{tf}$ to be directly related to the $\alpha$-Harder-Narasimhan filtration of $\mathcal{E}$.  However we have the following basic inequality.  

\begin{lemma} \label{lemm:furtherpullback}
Let $\psi: \widetilde{X} \to X$ be a birational morphism of smooth projective varieties.  Suppose that $\mathcal{E}$ is a torsion-free sheaf on $X$ and $\widetilde{\alpha} \in \Nef_{1}(\widetilde{X})$.  Define $\alpha = \psi_{*}\widetilde{\alpha}$ and $\widetilde{\mathcal{E}} = (\psi^{*}\mathcal{E})_{tf}$.  Then we have
\begin{equation*}
\mu_{\widetilde{\alpha}}(\widetilde{\mathcal{E}}) \leq \mu_{\alpha}(\mathcal{E}).
\end{equation*}
\end{lemma}

\begin{proof}
We know that $\psi^{*}c_{1}(\mathcal{E}) - c_{1}((\psi^{*}\mathcal{E})_{tf})$ is an effective divisor.  Thus
\begin{align*}
c_{1}((\psi^{*}\mathcal{E})_{tf}) \cdot \widetilde{\alpha} & \leq \psi^{*}c_{1}(\mathcal{E}) \cdot \widetilde{\alpha} \\
& = c_{1}(\mathcal{E}) \cdot \alpha
\end{align*}
leading to the desired inequality.
\end{proof}

In particular, this implies:

\begin{corollary} \label{coro:stillhnfilt}
Let $\psi: \widetilde{X} \to X$ be a birational morphism of smooth projective varieties.  Suppose that $\mathcal{E}$ is a locally free sheaf on $X$ and $\widetilde{\alpha} \in \Nef_{1}(\widetilde{X})$.  Define $\alpha = \psi_{*}\widetilde{\alpha}$.  If $\phi^{*}\mathcal{E}$ is $\widetilde{\alpha}$-semistable, then $\mathcal{E}$ is $\alpha$-semistable.
\end{corollary}

\begin{proof}
Let $\mathcal{E} \to \mathcal{Q}$ denote any torsion-free quotient.  Then by Lemma \ref{lemm:furtherpullback}
\begin{equation*}
\mu_{\alpha}(\mathcal{Q}) \geq \mu_{\widetilde{\alpha}}((\phi^{*}\mathcal{Q})_{tf}) \geq \mu_{\widetilde{\alpha}}(\phi^{*}\mathcal{E}) = \mu_{\alpha}(\mathcal{E}).
\end{equation*}
\end{proof}

\subsection{Bounds on birationally destabilizing sheaves}

Our first technical result is the following.  Loosely speaking, it shows that if the nef class $\alpha' \in \Nef_{1}(X')$ is ``more destabilizing'' for $\phi^{*}\mathcal{E}$ than $\phi_{*}\alpha'$ is for $\mathcal{E}$, we obtain a lower bound on the relative canonical degree of $\alpha'$.

\begin{theorem} \label{theo:birationaldegreebound3}
Let $X$ be a smooth projective variety and let $\mathcal{E}$ be a non-zero torsion-free sheaf of rank $n$ on $X$.  Fix a curve class $\beta \in N_{1}(X)_{\mathbb{Z}}$ in the interior of $\Nef_{1}(X)$.  Fix a positive integer $r$.  Fix a surjection $\mathcal{O}_{X}(L)^{\oplus d} \to \mathcal{E}$ where $L$ is a Cartier divisor and $d$ is a positive integer.  There is a bounded family $\mathcal{G}$ of torsion-free quotients of $\mathcal{E}$ of rank $r$ and a positive constant $\rho$ (which depend on $X, \mathcal{E}, \beta, r, L, d$) such that the following property holds.

Suppose we fix a non-zero curve class $\alpha \in \Nef_{1}(X)$. Then every surjection $\mathcal{E} \to \mathcal{Q}$ onto a torsion-free sheaf of rank $r$ satisfies one of the following properties:
\begin{enumerate}
\item $\mathcal{Q}$ is parametrized by the bounded family of torsion-free quotients $\mathcal G$, or
\item $\left(c_{1}(\mathcal{Q}) -  \frac{r}{n} c_{1}(\mathcal{E}) \right)   \cdot \beta > 0$ and for every birational morphism $\phi: X' \to X$ from a smooth projective $X'$, for every constant $R$, and for every $\alpha' \in \Nef_{1}(X')$ satisfying $\phi_{*}\alpha' = \alpha$ such that $\frac{1}{r} c_{1}((\phi^{*}\mathcal{Q})_{tf}) \cdot \alpha' \leq \frac{1}{n} c_{1}(\mathcal{E}) \cdot \alpha - \frac{1}{r}R$
we have
\begin{equation*}
K_{X'/X} \cdot \alpha' \geq  \rho \frac{  \left(c_{1}(\mathcal{Q}) -  \frac{r}{n} c_{1}(\mathcal{E}) \right) \cdot \alpha + R}{\left(c_{1}(\mathcal{Q}) -  \frac{r}{n} c_{1}(\mathcal{E}) \right)   \cdot \beta }
\end{equation*}
\end{enumerate}
\end{theorem}

\begin{proof}
Since $\beta$ lies in the interior of $\Nef_{1}(X)$, the set of torsion-free quotients $\mathcal{Q}$ of rank $r$ such that
\begin{equation*}
c_{1}(\mathcal{Q}) \cdot \beta \leq  \left( 2 \frac{r}{n}c_{1}(\mathcal{E}) - rL \right) \cdot \beta
\end{equation*}
is a bounded family by Theorem \ref{theo:grothendieckbounded}.
We let $\mathcal{G}$ denote this bounded family of quotients.  Also, we apply Theorem \ref{theo:lctestimate} to our fixed curve class $\beta$ to obtain a positive constant $\upsilon$.

Suppose $\mathcal{Q}$ is a torsion-free rank $r$ quotient of $\mathcal{E}$ not parametrized by $\mathcal{G}$.  Corollary \ref{coro:lctforfitting} shows that
\begin{equation*}
\lct(\Fit_{r}(\mathcal{Q})) \geq \frac{1}{\upsilon (d-r)r ( (c_{1}(\mathcal{Q}) - rL)  \cdot \beta)}
\end{equation*}
Suppose that $\phi: X' \to X$ is a birational morphism from a smooth projective variety, $\alpha' \in \Nef_{1}(X')$ pushes forward to $\alpha \in N_{1}(X)$, and $\mathcal E \to \mathcal Q$ is a rank $r$ quotient satisfying the intersection inequality $\frac{1}{r} c_{1}((\phi^{*}\mathcal{Q})_{tf}) \cdot \alpha' \leq \frac{1}{n} c_{1}(\mathcal{E}) \cdot \alpha - \frac{1}{r}R$. By Lemma \ref{lemm:furtherpullback} if we replace $X'$ by any higher birational model the slope of the torsion-free pullback of $\mathcal{Q}$ still satisfies this condition. Furthermore if we replace $\alpha'$ by its pullback to a higher birational model the intersection against the relative canonical divisor does not change.  Thus we may replace $X'$ by any higher birational model.  In particular, we may assume $\mathcal{Q}'$ is locally free, or equivalently, $\phi$ resolves the $r$th Fitting ideal of $\mathcal{Q}$.  Recall that a torsion-free sheaf is locally free in codimension $1$ and thus $\Fit_{r}(\mathcal{Q})$ defines a closed subscheme of $X$ of codimension $\geq 2$.  Letting $D$ denote the divisor defined by the inverse image ideal sheaf of $\Fit_{r}(\mathcal{Q})$, Lemma \ref{lemm:lctvskx} shows
\begin{align*}
K_{X'/X} \cdot \alpha' \geq \frac{1}{2}\lct(\Fit_{r}(\mathcal{Q})) (D \cdot \alpha') \geq \frac{D \cdot \alpha'}{2\upsilon (d-r) r ( (c_{1}(\mathcal{Q}) - rL)  \cdot \beta)}
\end{align*}
Lemma \ref{lemm:fittidealcomputation} shows that
\begin{align*}
D \cdot \alpha' & = (\phi^{*}c_{1}(\mathcal{Q}) - c_{1}(\mathcal{Q}') ) \cdot \alpha' \\
& \geq \left(c_{1}(\mathcal{Q}) -  \frac{r}{n} c_{1}(\mathcal{E}) \right) \cdot \alpha + R
\end{align*}
where the second line follows from our assumption on $\mathcal{Q}$.
Thus
\begin{align*}
K_{X'/X} \cdot \alpha' & \geq \frac{\left(c_{1}(\mathcal{Q}) -  \frac{r}{n} c_{1}(\mathcal{E}) \right) \cdot \alpha + R}{2\upsilon (d-r) r ( (c_{1}(\mathcal{Q}) - rL)  \cdot \beta)} 
\end{align*}

Note that the denominator is positive since $\beta$ is in the interior of the nef cone of curves and $c_{1}(\mathcal{Q}) - rL$ is pseudo-effective by Lemma~\ref{lemma:quotient_effective}.  By the inequality of intersection numbers we used to define $\mathcal{G}$, we know that 
\begin{equation*}
\left( c_{1}(\mathcal{Q}) - rL \right)   \cdot \beta \leq 2 \left(c_{1}(\mathcal{Q}) -  \frac{r}{n} c_{1}(\mathcal{E}) \right)   \cdot \beta
\end{equation*}
and in particular the right-hand side is positive.
Setting $\rho = \frac{1}{4\upsilon (d-r) r}$ we conclude that
\begin{align*}
K_{X'/X} \cdot \alpha' \geq \rho \frac{  \left(c_{1}(\mathcal{Q}) -  \frac{r}{n} c_{1}(\mathcal{E}) \right) \cdot \alpha + R}{\left(c_{1}(\mathcal{Q}) -  \frac{r}{n} c_{1}(\mathcal{E}) \right)   \cdot \beta }
\end{align*}
\end{proof}

In order to apply Theorem \ref{theo:birationaldegreebound3} in practice, we must control the relationship between the intersection number $\left(c_{1}(\mathcal{Q}) -  \frac{r}{n} c_{1}(\mathcal{E}) \right)  \cdot \alpha$ and the intersection number $\left(c_{1}(\mathcal{Q}) -  \frac{r}{n} c_{1}(\mathcal{E}) \right) \cdot \beta$.  In particular, if we assume that $\alpha$ is ``not too close'' to the boundary of $\Nef_{1}(X)$ then we can find a universal lower bound on the ratio between these two quantities.  The following theorem is the result of this computation.

\begin{theorem} \label{theo:semistabletwotypes}
Let $X$ be a smooth projective variety and let $\mathcal{E}$ be a non-zero torsion-free sheaf on $X$ of rank $n$.  Fix a big and nef divisor $H$ on $X$.  Fix a closed cone $\mathcal{C} \subset N_{1}(X)_{\mathbb{R}}$ such that $\mathcal{C} \backslash \{ 0 \}$ is contained in the interior of $\Nef_{1}(X)$ and $\mathcal{C}$ contains a class $\beta \in N_{1}(X)_{\mathbb{Z}}$ that is in the interior of $\Nef_{1}(X)$.  There is a bounded family $\mathcal{G}$ of torsion-free quotients of $\mathcal{E}$ and  a positive constant $\rho$ (which depend on $X,\mathcal{E},H,\mathcal{C},\beta$) with the following property.

Fix a non-zero curve class $\alpha \in \mathcal{C}$.  
Suppose that there is a birational model $\phi: X' \to X$ from a smooth projective variety $X'$, a class $\alpha' \in \Nef_{1}(X)_{\mathbb{Z}}$ satisfying $\phi_{*}\alpha' = \alpha$, and an $\alpha'$-destabilizing quotient $(\phi^{*}\mathcal{E})_{tf} \to \mathcal{Q}'$.  If we denote by $\mathcal{Q}$ the torsion-free quotient of $\mathcal{E}$ obtained by applying Lemma \ref{lemm:birationalconstruction} to $\mathcal{Q}'$, then either
\begin{enumerate}
\item $\mathcal{Q}$ is parametrized by the bounded family of  torsion-free quotients $\mathcal{G}$, or
\item $K_{X'/X} \cdot \alpha' \geq \rho (H \cdot \alpha)$.
\end{enumerate}
\end{theorem}

Note that if $\mathcal{E}$ is not semistable, then every destabilizing quotient on $X$ will be included in the bounded family $\mathcal{G}$.  In practice Theorem \ref{theo:semistabletwotypes} is most useful when $\mathcal{E}$ is $\alpha$-semistable since in this case the condition on the slopes of $\mathcal{Q}'$ is more restrictive.

\begin{proof}
Fix a surjection $\mathcal{O}_{X}(L)^{\oplus d} \to \mathcal{E}$.
For each possible rank $r = 1, 2, \ldots, n$, applying Theorem \ref{theo:birationaldegreebound3} to our chosen data with $R=0$ yields a bounded family of sheaves $\mathcal{G}_{r}$ and a constant $\rho_{r}$.  We let $\mathcal{G}$ denote the union of the bounded families $\cup_{r=1}^{\rk(\mathcal{E})} \mathcal{G}_{r}$ and set $\rho' = \inf_{r} \rho_{r}$.

We next enlarge $\mathcal{C}$ slightly: choose a closed convex full-dimensional cone $\mathcal{C}' \subset N_{1}(X)_{\mathbb{R}}$ such that $\mathcal{C'} \backslash \{ 0 \}$ is contained in the interior of $\Nef_{1}(X)$ and $\mathcal{C} \backslash \{0\}$ is contained in the interior of $\mathcal{C}'$.  Let $\mathcal{T}'$ denote all the classes $\gamma \in \mathcal{C}'$ satisfying $H \cdot \gamma = 1$.  
Since this set is compact, $c_{1}(\mathcal{E}) \cdot \gamma$ achieves its maximum as we vary $\gamma \in \mathcal{T}'$ and we define $V = \sup\{0,\max_{\gamma \in \mathcal{T}'} c_{1}(\mathcal{E}) \cdot \gamma\}$.  By applying Theorem \ref{theo:grothendieckbounded} we see that there is a bounded family of quotients $\mathcal{Q}$ of $\mathcal{E}$ satisfying $c_{1}(\mathcal{Q}) \cdot \gamma \leq V$ for some non-zero $\gamma \in \mathcal{T}'$. 
We enlarge $\mathcal{G}$ by including such quotients.  Note that for every $\gamma \in \mathcal{T}'$ and every $\mathcal{Q}$ not parametrized by $\mathcal{G}$ we have $\frac{1}{\rk(\mathcal{Q})} c_{1}(\mathcal{Q}) \cdot \gamma > \frac{1}{n} V \geq \frac{1}{n}c_{1}(\mathcal{E}) \cdot \gamma$. 
Since the slope is homogeneous with respect to rescaling $\gamma$, we conclude that $\mu_{\gamma}(\mathcal{Q}) > \mu_{\gamma}(\mathcal{E})$ for every non-zero $\gamma \in \mathcal{C}'$ and every quotient $\mathcal{Q}$ not parametrized by $\mathcal{G}$.

Let $\mathcal{T} = \mathcal{T}' \cap \mathcal{C}$. 
There is some positive constant $q$ such that $\gamma - q\beta \in \mathcal{C}'$ for every $\gamma \in \mathcal{T}$.  In particular, this means that
\begin{equation*}
\left( c_{1}(\mathcal{Q}) - \frac{\rk(\mathcal{Q})}{n} c_{1}(\mathcal{E}) \right) \cdot \gamma \geq q \left( c_{1}(\mathcal{Q}) - \frac{\rk(\mathcal{Q})}{n} c_{1}(\mathcal{E})  \right) \cdot \beta
\end{equation*}
for every $\gamma \in \mathcal{T}$ and every quotient $\mathcal{Q}$ not parametrized by $\mathcal{G}$.  Rescaling to allow $H \cdot \gamma$ to be arbitrary, we see that for every non-zero $\gamma \in \mathcal{C}$ and every quotient $\mathcal{Q}$ not parametrized by $\mathcal{G}$ we have
\begin{equation*}
\frac{\left( c_{1}(\mathcal{Q}) - \frac{\rk(\mathcal{Q})}{n} c_{1}(\mathcal{E}) \right)  \cdot \gamma}{\left( c_{1}(\mathcal{Q}) - \frac{\rk(\mathcal{Q})}{n} c_{1}(\mathcal{E}) \right)  \cdot \beta} \geq q (H \cdot \gamma).
\end{equation*}

Suppose we have an $\alpha'$-destabilizing torsion-free quotient $\mathcal{Q}'$ of $(\phi^{*}\mathcal{E})_{tf}$ of rank $r$ as in the statement of the theorem. In particular this means that
\begin{align*}
\mu_{\alpha'}(\mathcal{Q}') & < \mu_{\alpha'}((\phi^{*}\mathcal{E})_{tf}) \\
& \leq \mu_{\alpha}(\mathcal{E})
\end{align*}
where the second inequality is a consequence of Lemma \ref{lemm:furtherpullback}.
Applying Theorem \ref{theo:birationaldegreebound3} with $R=0$ to the quotient $\mathcal{E} \to \mathcal{Q}$ corresponding to $\mathcal{Q}'$ as in Lemma \ref{lemm:birationalconstruction}, we find that either:
\begin{enumerate}[(a)]
\item $\mathcal{Q}$ is parametrized by our bounded family of torsion-free quotients $\mathcal{G}$, or
\item we have
\begin{align*}
K_{X'/X} \cdot \alpha' & \geq \rho' \frac{  \left(c_{1}(\mathcal{Q}) -  \frac{r}{n} c_{1}(\mathcal{E}) \right) \cdot \alpha}{\left(c_{1}(\mathcal{Q}) -  \frac{r}{n} c_{1}(\mathcal{E}) \right)   \cdot \beta } \\
& \geq \rho' q (H \cdot \alpha)
\end{align*}
\end{enumerate}
We conclude by setting $\rho = \rho' q$.
\end{proof}

\subsection{Controlling bounded families of quotients}
In both Theorem \ref{theo:birationaldegreebound3} and Theorem \ref{theo:semistabletwotypes} there is a bounded family of quotients $\mathcal{E} \to \mathcal{Q}$ to which our bounds do not apply.  We next prove several lemmas intended to handle this bounded family.  The key result is Corollary \ref{coro:singlee} which shows that for any bounded family of quotients there exists a single exceptional divisor $E$ on a family of birational models of $X$ that controls the birational difference in slope panels.

Suppose we have a birational morphism $\phi: X' \to X$ of smooth projective varieties and a sequence of quotients on $X$:
\begin{equation*}
\mathcal{E} = \mathcal{Q}_{s} \to \mathcal{Q}_{s-1} \to \ldots \to \mathcal{Q}_{1} \to \mathcal{Q}_{0} = 0
\end{equation*}
Let $(\mathcal{Q}'_{\bullet})$ denote the birational transform on $X'$.  As discussed earlier, the graded pieces of $(\mathcal{Q}'_{\bullet})$ on $X'$ may not be the the torsion-free parts of the pullbacks of the graded pieces of $(\mathcal{Q}_{\bullet})$ on $X$.  However we can bound the differences in slope between these two different constructions.

\begin{lemma} \label{lemm:birationaltransslopebound}
Let $X$ be a smooth projective variety and let $\mathcal{E}$ denote a  non-zero torsion-free sheaf on $X$.  Suppose that $\phi: X' \to X$ is a birational map from a smooth projective $X'$ and that $\alpha' \in \Nef_{1}(X')$ is non-zero.  Set $\alpha = \phi_{*}\alpha'$ and $\mathcal{E}' = (\phi^{*}\mathcal{E})_{tf}$.

Let $(\mathcal{Q}'_{\bullet})_{i=1}^{s}$ denote a sequence of quotients of $\mathcal{E}'$ and let $(\mathcal{Q}_{\bullet})_{i=1}^{s}$ denote the birational pushforward sequence on $X$.  Define
\begin{itemize}
\item $\mathcal{T}'_{i} = \ker(\mathcal{Q}'_{i} \to \mathcal{Q}'_{i-1})$,
\item $\mathcal{T}_{i} = \ker(\mathcal{Q}_{i} \to \mathcal{Q}_{i-1})$,
\item $\widetilde{\mathcal{T}}_{i}$ to be the birational transform of $\mathcal{T}_{i}$ on $X'$.
\end{itemize}
Then
\begin{equation*}
\Vert \SP_{X,\alpha}(\mathcal{E};\mathcal{Q}_{\bullet}) - \SP_{X',\alpha'}(\mathcal{E}';\mathcal{Q}'_{\bullet}) \Vert_{sup} \leq \sum_{i=1}^{s} \left( \prod_{k=i}^{s-1} \rk(\mathcal{Q}_{k})\right)\rk(\mathcal{T}_{i}) \cdot | \mu_{\alpha}(\mathcal{T}_{i}) - \mu_{\alpha'}(\widetilde{\mathcal{T}}_{i})|.
\end{equation*}
where the empty product is taken to be $1$.
\end{lemma}

Note that the left-hand side of the equation above is the same as $\sup_{i} \{ | \mu_{\alpha}(\mathcal{T}_{i}) - \mu_{\alpha'}(\mathcal{T}'_{i}) | \}$ so that the theorem does indeed control the difference in slopes between the graded pieces $\mathcal{T}'_{i}$ of the birational transform and the birational transforms $\widetilde{\mathcal{T}}_{i}$ of the graded pieces.

\begin{proof}
We prove this by induction on $s$.  For the base case $s=1$ we have $\mathcal{T}_{1}' = (\phi^{*}\mathcal{E})_{tf} = \widetilde{\mathcal{T}}_{1}$ and the statement is immediately true. 
For a fixed $s > 1$, we prove by induction that for $1 \leq j \leq s$ we have
\begin{equation*}
| \mu_{\alpha'}(\mathcal{T}'_{j}) - \mu_{\alpha}(\mathcal{T}_{j}) |  \leq \sum_{i=1}^{j} \left( \prod_{k=i}^{j-1} \rk(\mathcal{Q}_{k})\right) \rk(\mathcal{T}_{i}) \cdot  | \mu_{\alpha}(\mathcal{T}_{i}) - \mu_{\alpha'}(\widetilde{\mathcal{T}}_{i})|.
\end{equation*}
These inequalities for $1 \leq j \leq s$ collectively imply the desired statement. 
 For the base case $j=1$ we have $\mathcal{T}_{1}' = \mathcal{Q}'_{1} = \widetilde{\mathcal{T}}_{1}$ implying the statement.

In general, consider the diagram with exact rows
\begin{equation*}
\xymatrix
   { & \phi^{*}\mathcal{T}_{j} \ar[d]_{\psi} \ar[r] & \phi^{*}\mathcal{Q}_{j} \ar[d] \ar[r] & \phi^{*}\mathcal{Q}_{j-1} \ar[d] \ar[r] & 0  \\
     0 \ar[r] & \mathcal{T}'_{j}  \ar[r] & (\phi^{*}\mathcal{Q}_{j})_{tf} \ar[r] & (\phi^{*}\mathcal{Q}_{j-1})_{tf} \ar[r] & 0 
   }
\end{equation*}
Note that:
\begin{enumerate}
\item $\mathcal{T}_{j}$, $\mathcal{T}'_{j}$, and $\widetilde{\mathcal{T}}_{j}$ all have the same rank (since they are isomorphic over an open subset).
\item The map $\psi$ factors through $(\phi^{*}\mathcal{T}_{j})_{tf} \cong \widetilde{\mathcal{T}}_{j}$.  Since the sheaves $\widetilde{\mathcal{T}}_{j}$ and $\mathcal{T}'_{j}$ have the same rank and are torsion-free, the induced $(\phi^{*}\mathcal{T}_{j})_{tf} \to \mathcal{T}'_{j}$ is injective.
\end{enumerate}
By combining the logic above with the Snake Lemma, we conclude that the cokernel of $\widetilde{\mathcal{T}}_{j} \to \mathcal{T}'_{j}$ is a torsion sheaf which admits a surjection from $(\phi^{*}\mathcal{Q}_{j-1})_{tors}$.  Thus
\begin{align*}
| \mu_{\alpha'}(\mathcal{T}'_{j}) - \mu_{\alpha}(\mathcal{T}_{j}) | & \leq | \mu_{\alpha'}(\mathcal{T}'_{j}) - \mu_{\alpha'}(\widetilde{\mathcal{T}}_{j}) | + |\mu_{\alpha'}(\widetilde{\mathcal{T}}_{j}) - \mu_{\alpha}(\mathcal{T}_{j}) | \\
& \leq \frac{\rk(\mathcal{Q}_{j-1})}{\rk(\mathcal{T}_{j})} \cdot | \mu_{\alpha'}(\mathcal{Q}'_{j-1}) - \mu_{\alpha}(\mathcal{Q}_{j-1}) |  + |\mu_{\alpha'}(\widetilde{\mathcal{T}}_{j}) - \mu_{\alpha}(\mathcal{T}_{j}) |
\end{align*}
By applying the combinatorial Lemma \ref{lemm:mediantcomputation} where $a_{i} = c_{1}(\mathcal{T}_{i}) \cdot \alpha$, $a'_{i} = c_{1}(\mathcal{T}'_{i}) \cdot \alpha'$, and $b_{i} = \rk(\mathcal{T}_{i})$, we obtain $| \mu_{\alpha'}(\mathcal{Q}'_{j-1}) - \mu_{\alpha}(\mathcal{Q}_{j-1}) | \leq  \Vert \SP_{X',\alpha'}(\mathcal{Q}'_{j-1}, \mathcal{Q}'_{\bullet}) - \SP_{X,\alpha}(\mathcal{Q}_{j-1}, \mathcal{Q}_{\bullet}) \Vert_{sup}$.  Continuing the chain of inequalities:
\begin{align*}
| \mu_{\alpha'}(\mathcal{T}'_{j}) - \mu_{\alpha}(\mathcal{T}_{j}) | & \leq \frac{ \rk(\mathcal{Q}_{j-1})}{\rk(\mathcal{T}_{j})} \cdot\Vert \SP_{X',\alpha'}(\mathcal{Q}'_{j-1}, \mathcal{Q}'_{\bullet}) - \SP_{X,\alpha}(\mathcal{Q}_{j-1}, \mathcal{Q}_{\bullet}) \Vert_{sup}  + |\mu_{\alpha'}(\widetilde{\mathcal{T}}_{j}) - \mu_{\alpha}(\mathcal{T}_{j}) | \\
& \leq \frac{1}{\rk(\mathcal{T}_{j})}  \sum_{i=1}^{j} \left( \prod_{k=i}^{j-1} \rk(\mathcal{Q}_{k})\right) \rk(\mathcal{T}_{i}) \cdot | \mu_{\alpha}(\mathcal{T}_{i}) - \mu_{\alpha'}(\widetilde{\mathcal{T}}_{i})|.
\end{align*}
where the last line follows from the induction assumption on $s$.  Since $\rk(\mathcal{T}_{j}) \geq 1$, this finishes the inductive step.
\end{proof}

\begin{lemma} \label{lemm:mediantcomputation}
Suppose we have sets $\{ (a_{i},b_{i}) \}_{i=1}^{t}$ and $\{ (a'_{i},b_{i}) \}_{i=1}^{t}$ where the $a_{i}, a_{i}'$ are integers and the $b_{i}$ are positive integers.  Then
\begin{equation*}
\left| \frac{\sum_{i=1}^{t} a_{i}' - \sum_{i=1}^{t} a_{i}}{\sum_{i=1}^{t} b_{i}} \right| \leq \sup_{i=1,\ldots,t} \frac{|a_{i}'-a_{i}|}{b_{i}}.
\end{equation*}
\end{lemma}

\begin{proof}
Choose the index $k \in \{1,\ldots,t\}$ which maximizes the value of $\frac{|a_{k}'-a_{k}|}{b_{k}}$.  This implies that for every $i$ we have $b_{k}|a_{i}'-a_{i}| \leq b_{i}|a_{k}' - a_{k}|$.  Thus
\begin{align*}
b_{k} \left| \sum_{i=1}^{t} a_{i}' - \sum_{i=1}^{t} a_{i}  \right| & \leq \sum_{i=1}^{t} b_{k} |a_{i}' - a_{i}| \\
& \leq \sum_{i=1}^{t} b_{i} |a_{k}' - a_{k}|
\end{align*}
Rearranging gives
\begin{equation*}
\left| \frac{\sum_{i=1}^{t} a_{i}' - \sum_{i=1}^{t} a_{i}}{\sum_{i=1}^{t} b_{i}} \right| \leq \frac{|a_{k}'-a_{k}|}{b_{k}} = \sup_{i=1,\ldots,t} \frac{|a_{i}'-a_{i}|}{b_{i}}.
\end{equation*}
\end{proof}

Next suppose we take a sheaf $\mathcal{E}$ and its birational transform $\mathcal{E}'$.  The following lemma allows us to control the difference in slope panels between $\mathcal{E},\mathcal{E}'$ using the difference in slope panels of the graded pieces of their Harder-Narasimhan filtrations.

\begin{lemma} \label{lemm:allintvertical}
Let $X$ be a smooth projective variety and let $\mathcal{E}$ denote a non-zero torsion-free sheaf on $X$ of rank $n$.  Suppose that $\phi: X' \to X$ is a birational map from a smooth projective $X'$ and that $\alpha' \in \Nef_{1}(X')$ is non-zero.  Set $\alpha = \phi_{*}\alpha'$ and $\mathcal{E}' = (\phi^{*}\mathcal{E})_{tf}$.

Let $(\mathcal{Q}'_{\bullet})_{i=1}^{s}$ denote the $\alpha'$-Harder-Narasimhan filtration of $\mathcal{E}'$ and $(\mathcal{Q}_{\bullet})$ denote the birational pushforward quotient sequence on $X$.  Let $(\mathcal{R}_{\bullet})_{i=1}^{t}$ denote the $\alpha$-Harder-Narasimhan filtration of $\mathcal{E}$ and let $(\mathcal{R}'_{\bullet})$ denote the birational transform quotient sequence on $X'$.  Define
\begin{equation*}
N = \sup \left\{ \Vert \SP_{X,\alpha}(\mathcal{E};\mathcal{Q}_{\bullet}) - \SP_{X',\alpha'}(\mathcal{E}';\mathcal{Q}'_{\bullet}) \Vert_{sup}, \Vert \SP_{X,\alpha}(\mathcal{E};\mathcal{R}_{\bullet}) - \SP_{X',\alpha'}(\mathcal{E}';\mathcal{R}'_{\bullet}) \Vert_{sup} \right\}.
\end{equation*}
Then
\begin{equation*}
\Vert \SP_{X,\alpha}(\mathcal{E}) - \SP_{X',\alpha'}(\mathcal{E}') \Vert_{sup} \leq  \left( 2^{n}-1 \right) N.
\end{equation*}
\end{lemma}

\begin{proof}
For convenience we define
\begin{align*}
(a_{1},\ldots,a_{n}) & = \SP_{X,\alpha}(\mathcal{E}, \mathcal{R}_{\bullet})  \\
(b_{1},\ldots,b_{n}) & = \SP_{X,\alpha}(\mathcal{E};\mathcal{Q}_{\bullet}) \\
(a'_{1},\ldots,a'_{n}) & = \SP_{X',\alpha'}(\mathcal{E}',\mathcal{R}'_{\bullet})  \\
(b'_{1},\ldots,b'_{n}) & = \SP_{X',\alpha'}(\mathcal{E}';\mathcal{Q}'_{\bullet})
\end{align*}
We prove inductively that for any $1 \leq j \leq n$ we have
\begin{equation*}
|a_{j} - b_{j}'| \leq  \left(2^{j} - 1 \right) N
\end{equation*}
To prove an upper bound on $a_{j} - b'_{j}$, we apply Lemma \ref{lemm:hnmaininequality} on $X'$ to conclude that
\begin{equation*}
\sum_{i=1}^{j} (b'_{i} - a'_{i}) \geq 0.
\end{equation*}
Rearranging, we obtain
\begin{align*}
a'_{j} - b'_{j} \leq \sum_{i=1}^{j-1} (b'_{i} - a'_{i}) & =  \sum_{i=1}^{j-1}  (b'_{i}-a_{i}) + \sum_{i=1}^{j-1} (a_{i} - a'_{i}) \\
& \leq \left( \sum_{i=1}^{j-1}  (b'_{i}-a_{i}) \right) + (j-1)N
\end{align*}
Combining with the induction assumption, we see that
\begin{align*}
a_{j} - b'_{j} \leq N + a'_{j} - b'_{j} \leq jN + \sum_{i=1}^{j-1} (2^{i}-1) N = (2^{j}-1) N 
\end{align*}
To prove a lower bound on $a_{j} - b'_{j}$, we apply Lemma \ref{lemm:hnmaininequality} on $X$ to conclude that
\begin{equation*}
\sum_{i=1}^{j} (a_{i} - b_{i}) \geq 0.
\end{equation*}
Arguing as before, we have
\begin{align*}
a_{j} - b_{j} \geq \sum_{i=1}^{j-1}  (b_{i}-a_{i}) & = \sum_{i=1}^{j-1}  (b'_{i}-a_{i}) - \sum_{i=1}^{j-1} (b'_{i} - b_{i}) \\
& \geq \left( \sum_{i=1}^{j-1}  (b'_{i}-a_{i}) \right) - (j-1)N
\end{align*}
 Thus
\begin{align*}
a_{j} - b'_{j} \geq a_{j} - b_{j} - N \geq -jN - \sum_{i=1}^{j-1} (2^{i}-1) N  =  -(2^{j}-1) N 
\end{align*}
\end{proof}

Combining Lemma \ref{lemm:birationaltransslopebound} and Lemma \ref{lemm:allintvertical}, we obtain a result which allows us to control the changes in slopes of Harder-Narasimhan filtrations corresponding to any bounded family of quotients of $\mathcal{E}$.  As mentioned before, the key is the existence of the single divisor $E$ controlling the differences on slope panels.

\begin{corollary} \label{coro:singlee}
Let $X$ be a smooth projective variety and let $\mathcal{E}$ be a non-zero torsion-free sheaf on $X$.  Fix a bounded family of sequences $(\mathcal{Q}_{\bullet})$ of torsion-free quotients of $\mathcal{E}$ and a finite set $\{ (\mathcal{R}^{\ell}_{\bullet}) \}_{\ell \in L}$ of sequences of torsion-free quotients of $\mathcal{E}$.  Let $\psi: \mathfrak{Y} \to X \times S$ denote the family of birational models obtained by applying Construction \ref{cons:resolvingfamiliesofsequences} to construct locally free birational transforms of the sheaves in the bounded family $(\mathcal{Q}_{\bullet})$  (where we use the shorthand $S = \sqcup S_{k}$ and $\mathfrak{Y} = \sqcup \mathfrak{Y}_{k}$).   Then there is an effective $\psi$-exceptional divisor $E$ on $\mathfrak{Y}$ satisfying the following property.

Suppose that for some closed point $s \in S$ there is a non-zero nef class $\alpha_{s} \in \Nef_{1}(\mathfrak{Y}_{s})$ with the following properties:
\begin{itemize}
\item The $\alpha_{s}$-Harder-Narasimhan filtration of $(\psi_{s}^{*}\mathcal{E})_{tf}$ is given by the birational transform $(\mathcal{Q}'_{s,\bullet})$ of a quotient sequence $(\mathcal{Q}_{\bullet})$ in our bounded family.
\item Setting $\alpha = \psi_{s*}\alpha_{s}$, the $\alpha$-Harder-Narasimhan filtration of $\mathcal{E}$ is given by one of the sequences $(\mathcal{R}^{\ell}_{\bullet})$.
\end{itemize}
Then we have
\begin{align*}
\Vert \SP_{X,\alpha}(\mathcal{E}) - \SP_{\mathfrak{Y}_{s},\alpha_{s}}((\psi_{s}^{*}\mathcal{E})_{tf}) \Vert_{sup} \leq E|_{\mathfrak{Y}_{s}} \cdot \alpha_{s}.
\end{align*}
\end{corollary}

\begin{proof}
Let $\mathcal{E}' = (\psi^{*}\pi_{1}^{*}\mathcal{E})_{tf}$ and let $(\mathcal{R}'^{\ell}_{\bullet})$ denote the birational transform on $\mathfrak{Y}$ of the pullback $(\pi_{1}^{*}\mathcal{R}^{\ell}_{\bullet})$ under the projection map $\pi_{1}: X \times S \to X$.  Retaining the notation of Construction \ref{cons:resolvingfamiliesofsequences}, define
\begin{itemize}
\item $\mathcal{T}_{i}^{k} := \ker(\mathcal{Q}_{i}^{k} \to \mathcal{Q}_{i-1}^{k})$ on $X \times S_{k}$,
\item $\widetilde{\mathcal{T}}_{i}^{k} := (\psi^{*}\mathcal{T}_{i}^{k})_{tf}$ on $\mathfrak{Y}_{k}$,
\item $\mathcal{H}^{\ell}_{i} =  \ker(\pi_{1}^{*}\mathcal{R}^{\ell}_{i} \to \pi_{1}^{*}\mathcal{R}^{\ell}_{i-1})$ on $X \times S$,
\item $\widetilde{\mathcal{H}}^{\ell}_{i} =  (\psi^{*}\mathcal{H}^{\ell}_{i})_{tf}$ on $\mathfrak{Y}$.
\end{itemize}
As we vary $i,k$ the divisors $\psi^{*}c_{1}(\mathcal{T}_{i}^{k}) - c_{1}(\widetilde{\mathcal{T}}_{i}^{k})$ form a finite set of effective $\psi$-exceptional divisors on $\mathfrak{Y}$ which we denote by $\{ F_{j} \}$.  Similarly, as we vary $\ell,i$ the divisors $\psi^{*}c_{1}(\mathcal{H}^{\ell}_{i})  - c_{1}(\widetilde{\mathcal{H}}^{\ell}_{i})$ form a finite set of effective $\psi$-exceptional divisors on $\mathfrak{Y}$ which we denote by $\{ G_{p} \}$.

By Lemma \ref{lemm:allintvertical} we have 
\begin{align*}
\Vert \SP_{X,\alpha}(\mathcal{E}) - \SP_{\mathfrak{Y}_{s},\alpha_{s}}((\psi_{s}^{*}\mathcal{E})_{tf}) \Vert_{sup} \leq (2^{\rk(\mathcal{E})} - 1)N
\end{align*}
where
\begin{equation*}
N = \sup \left\{ \Vert \SP_{X,\alpha}(\mathcal{E};\mathcal{Q}_{\bullet}) - \SP_{\mathfrak{Y}_{s},\alpha_{s}}(\mathcal{E}'_{s};\mathcal{Q}'_{s,\bullet}) \Vert_{sup}, \Vert \SP_{X,\alpha}(\mathcal{E};\mathcal{R}^{\ell}_{\bullet}) - \SP_{\mathfrak{Y}_{s},\alpha_{s}}(\mathcal{E}'_{s};\mathcal{R}'^{\ell}_{s,\bullet}) \Vert_{sup} \right\}
\end{equation*}
In turn, Lemma \ref{lemm:birationaltransslopebound} shows that $N$ is bounded above by intersections of $\alpha_{s}$ against certain linear combinations of $F_{j}|_{\mathfrak{Y}_{s}}$ and $G_{p}|_{\mathfrak{Y}_{s}}$.  Since $\alpha_{s}$ is nef, it suffices to choose $E$ to be an effective $\psi$-exceptional divisor that is more effective than all of the finitely many possible linear combinations of the various divisors $F_{j}$ and $G_{p}$ obtained by combining the results of Lemma \ref{lemm:birationaltransslopebound} and Lemma \ref{lemm:allintvertical}.
\end{proof}

\subsection{Canonical degree bound}
We are finally prepared to prove the main theorem which gives us control of the relative canonical degree of a nef class $\alpha'$ for which the pullback of $\mathcal{E}$ admits a destabilizing quotient.

\begin{theorem} \label{theo:kxbound}
Let $X$ be a smooth projective variety and let $\mathcal{E}$ be a non-zero torsion-free sheaf on $X$.  Fix a big and nef divisor $H$.  Fix a closed convex cone $\mathcal{C}$ such that $\mathcal{C} \backslash \{0\}$ is contained in the interior of $\Nef_{1}(X)$ and $\mathcal{C}$ contains a class $\beta \in N_{1}(X)_{\mathbb{Z}}$.  There are positive constants $\rho, \eta$ which satisfy the following property.  

Let $\phi: X' \to X$ denote a birational morphism from a smooth projective variety and let $\alpha' \in \Nef_{1}(X')$ be a non-zero nef class such that $\alpha := \phi_{*}\alpha'$ lies in $\mathcal{C}$.  Define the constant $\mu$ as
\begin{equation*}
\mu = \Vert \SP_{X,\alpha}(\mathcal{E}) - \SP_{X',\alpha'}((\phi^{*}\mathcal{E})_{tf}) \Vert_{sup}.
\end{equation*}
Then 
\begin{equation*}
K_{X'/X} \cdot \alpha' \geq \min\{ \rho(H \cdot \alpha), \eta  \mu \}.
\end{equation*}
\end{theorem}

\begin{proof}
By enlarging $\mathcal{C}$ slightly we may suppose that it is a rational polyhedral cone that is still interior to $\Nef_{1}(X)$.  By \cite[Lemma 3.3.3]{Neumann10} there is a finite set of filtrations of $\mathcal{E}$ which can occur as Harder-Narasimhan filtrations as we vary $\alpha \in \mathcal{C}$.  We denote the corresponding finite set of quotient sequences by $\{ (\mathcal{R}^{\ell}_{\bullet}) \}_{\ell \in L}$.

Apply Theorem \ref{theo:semistabletwotypes} to $\mathcal{E}$ to obtain a positive constant $\rho$ and a bounded family of sheaves $\mathcal{G}$. 
According to the conclusion of Theorem \ref{theo:semistabletwotypes}, one of the three following possibilities must hold:
\begin{itemize}
\item $(\phi^{*}\mathcal{E})_{tf}$ is $\alpha'$-semistable,
\item if $(\phi^{*}\mathcal{E})_{tf} \to \mathcal{Q}'$ is an $\alpha'$-destabilizing quotient, then $\mathcal{Q}'$ is the birational transform of a quotient parametrized by $\mathcal{G}$, or
\item $K_{X'/X} \cdot \alpha' \geq \rho (H \cdot \alpha)$.
\end{itemize}
We add the trivial quotient $\mathcal{E} \to 0$ to the bounded family of sheaves $\mathcal{G}$.  Referring to the three possible outcomes of Theorem \ref{theo:semistabletwotypes}, we see that either:
\begin{enumerate}
\item every quotient $(\mathcal{Q}'_{\bullet})$ in the $\alpha'$-Harder-Narasimhan sequence of $(\phi^{*}\mathcal{E})_{tf}$ is the birational transform of a quotient parametrized by the bounded family $\mathcal{G}$, or
\item $K_{X'/X} \cdot \alpha' \geq \rho (H \cdot \alpha)$.
\end{enumerate}
It suffices to prove the desired linear lower bound on $K_{X'/X} \cdot \alpha'$ separately in each of these two cases.  This statement is transparently true in case (2).

In case (1), we apply Lemma \ref{lemm:boundedfrombounded} to the bounded family of sheaves $\mathcal{G}$ to get a bounded family of quotient sequences of $\mathcal{E}$ whose terms lie in $\mathcal{G}$.   We then apply Construction \ref{cons:resolvingfamiliesofsequences} to this bounded family of quotient sequences to get a family of smooth birational models $\psi: \mathfrak{Y} \to X \times S$ (where we use the shorthand $S = \sqcup S_{k}$ and $\mathfrak{Y} = \sqcup \mathfrak{Y}_{k}$).  Recall that for each closed point $s \in S$ if we consider the restriction $(\mathcal{Q}_{i}^{k}|_{X_{s}})$ and take the birational transform to $\mathfrak{Y}_{s}$ the terms $\mathfrak{Q}_{i}^{k,\prime}|_{\mathfrak{Y}_{s}}$ are locally free sheaves and the successive kernels $\ker(\mathfrak{Q}_{i}^{k,\prime}|_{\mathfrak{Y}_{s}} \to \mathfrak{Q}_{i-1}^{k,\prime}|_{\mathfrak{Y}_{s}})$ are also locally free.  Furthermore, by applying Lemma \ref{lemm:resolvingsingsofsequence} to pass to a higher birational model, we may also assume that the birational transforms of every quotient $\mathcal{R}^{\ell}_{i}$ and every successive kernel of the birational transforms are locally free.

Suppose $\phi: X' \to X$ is a birational map and $\alpha' \in \Nef_{1}(X')$ as in the statement of the theorem and that the $\alpha'$-Harder-Narasimhan filtration of $(\phi^{*}\mathcal{E})_{tf}$ is given by some sequence of quotients $(\mathcal{Q}_{i})_{i=1}^{s}$
parametrized by our bounded family.  In particular the quotient sequence $(\mathcal{Q}_{\bullet})$ corresponds to a closed point $s$ in our parameter space $S$.  Let $\psi_{s}: \mathfrak{Y}_{s} \to X$ denote the corresponding birational map and choose a birational model $\widehat{X}$ with birational maps $g_{1}: \widehat{X} \to X'$ and $g_{2}: \widehat{X} \to \mathfrak{Y}_{s}$.

For notational convenience we define $\mathcal{E}' = (\phi^{*}\mathcal{E})_{tf}$ and define the nef curve class $\alpha_{s} = g_{2*}g_{1}^{*}\alpha'$.  We claim that the $\alpha_{s}$-Harder-Narasimhan filtration of $(\psi_{s}^{*}\mathcal{E})_{tf}$ is given by the locally free sheaves $\mathfrak{Q}_{i}^{k,\prime}|_{\mathfrak{Y}_{s}}$ coming from Construction \ref{cons:resolvingfamiliesofsequences}.  Indeed, recall that all the $\mathfrak{Q}_{i}^{k,\prime}|_{\mathfrak{Y}_{s}}$ as well as their successive kernels are locally free.  This means that every term and successive kernel of the filtration of the birational transform of $(\mathcal{Q}_{\bullet})$ to $\widehat{X}$ is pulled back from $\mathfrak{Y}_{s}$.  Furthermore, we know that this birational transform sequence on $\widehat{X}$ is the $g_{1}^{*}\alpha'$-Harder-Narasimhan filtration of $g_{1}^{*}(\mathcal{E}')_{tf}$ by \cite[Proposition 2.8]{GKP16}.  By applying Corollary \ref{coro:stillhnfilt} to the birational map $g_{2}$ we see that the successive kernels of $\mathfrak{Q}_{i}^{k,\prime}|_{\mathfrak{Y}_{s}}$ are $\alpha_{s}$-semistable and have decreasing $\alpha_{s}$-slope, proving the claim.

As discussed above, \cite[Proposition 2.8]{GKP16} shows that stability is compatible with pullbacks so that $\SP_{X',\alpha'}(\mathcal{E}') = \SP_{\widehat{X},g_{1}^{*}\alpha'}(g_{1}^{*}(\mathcal{E}')_{tf})$.  In turn, since the $g_{1}^{*}\alpha'$-Harder-Narasimhan filtration of $g_{1}^{*}(\mathcal{E}')_{tf}$ is pulled back from the corresponding filtration on $\mathfrak{Y}_{s}$ we see that $\SP_{\widehat{X},g_{1}^{*}\alpha'}(g_{1}^{*}(\mathcal{E}')_{tf}) = \SP_{\mathfrak{Y}_{s},\alpha_{s}}((\psi_{s}^{*}\mathcal{E})_{tf})$.
Thus we have
\begin{align*}
\mu & = \Vert \SP_{X,\alpha}(\mathcal{E}) - \SP_{X',\alpha'}(\mathcal{E}') \Vert_{sup} \\
& = \Vert \SP_{X,\alpha}(\mathcal{E}) - \SP_{\mathfrak{Y}_{s},\alpha_{s}}((\psi_{s}^{*}\mathcal{E})_{tf}) \Vert_{sup}
\end{align*}

Since we have constructed the family $\mathfrak{Y}$ by using Construction \ref{cons:resolvingfamiliesofsequences}, we may apply Corollary \ref{coro:singlee} to see that there is an effective $\psi$-exceptional divisor $E$ (not depending on $\phi: X' \to X$) such that
\begin{align*}
\Vert \SP_{X,\alpha}(\mathcal{E}) - \SP_{\mathfrak{Y}_{s},\alpha_{s}}((\psi_{s}^{*}\mathcal{E})_{tf}) \Vert_{sup} \leq E|_{\mathfrak{Y}_{s}} \cdot \alpha_{s}.
\end{align*}
Next we choose a positive constant $\zeta$ such that $\zeta K_{\mathfrak{Y}/X \times S} \geq E$.  Thus
\begin{align*}
\mu & \leq \zeta K_{\mathfrak{Y}_{s}/X} \cdot \alpha_{s} \\
& = \zeta g_{2}^{*}K_{\mathfrak{Y}_{s}/X} \cdot g_{1}^{*}\alpha' \\
& \leq  \zeta K_{\widehat{X}/X} \cdot g_{1}^{*}\alpha' \\
& = \zeta K_{X'/X} \cdot \alpha'
\end{align*}
We obtain the desired statement by setting $\eta = 1/\zeta$.
\end{proof}

\section{Jumping loci for dominant families with connected fibers}

Let $X$ be a smooth projective variety, let $\mathcal{E}$ be a locally free sheaf on $X$, and let $M$ be an irreducible component of $\mathcal{M}_{g,0}(X)$. Given a variety $W$ equipped with a generically finite morphism $W \to M$, our goal is to give a lower bound on the codimension of the image of $W$ based on the positivity of $s^{*}\mathcal{E}$ where $s: C \to X$ is a general curve parametrized by $W$.

\subsection{Grauert-M\" ulich}

\begin{definition}
Let $Y$ be a variety and let $\mathcal{E}$ be a globally generated vector bundle on $Y$.  The syzygy bundle $M_{\mathcal{E}}$ is the kernel of the evaluation map $\OO_Y\otimes H^0(\mathcal{E})\to \mathcal{E}$.
\end{definition}

We will use the following result of \cite{Butler94} describing syzygy bundles on curves. 

\begin{theorem}[\cite{Butler94}] \label{theo:butler}
Let $\mathcal{E}$ be a globally generated locally free sheaf on a curve $C$ of genus $g$ and let $M_{\mathcal{E}}$ be its syzygy bundle.
\begin{enumerate}
\item Fix a constant $k \geq 2$.  If $\mu^{min}(\mathcal{E}) \geq kg$ then $\mu^{min}(M_{\mathcal{E}}) \geq -\frac{k}{k-1}$. 
\item If $\mu^{min}(\mathcal{E}) < 2g$ then $\mu^{min}(M_{\mathcal{E}}) \geq -2g \rk(\mathcal{E}) - 2$.
\end{enumerate}
\end{theorem}

\begin{proof}
\cite[Theorem 6.8]{LRT23} carefully explains how Butler's work implies (2) and implies (1) when either $g=0$ or when $k = 2$.  \cite[1.3 Corollary]{Butler94} proves case (1) in the remaining situations $g \geq 1$ and $k > 2$.
\end{proof}

The Grauert-M\" ulich theorem of \cite{LRT23} bounds the difference between the expected and actual values of the Harder-Narasimhan filtration of $\mathcal{E}|_{C}$.  It relies on several geometric assumptions about the evaluation map: it should be dominant with connected fibers and it should be flat along a general curve.

\begin{theorem}[{\cite[Corollary 6.6]{LRT23}}] \label{theo:hnfgm}
Let $X$ be a smooth projective variety and let $\mathcal{E}$ be a non-zero torsion-free sheaf on $X$ of rank $r$.  Let $W$ be a variety equipped with a generically finite morphism $W \to {\mathcal M}_{g,0}(X)$.  Let $p: U^{\nu} \to W$ denote the normalization of the universal family over $W$ with evaluation map $ev^{\nu}: U^{\nu} \to X$.  Assume that $ev^{\nu}$ is dominant with connected fibers and that a general fiber of $p$ is contained in the locus where $ev^{\nu}$ is flat.

Let $C$ denote a general fiber of $U^{\nu} \to W$ equipped with the induced morphism $s: C \to X$.  Let $t$ be the length of the torsion part of the normal sheaf $N_{s}$, let $\mathcal{G}$ be the subsheaf of $(N_{s})_{tf}$ generated by global sections, and let $V$ be the tangent space to $W$ at $s$. Let $q$ be the dimension of the cokernel of the composition
\begin{equation*}
V \to T_{{\mathcal M}_{g,0}(X),s} = H^{0}(C,N_{s})  \to H^{0}(C,(N_{s})_{tf}).
\end{equation*}
Then we have
\begin{equation*}
\Vert \SP_{X,[C]}(\mathcal{E}) - \SP_{C}(s^{*}\mathcal{E}) \Vert_{sup} \leq \frac{1}{2} \left( (q+1)\mu^{max}(M_{\mathcal{G}}^{\vee}) + t \right)  (\rk(\mathcal{E}) - 1).
\end{equation*}
\end{theorem}

\subsection{Codimension of flat families of curves} \label{sect:gmcodim}

Our next goal is to bound the codimension of families which satisfy an extra flatness assumption.  We first need a couple lemmas. 

\begin{lemma} \label{lemm:expdimmbar}
Let $X$ be a smooth projective variety.  Let $M$ be an irreducible component of ${\mathcal{M}}_{g,0}(X)$ parametrizing a dominant family of maps.  Then
\begin{equation*}
-K_{X} \cdot C + (\dim(X)-3)(1-g) \leq \dim(M) \leq -K_{X} \cdot C + \dim(X) + 2g - 3. 
\end{equation*}
\end{lemma}

\begin{proof}
In general the expected dimension of $M$ is $\chi(N_{s})$ and this gives a lower bound on $\dim(M)$.  An upper bound on $\dim(M)$ is given by $h^{0}(C,N_{s})$.  Since $M$ parametrizes a dominant family we know that $N_{s}$ is generically globally generated.  Thus \cite[Lemma 2.8]{LRT23} shows that
\begin{equation*}
h^{1}(C,N_{s}) \leq g (\dim(X)-1)
\end{equation*}
leading to the result.
\end{proof}

\begin{lemma} \label{lemm:torsbound}
Let $X$ be a smooth projective variety.  Let $W$ be a variety equipped with a generically finite morphism $W \to {\mathcal M}_{g,0}(X)$.  Let $p: {U}^{\nu} \to W$ denote the normalization of the universal family over $W$ with evaluation map $ev^{\nu}: U^{\nu} \to X$.  Assume that $ev^{\nu}$ is dominant and that the general morphism $s: C \to X$ parametrized by $W$ is birational onto its image.  Then the torsion part of $N_{s}$ satisfies
\begin{equation*}
h^{0}(C,(N_{s})_{tors}) \leq g (\dim(X)-1) + \codim(W).
\end{equation*}
\end{lemma}

\begin{proof}
Since $s: C \to X$ is a general point of $W$, we know that $W$ is smooth at that point.  By \cite[Corollary 6.11]{AC81}, the image of the map
\begin{equation*}
T_{W,[s]} \to T_{{\mathcal M}_{g,0}(X),[s]} = H^{0}(C,N_{s})
\end{equation*}
has zero intersection with the subgroup $H^{0}(C,(N_{s})_{tors})$.  Letting $M$ denote an irreducible component of ${\mathcal M}_{g,0}(X)$ containing the image of $W$, we have
\begin{align*}
h^{0}(C,(N_{s})_{tors}) & \leq \dim( T_{{\mathcal M}_{g,0}(X),[s]}) - \dim(T_{W,[s]}) \\
& = \left( \dim(T_{{\mathcal M}_{g,0}(X),[s]}) - \dim(M) \right) + (\dim(M) - \dim(W)) \\
& \leq h^{1}(C,N_{s}) + \codim_{M}(W).
\end{align*}
Since $W$ defines a dominant family we know that $N_{s}$ is generically globally generated.  Then \cite[Lemma 2.8]{LRT23} shows that $h^{1}(C,N_{s}) \leq g (\dim(X)-1)$.
\end{proof}

We can now prove the main theorem in this subsection.

\begin{theorem} \label{theo:tangentgaps}
Let $X$ be a smooth projective variety.  Let $W$ be a variety equipped with a generically finite morphism $W \to{\mathcal M}_{g,0}(X)$.  Let $p: {U}^{\nu} \to W$ denote the normalization of the universal family over $W$ with evaluation map $ev^{\nu}: U^{\nu} \to X$.  Assume that
\begin{itemize}
\item $ev^{\nu}$ is dominant with connected fibers,
\item a general map parametrized by $W$ is birational onto its image,
and
\item a general fiber of $p$ is contained in the locus where $ev^{\nu}$ is flat.
\end{itemize}
Let $\mathcal{E}$ be a non-zero torsion-free sheaf on $X$.  Define
\begin{equation*}
\mu = \Vert \SP_{X,s_{*}[C]}(\mathcal{E}) - \SP_{C}((s^{*}\mathcal{E})_{tf}) \Vert_{sup}
\end{equation*}
where $s: C \to X$ is a general map parametrized by $W$.  Letting $M$ denote an irreducible component of ${\mathcal M}_{g,0}(X)$ containing the image of $W$, we have
\begin{equation*}
\codim_{M}(W) \geq \frac{\mu}{(\gamma+1) (\rk(\mathcal{E})-1)} - \gamma
\end{equation*}
where $\gamma = g (\dim(X)-1) + 1$. 
\end{theorem}

\begin{proof}
Letting $t$ denote the length of $(N_{s})_{tors}$, Lemma \ref{lemm:torsbound} implies that
\begin{equation*}
t \leq g (\dim(X)-1) + \codim_{M}(W).
\end{equation*}

For $s$ general in $W$ the normal sheaf $N_{s}$ is generically globally generated.  Thus \cite[Lemma 2.8]{LRT23} shows that $h^{1}(C,N_{s}) \leq g (\dim(X)-1)$.  
We conclude that $\dim(M) \geq \dim(T_{M,[s]}) - g (\dim(X)-1)$.
Since a general map $s: C \to W$ will define a smooth point of $W$, we have
\begin{equation*}
\dim(M) - \dim(W) + g (\dim(X)-1) \geq \dim \mathrm{coker}(T_{W,[s]} \to T_{M,[s]}).
\end{equation*}
Since $H^{0}(C,N_{s})  \to H^{0}(C,(N_{s})_{tf})$ is surjective, the dimension of the cokernel of the composed map
\begin{equation*}
T_{W,[s]} \to T_{ {\mathcal M}_{g,0}(X),s} = H^{0}(C,N_{s})  \to H^{0}(C,(N_{s})_{tf})
\end{equation*}
is also bounded above by the sum $\codim_{M}(W) + g (\dim(X)-1)$.

Finally, letting $\mathcal{G}$ denote the globally generated subsheaf of $(N_{s})_{tf}$, Theorem \ref{theo:butler} shows that $\mu^{min}(M_{\mathcal{G}}) \geq -2\gamma$ (regardless of whether we are in case (1) or case (2) of Theorem \ref{theo:butler}).
Putting all the bounds together, Theorem \ref{theo:hnfgm} shows that $\Vert \SP_{X,s_{*}[C]}(\mathcal{E}) - \SP_{C}(s^{*}\mathcal{E}) \Vert$ is bounded above by
\begin{equation*}
 \frac{1}{2}   \Bigl( (\codim_{M}(W)   + g (\dim(X)-1) +1)(2\gamma)
  + g (\dim(X)-1) + \codim_{M}(W)  \Bigl)  (\rk(\mathcal{E})-1).
\end{equation*}
Rearranging shows that
\begin{equation*}
\codim_{M}(W) \geq \frac{2\mu}{(2\gamma+1) (\rk(\mathcal{E})-1)} - \frac{2\gamma^{2}+\gamma - 1}{2\gamma+1}
\end{equation*}
which is slightly stronger than the given statement.
\end{proof}

 If we add a few more assumptions on the curves, we can get much nicer bounds:

\begin{corollary} \label{coro:cleantangentgaps}
Suppose we are in the setting of Theorem \ref{theo:tangentgaps}.  Furthermore assume that a general $s: C \to X$ parametrized by $W$ satisfies
\begin{itemize}
\item $s$ is an immersion, and
\item the normal sheaf $N_{s}$ has $\mu^{min}(N_{s}) \geq kg$ for some constant $k \geq 2$.
\end{itemize}
Then
\begin{equation*}
\codim_{M}(W) \geq \frac{2(k-1)}{(\rk(\mathcal{E})-1)k} \cdot \mu - 1.
\end{equation*}
\end{corollary}

\begin{proof}
Note that our two additional assumptions will also hold true for the general map parametrized by $M$.  In particular, the degree assumption on $N_{s}$ implies that $h^{1}(C,N_{s}) = 0$ and thus $M$ must be generically smooth of the expected dimension.  By assumption $N_{s}$ is torsion-free.  By repeating the argument in Theorem \ref{theo:tangentgaps} and appealing to Theorem \ref{theo:butler} for an improved bound on the slope of the syzygy bundle we obtain the desired statement.
\end{proof}

\subsection{Codimension of nonflat families of curves} 

By combining Grauert-M\" ulich with Theorem \ref{theo:kxbound}, we can prove our most general statement addressing the codimension of a family of curves based on slope panels.

\begin{construction}[Flattening construction] \label{cons:flatteningfamilyofcurves}
Let $X$ be a smooth projective variety.  Let $W$ be a variety equipped with a generically finite morphism $W \to \mathcal{M}_{g,0}(X)$.  Let $p: {U}^{\nu} \to W$ denote the normalization of the universal family over $W$ with evaluation map $ev^{\nu}: U^{\nu} \to X$.  Assume that $ev^{\nu}$ is dominant.

Then there is a birational map $\phi: X' \to X$ from a smooth projective variety $X'$ and a non-empty open subset $W' \subset W$ such that the preimage $U^{\nu,\circ} := p^{-1}W'$ admits a flat morphism $ev': U^{\nu,\circ} \to X'$ satisfying $ev^{\nu}|_{U^{\nu,\circ}} = \phi \circ ev'$.  This follows from standard facts about flattening birational maps (see \cite[Construction 4.2]{LRT24} for details).
\end{construction}

\begin{theorem} \label{theo:maincodim}
Let $X$ be a smooth projective variety and let $\mathcal{E}$ be a non-zero torsion-free sheaf on $X$.  Fix a big and nef divisor $H$.  Fix a closed convex cone $\mathcal{C}$ such that $\mathcal{C} \backslash \{0\}$ is contained in the interior of $\Nef_{1}(X)$ and $\mathcal{C}$ contains a class $\beta \in N_{1}(X)_{\mathbb{Z}}$.  There are affine linear functions $L,S: \mathbb{R} \to \mathbb{R}$ with positive leading coefficients which satisfy the following property.

Let $W$ be a variety equipped with a generically finite morphism $W \to {\mathcal M}_{g,0}(X)$.  Let $p: {U}^{\nu} \to W$ denote the normalization of the universal family over $W$ with evaluation map $ev^{\nu}: U^{\nu} \to Z$. Let $s: C \to X$ be a general map parametrized by $W$.  Assume that
\begin{itemize}
\item $ev^{\nu}$ is dominant with connected fibers,
\item $s$ is birational onto its image, and
\item $s_*[C]$ is contained in the cone $\mathcal{C}$.
\end{itemize}
Define
\begin{equation*}
\mu = \Vert \SP_{X,s_{*}[C]}(\mathcal{E}) - \SP_{C}((s^{*}\mathcal{E})_{tf}) \Vert_{sup}.
\end{equation*}
Letting $M$ denote an irreducible component of ${\mathcal M}_{g,0}(X)$ containing the image of $W$, we have
\begin{equation*}
\codim_{M}(W) \geq \min\{  L(H \cdot C) ,  S(\mu) \}.
\end{equation*}
\end{theorem}

\begin{proof}
Let $\phi: X' \to X$ be a birational morphism from a smooth projective variety $X'$ that flattens the family parametrized by $W$ (in the sense of Construction \ref{cons:flatteningfamilyofcurves}).  Let $W' \subset W$ be an open subset parametrizing the strict transforms $s': C \to X'$ of the general maps parametrized by $W$.

The first step is to apply Grauert-M\" ulich to the family of curves $W'$.  Let $M' \subset {\mathcal M}_{g,0}(X')$ denote an irreducible component which contains $W'$ and maps to $M$ under pushforward.  Set
\begin{equation*}
\mu' = \Vert \SP_{X',s'_{*}[C]}((\phi^{*}\mathcal{E})_{tf}) - \SP_{C}(s'^{*}(\phi^{*}\mathcal{E})_{tf}) \Vert_{sup}.
\end{equation*}
By Theorem \ref{theo:tangentgaps} we have
\begin{equation*}
\codim_{M'}(W') \geq \frac{\mu'}{(\gamma+1) (\rk(\mathcal{E})-1)} - \gamma
\end{equation*}
where $\gamma = g (\dim(X)-1) + 1$.  Set
\begin{equation*}
\mu'' = \Vert \SP_{X,s_{*}[C]}(\mathcal{E}) - \SP_{X',s'_{*}[C]}((\phi^{*}\mathcal{E})_{tf}) \Vert_{sup}.
\end{equation*}
Theorem \ref{theo:kxbound} shows that under the conditions of the theorem, there are positive constants $\rho, \eta$ such that
\begin{equation*}
K_{X'/X} \cdot s'_{*}C \geq \min\{ \rho(H \cdot \alpha), \eta \cdot \mu'' \}.
\end{equation*}
Combining, we see that
\begin{align*}
\codim_{M}(W) & = \codim_{M'}(W') + \dim(M) - \dim(M') \\
& \geq \left( \frac{\mu'}{(\gamma+1) (\rk(\mathcal{E})-1)} - \gamma \right) +   \dim(M) - \dim(M')
\end{align*}
Applying Lemma \ref{lemm:expdimmbar} to estimate $\dim(M)$ and $\dim(M')$, we find
\begin{align*}
\codim_{M}(W) & \geq \left( \frac{\mu'}{(\gamma+1) (\rk(\mathcal{E})-1)} - \gamma \right) + K_{X'/X} \cdot s'_{*}C - (\dim X -1)g.
\end{align*}
Note that there is a natural surjection $\phi^{*}\mathcal{E} \to (\phi^{*}\mathcal{E})_{tf}$ whose kernel $\mathcal{K}$ is torsion.  Since the curves $s': C \to X'$ form a dominant family, a general curve $s'$ will not be contained in the support of $\mathcal{K}$ and thus $(s'^{*}\phi^{*}\mathcal{E})_{tf} \cong (s'^{*}(\phi^{*}\mathcal{E})_{tf})_{tf}$.  If $s: C \to X$ denotes the map corresponding to $s'$, then it also follows that $(s^{*}\mathcal{E})_{tf} \cong (s'^{*}(\phi^{*}\mathcal{E})_{tf})_{tf}$.  Thus the triangle inequality shows that $\mu \leq \mu' + \mu''$.  We conclude
\begin{align*}
\codim_{M}(W) & \geq \min \left\{ \rho(H \cdot \alpha) ,  \min \left\{ \eta , \frac{1}{(\gamma+1) (\rk(\mathcal{E})-1)} \right\} \mu \right\} - \left( (\dim X -1)g + \gamma \right).
\end{align*}
The two possibilities for the outermost minimum yield the two linear functions $L(H \cdot C), S(\mu)$.
\end{proof}

\begin{proof}[Proof of Theorem \ref{theo:maintheorem1}]
It is the special case of Theorem \ref{theo:maincodim} where $\mathcal{E}$ is locally free.
\end{proof}

\section{Fano varieties}

In this section we combine our main results with the theory developed by \cite{LRT23} to prove stronger theorems for Fano varieties.  The main result of this section is Theorem \ref{theo:mainfanocodim}, which is a more general version of Theorem \ref{theo:intrononfree}.  Compared against Theorem \ref{theo:maincodim}, the main advantage of Theorem \ref{theo:mainfanocodim} is that it places minimal restrictions on the family of curves (at the cost of requiring $X$ to be a Fano variety).

\subsection{Fujita invariant and accumulating maps}
The Fujita invariant is one of the basic tools for understanding the geometry of curves on Fano varieties.

\begin{definition}
\label{defi:a-invariant}
Let $X$ be a smooth projective variety over a field of characteristic $0$ and let $L$ be a big and nef $\mathbb{Q}$-Cartier divisor on $X$.  The Fujita invariant of $(X,L)$ is
\begin{equation*}
a(X,L) = \min \{ t\in \mathbb{R} \mid  K_X + tL \textrm{ is pseudo-effective }\}.
\end{equation*}
If $L$ is nef but not big, we formally set $a(X,L) = \infty$.  If $X$ is singular, choose a resolution of singularities $\phi: X' \to X$ and define $a(X,L)$ to be $a(X',\phi^{*}L)$.  (The choice of resolution does not affect the value by \cite[Proposition 2.7]{HTT15}.)
\end{definition}

\begin{definition} \label{defi:accumulatingmap}
Let $X$ be a smooth projective Fano variety.  Suppose $f: Y \to X$ is a generically finite morphism that is not birational.  We say that $f$ is an accumulating map if $a(Y,-f^{*}K_{X}) \geq a(X,-K_{X})$.
\end{definition}

The following modification of \cite[Theorem 11.5]{LRT23} connects Fujita invariants to the geometry of curves.

\begin{theorem} \label{theo:nonconnectedfiberscase}
Let $X$ be a smooth projective Fano variety.  Fix a genus $g$.  There is a linear function $R(d)$ whose leading coefficient is a positive number depending only on $\dim(X)$ and $g$ such that the following property holds.

Let $W$ be a variety equipped with a generically finite morphism $W \to {\mathcal M}_{g,0}(X)$.  Let $p: U^{\nu} \to W$ denote the normalization of the universal family over $W$ with evaluation map $ev^{\nu}: U^{\nu} \to X$.  

Let $M$ denote an irreducible component of $\mathcal{M}_{g,0}(X)$ with anticanonical degree $d$ containing the image of $W$.  Then either:
\begin{enumerate}
\item the codimension of $W$ in $M$ is at least $R(d)$,
\item $ev^{\nu}$ is dominant with connected fibers, or
\item $ev^{\nu}$ factors through an accumulating morphism.
\end{enumerate}
\end{theorem}

\begin{proof}
Assume that (2) does not hold for $ev^{\nu}$.  Thus there is a generically finite morphism $f: Y \to X$ that is not birational such that $ev^{\nu}$ factors rationally through $f$.  We construct a linear function $R(d)$ so that if the codimension of $W$ in $M$ is smaller than $R(d)$ then $a(Y,-f^{*}K_{X}) \geq a(X,-K_{X})$.

Starting from $W$, we can find a fixed genus $g$ curve $B$ and an irreducible subscheme $\widetilde{W} \subset \Mor(B,X)$ such that the maps parametrized by $\widetilde{W}$ are also parametrized by $W$ and $\dim(\widetilde{W}) \geq \dim(W) - (3g-3)$.  We can equally well think of $\widetilde{W}$ as a parameter space for sections of the projection map $\pi: X \times B \to B$.  We now apply \cite[Theorem 11.5]{LRT23} to the Fano fibration $\pi: X \times B \to B$ and the sublocus $\widetilde{W} \subset \Sec(X \times B/B)$.  The conclusion is the existence of a linear function $R'(d)$ such that if the codimension of $\widetilde{W} \subset \Sec(X \times B/B)$ is smaller than $R'(d)$ then $a(Y,-f^{*}K_{X}) \geq a(X,-K_{X})$.  Since the difference between the codimension of $\widetilde{W}$ and the codimension of $W$ is bounded by a constant depending only on $g$ we obtain the desired statement.
\end{proof}

\subsection{Codimension bounds for Fano varieties}

Combining results, we obtain a comprehensive statement for Fano varieties.  Theorem \ref{theo:mainfanocodim} improves on Theorem \ref{theo:maincodim} by addressing families of curves whose evaluation map (1) fails to dominate $X$ or (2) dominates $X$ but fails to have connected fibers.

\begin{theorem} \label{theo:mainfanocodim}
Let $X$ be a smooth projective Fano variety and let $\mathcal{E}$ be a non-zero torsion-free sheaf on $X$.  Fix a closed cone $\mathcal{C}$ such that $\mathcal{C} \backslash \{0\}$ is contained in the interior of $\Nef_{1}(X)$ and $\mathcal{C}$ contains a class $\beta \in N_{1}(X)_{\mathbb{Z}}$.  There are affine linear functions $T,S: \mathbb{R} \to \mathbb{R}$ with positive leading coefficients which satisfy the following property.

Let $W$ be a variety equipped with a generically finite morphism $W \to {\mathcal M}_{g,0}(X)$ and let $p: {U}^{\nu} \to W$ denote the normalization of the universal family over $W$ with evaluation map $ev^{\nu}: {U}^{\nu} \to X$.  Assume that
\begin{itemize}
\item a general map parametrized by $W$ is birational onto its image, and
\item the class of the general map is contained in the cone $\mathcal{C}$.
\end{itemize}
Define
\begin{equation*}
\mu = \Vert \SP_{X,s_{*}[C]}(\mathcal{E}) - \SP_{C}((s^{*}\mathcal{E})_{tf}) \Vert_{sup}
\end{equation*}
where $s: C \to X$ is a general map parametrized by $W$.  Then one of the following properties holds:
\begin{enumerate}
\item Letting $M$ denote an irreducible component of ${\mathcal M}_{g,0}(X)$ containing the image of $W$, we have
\begin{equation*}
\codim_{M}(W) \geq \min\{  T(-K_{X} \cdot C) ,  S(\mu)  \}.
\end{equation*}
\item The evaluation map $ev^{\nu}$ factors rationally through an accumulating morphism $f: Y \to X$.
\end{enumerate}
\end{theorem}

\begin{proof}
If $ev^{\nu}$ is dominant with connected fibers, Theorem \ref{theo:maincodim} (applied with $H=-K_{X}$) gives linear functions $L,S$ such that
\begin{equation*}
\codim_{M}(W) \geq \min\{ L(-K_{X} \cdot C) ,  S(\mu)  \}.
\end{equation*}
If $ev^{\nu}$ is not dominant or does not have connected fibers, Theorem \ref{theo:nonconnectedfiberscase} gives a linear function $R(-K_{X} \cdot C)$ such that either (2) holds or
\begin{equation*}
\codim_{M}(W) \geq   R(-K_{X} \cdot C).
\end{equation*}
Combining we find the result.
\end{proof}

We expect that one can do a little better:

\begin{question}
Let $X$ be a smooth projective Fano variety.  Can one prove a version of Theorem \ref{theo:mainfanocodim} that holds for the entire nef cone and not just classes in the smaller cone $\mathcal{C}$?
\end{question}

Our next result requires the following definition:

\begin{definition} \label{defi:mfree}
Let $X$ be a smooth projective variety and let $C$ be a smooth projective curve.  We say that $s: C \to X$ is $m$-free if $\mu^{min}(s^{*}T_{X}) \geq 2g+m$.
\end{definition}

As an application of Theorem \ref{theo:mainfanocodim}, we prove Theorem \ref{theo:intrononfree} showing that the locus of non-$m$-free curves in ${\mathcal{M}}_{g,0}(X)$ has codimension which increases linearly in the degree unless it is contained in the exceptional set.  One can view this result as supporting evidence for Peyre's formulation of Manin's Conjecture in \cite{Peyre17}.

\begin{proof}[Proof of Theorem \ref{theo:intrononfree}:]
By enlarging $\mathcal{C}$ slightly we may assume that it also contains a class $\beta \in N_{1}(X)_{\mathbb{Z}}$.

We first recall a few facts about the tangent bundle of $X$.  By \cite[Lemma 3.3.3]{Neumann10} there is a finite set of filtrations of $T_{X}$ such that for every $\alpha \in \Nef_{1}(X)$ the $\alpha$-Harder-Narasimhan filtration is a member of our set.   Furthermore \cite[Theorem 1.4]{Ou23} shows that every quotient of $T_{X}$ has non-zero pseudo-effective first Chern class.  Together, these imply that $\mu^{min}_{\alpha}(T_{X}): \Nef_{1}(X) \to \mathbb{R}$ is a non-negative piecewise linear function on $\Nef_{1}(X)$ that can only vanish along the boundary of the cone.  Since $\mathcal{C}$ is interior to $\Nef_{1}(X)$, we conclude that $\mu^{min}_{\alpha}(T_{X}): \mathcal{C} \to \mathbb{R}$ is bounded below by $\zeta (-K_{X} \cdot \alpha)$ for some positive constant $\zeta$.

We apply Theorem \ref{theo:mainfanocodim} to $T_{X},\mathcal{C}$ to obtain linear functions $T',S$.  Given a family $W$ as in the statement, define
\begin{equation*}
\mu_{C} = \zeta (-K_{X} \cdot s_{*}C) - (m+2g).
\end{equation*}
Since the general curve $s: C \to X$ parametrized by $W$ fails to be $m$-free, we have
\begin{equation*}
\Vert \SP_{X,s_{*}[C]}(T_{X}) - \SP_{C}(s^{*}T_{X}) \Vert_{sup} \geq \mu^{min}_{s_{*}[C]}(T_{X}) - (m+2g) \geq \mu_{C}.
\end{equation*}

First suppose that the general map parametrized by $W$ is birational onto its image.  Then Theorem \ref{theo:mainfanocodim} shows that either the evaluation map for $W$ factors rationally through an accumulating morphism or 
\begin{equation*}
\codim_{M}(W) \geq \min\{ T'(-K_{X} \cdot C),  S(\mu_{C}) \}.
\end{equation*}
(Since $S$ has positive leading coefficient, the inequality $S(\Vert \SP_{X,s_{*}[C]}(T_{X}) - \SP_{C}(s^{*}T_{X}) \Vert_{sup}) \geq S(\mu_{C})$ holds even if $\mu_{C}$ is negative.)
Since $\mu_{C}$ depends linearly on $-K_{X} \cdot s_{*}C$, we can finish the proof by choosing a linear function $T$ with positive leading coefficient which is bounded above by the two linear functions of $(-K_{X} \cdot C)$ in the previous equation.

Next suppose that the general curve $s: C \to X$ parametrized by $W$ factors as a finite morphism $C \to C'$ of degree $f \geq 2$ to a genus $g'$ curve followed by a morphism $s': C' \to X$ that is birational onto its image.  We write $e$ for $-K_{X} \cdot s'_{*}C'$.  The dimension of the space $\mathcal{M}_{g,0}(C',e)$ depends on the genus $g'$:
\begin{itemize}
\item If $g'=0$ then the dimension is $2f + 2g - 2$.
\item If $g'=1$ then the dimension is $2g-2$.
\item If $g'\geq 2$ then the dimension is $0$.
\end{itemize}
Let $Y \subset X$ denote the subvariety swept out by the images of the curves $s'$.  If $a(Y,-K_{X}|_{Y}) \geq 1$ then the inclusion $Y \to X$ is an accumulating map and we are in case (2) of the statement of the theorem.  Otherwise, we see that $K_{Y} - K_{X}|_{Y}$ has non-negative intersection against the family of curves $s'(C')$.  In particular the dimension of the space $\mathcal{M}_{g',0}(X)$ is bounded above by
\begin{equation*}
-K_{Y} \cdot s'_{*}C' + \dim(Y) + 2g' - 3 \leq -K_{X} \cdot s'_{*}C' + \dim(X) + 2g' - 3.  
\end{equation*}
Combining, we see that
\begin{align*}
\dim(W) & \leq (2f + 2g - 2) + (e + \dim(X) + 2g' - 3) \\
& = \left( \frac{2}{e} + \frac{1}{f} \right) (-K_{X} \cdot s_{*}C) + (\dim(X) + 2g + 2g' - 5).
\end{align*}
By assumption $f \geq 2$.  If the coefficient of $(-K_{X} \cdot s_{*}C)$ is $\geq 1$, then we see that we must have $e \leq 2$ or $f=2,e \leq 4$ or $f=3,e=3$.  Note the latter two cases only occur in low degree and thus can be accounted for by appropriately modifying the function $T$ in case (1) of the theorem statement.  If $e \leq 2$, then $C'$ is a rational curve of anticanonical degree at most $2$ as in case (3).  If the coefficient of $(-K_{X} \cdot s_{*}C)$ is $< 1$, then it is at most $11/12$.  By comparing against Lemma \ref{lemm:expdimmbar} we obtain a linear bound on the codimension of $W$ in terms of the anticanonical degree as in case (1).
\end{proof}

\section{Rank two bundles on the projective plane} \label{sect:r2p2}

Since the statements in the earlier portion of the paper are designed for maximal generality, it is not reasonable to expect them to give sharp bounds. However, the techniques and perspective can give much better bounds in specific examples.  In this section we illustrate this claim by focusing on what our results say for rank two vector bundles on $\mathbb{P}^2$. These have been studied in many cases in the literature (including \cite{Barth77, Hulek79} and many follow-up papers) and also are the original setting for the Grauert-M\" ulich theorem \cite{GM75}.

\begin{example} \label{exam:ramella}
The results of \cite{Ramella} imply that the tangent bundle of $\mathbb{P}^{2}$ behaves as nicely as possible when restricted to rational curves.   \cite{Ramella} shows that:
\begin{enumerate}
\item For the general degree $d$ rational curve $s: \mathbb{P}^{1} \to \mathbb{P}^{2}$, the pullback $s^{*}T_{\mathbb{P}^{2}}$ is as balanced as possible: $s^{*}T_{\mathbb{P}^{2}} = \mathcal{O}(\lceil \frac{3d}{2} \rceil) \oplus \mathcal{O}(\lfloor \frac{3d}{2} \rfloor)$.
\item For any constant $0 < c < \lfloor \frac{d}{2} \rfloor$, the sublocus of $\Mor(\mathbb{P}^{1},\mathbb{P}^{2})_{\deg d}$ parametrizing morphisms $s$ such that
\begin{equation*}
s^{*}T_{\mathbb{P}^{2}} = \mathcal{O}\left( \left\lceil \frac{3d}{2} \right\rceil + c \right) \oplus \mathcal{O}\left( \left\lfloor \frac{3d}{2} \right\rfloor - c\right)
\end{equation*}
is smooth and irreducible and has the expected codimension (which is either $2c$ or $2c-1$).
\end{enumerate}
\end{example}

For an arbitrary stable rank $2$ vector bundle $\mathcal{E}$ on $\mathbb{P}^{2}$, one cannot expect to get the nicest possible behavior as in the previous example.  However, we are able to precisely control the obstructions.  As discussed earlier in the paper, if $W$ parametrizes a dominant family of curves $s: C \to \mathbb{P}^{2}$, then either:
\begin{enumerate}
\item the codimension of $W$ is controlled by the behavior of $\mathcal{E}|_{C}$ via the Grauert-M\" ulich theorem (as in Example \ref{exam:ramella}),
\item the pullback of $\mathcal{E}$ to a birational model $\phi: X' \to X$ that flattens the family of curves parametrized by $W$ becomes unstable with respect to the strict transform of the curves (as in Example \ref{exam:conics}), or
\item the evaluation map over $W$ factors through a non-trivial finite morphism $f: Y \to X$ (as in Example \ref{exam:lines}).
\end{enumerate}
The following theorem identifies a mild set of assumptions which allow us to precisely control the codimension in these three cases.

\begin{theorem} \label{theo:goodboundsforp2}
Let $\mathcal{E}$ be a stable rank $2$ bundle on $\mathbb{P}^{2}$.

Let $W$ be a variety equipped with a generically finite morphism $W \to \mathcal{M}_{g,0}(\mathbb{P}^{2})$ parametrizing a dominant family of curves on $\mathbb{P}^{2}$ of degree $d$.  Suppose that:
\begin{itemize}
\item the general map $s: C \to \mathbb{P}^{2}$ parametrized by $W$ is a birational immersion, and
\item $\dim(W) \geq kg$ for some $k \geq 2$.
\end{itemize}
Then one of the following conditions holds:
\begin{enumerate}
\item Setting $\mu = \Vert \SP_{\mathbb{P}^{2},s_{*}[C]}(\mathcal{E}) - \SP_{C}(s^{*}\mathcal{E}) \Vert_{sup}$, we have $\codim(W) \geq \frac{2(k-1)}{k} \mu - 1$. 
\item There is a birational morphism $\phi: X' \to \mathbb{P}^{2}$ such that $\phi^{*}\mathcal{E}$ fails to be semistable with respect to the strict transform of the curves parametrized by $W$.  In this case there is a explicit constant $\zeta'$ depending only on the Chern classes of $\mathcal{E}$ (defined in Remark \ref{rema:chernonly}) such that $\codim(W) \geq \zeta'd$.
\item There is a non-trivial generically finite morphism $f: Y \to X$ such that the evaluation map for the normalization of the universal family over $W$ factors rationally through $f$.  Then we have $\codim(W) \geq d - g$.
\end{enumerate}
\end{theorem}

\begin{remark}
Note that when the codimension of $W$ is small compared to the degree of the curves then we must be in Case (1) and then we obtain the good bounds coming from the Grauert-M\" ulich theorem.
\end{remark}

\begin{remark}
Our techniques apply even when the two bulleted assumptions in Theorem \ref{theo:goodboundsforp2} are not satisfied.  However, the constants become somewhat worse.  One can find the most general versions in Remark \ref{rema:generalversion} and Remark \ref{rema:worseboundgenfinite}.
\end{remark}

Case (1) of Theorem \ref{theo:goodboundsforp2} follows from the results of Section \ref{sect:gmcodim}.  Case (2) will be addressed in Section \ref{sect:p2domconn} and Case (3) in Section \ref{sect:p2domnonconn}.  Finally, we discuss non-stable rank $2$ bundles at the end of Section \ref{sect:r2p2}.  
As discussed in Theorem \ref{theo:introp2}, for rational curves we obtain the best possible bounds using Theorem \ref{theo:goodboundsforp2}.

\begin{proof}[Proof of Theorem \ref{theo:introp2}:]
If the curves parametrized by $W$ are rational, then $\dim(W) \geq kg$ for every positive integer $k$.  Thus in Case (1) we conclude that the codimension is least the expected value $\sup\{ 2\mu - 1, 0\}$.   On the other hand, in the moduli stack of rank 2 vector bundles on $\mathbb{P}^1$ the locus parametrizing a given splitting type has exactly the expected codimension (essentially by definition).  Thus the preimage of this locus under the classifying map from $W$ will always have at most the expected codimension so long as it is non-empty.  Altogether we see that in Case (1) the codimension is exactly the expected value.

In Case (2) and Case (3), by setting $\zeta = \inf\{\zeta', 1\}$ we obtain the desired statement.
\end{proof}

\begin{example} \label{exam:lines}
In the setting of Theorem \ref{theo:goodboundsforp2} there need not be a linear relationship between $\mu$ and $\codim(W)$ when the family of curves factors through a generically finite map $f: Y \to \mathbb{P}^{2}$.  We briefly recall the examples of \cite{Schwarzenberger61} which exhibit this phenomenon.

Let $Y$ be a smooth quadric in $\mathbb{P}^{3}$ and let $f: Y \to \mathbb{P}^{2}$ be the projection map.  The branch divisor $B$ is a smooth conic in $\mathbb{P}^{2}$.  We let $W$ denote the set of tangent lines for $B$; note that the curves parametrized by $W$ are the $f$-images of the lines in $Y$.  Define the rank two bundle $\mathcal{E}_{p,q} = f_{*}\mathcal{O}(p,q)$ where $p \geq q+2$.  \cite[Proposition 2]{Schwarzenberger61} gives the exact sequence
\begin{equation*}
0 \to O(q-1)^{\oplus p-q-1} \to \mathcal{O}(q)^{\oplus p-q+1} \to \mathcal{E}_{p,q} \to 0
\end{equation*}
and in particular $\mathcal{E}_{p,q}(-q-1)$ has no global sections. 
We conclude that every subsheaf $\mathcal{O}(k) \to \mathcal{E}_{p,q}$ must have $k \leq q$ so that $\mathcal{E}_{p,q}$ is stable. 
\cite[Proposition 8]{Schwarzenberger61} shows that the jumping lines for $\mathcal{E}_{p,q}$ are parametrized by $W$ and that for any line $\ell$ parametrized by $W$ we have $\mathcal{E}_{p,q}|_{\ell} \cong \mathcal{O}(p-1) \oplus \mathcal{O}(q)$. In particular $W$ has codimension $1$ but the difference $\mu$ between the expected and actual slope panels can be arbitrarily large. 
\end{example}

\begin{example} \label{exam:conics}
Suppose $\mathcal{E}$ is a stable rank $2$ bundle on $\mathbb{P}^{2}$.  The papers \cite{Stromme84, Manaresi90, Vitter04, Marangone24} study the behavior of jumping conics for $\mathcal{E}$.  We briefly revisit this setting using our techniques.

Suppose that $W \subset \mathbb{P}H^{0}(\mathbb{P}^{2},\mathcal{O}(2))$ parametrizes a codimension $1$ family of conics such that the general conic $C$ is smooth.  We expect $\mathcal{E}|_{C}$ to be isomorphic either to $\mathcal{O}(a) \oplus \mathcal{O}(a)$ or $\mathcal{O}(a-1) \oplus \mathcal{O}(a+1)$.  Theorem \ref{theo:goodboundsforp2} shows that if the restriction of $\mathcal{E}$ to the general curve parametrized by a codimension $1$ locus $W$ does not have the expected behavior then there must be a birational morphism $\phi: X' \to \mathbb{P}^{2}$ such that $\phi^{*}\mathcal{E}$ fails to be stable with respect to the strict transforms of the conics $C$.  In particular $W$ must parametrize the conics through a fixed point.

This phenomenon is illustrated by the following family of examples considered in Section 7 of \cite{Marangone24}. Let $\mathcal{E}$ be the kernel of the map $\OO(-d)^2 \oplus \OO(-2d+1) \to \OO$ with coordinates $(x^d, y^d, z^{2d-1})$. We have the following diagram, where $Z$ is the vanishing locus of $(x^d, y^d)$.

\begin{equation*}
\xymatrix{
 & 0 & 0 & 0 & \\
0 \ar[r] & Q \ar[r] \ar[u] & \OO(-2d+1) \ar[r] \ar[u] & \OO_Z \ar[r] \ar[u] & 0 \\
0 \ar[r] & \mathcal{E} \ar[r] \ar[u] & \OO(-d)^2 \oplus \OO(-2d+1) \ar[r] \ar[u] & \OO \ar[r] \ar[u] & 0 \\
0 \ar[r] & S \ar[r] \ar[u] & \OO(-d)^2 \ar[r] \ar[u] & I_Z \ar[r] \ar[u]  & 0 \\
 & 0 \ar[u] & 0 \ar[u] & 0 \ar[u] &
}
\end{equation*}  

Since the slope of $\mathcal{E}$ is $-2d + \frac{1}{2}$, we see that $\mathcal{E}$ is stable.  There is a codimension 1 family of conics passing through the point $p = [0:0:1] = \Supp Z$, and for a general such conic $C$ we have $Q|_C = \OO(-2d+1-d) \oplus T$, where $T$ is torsion supported at $p$ of length $d$. Thus, for large $d$ we see that $\mathcal{E}|_C$ can become very unstable when restricted to a codimension 1 family of conics.

Returning to the general situation, suppose that $W \subset \mathbb{P}H^{0}(\mathbb{P}^{2},\mathcal{O}(2))$ has codimension 2.  If $\mathcal{E}|_{C}$ fails to have the expected behavior, Theorem \ref{theo:goodboundsforp2} identifies an additional possibility: there is a generically finite morphism $f: Y \to \mathbb{P}^{2}$ such that the conics are the images of curves on $Y$.  A quick argument based on Lemma \ref{lemm:fujinvsurf} shows that in this case our family must be obtained from the conics on a quadric hypersurface $Y$ via the projection map $f: Y \to \mathbb{P}^{2}$. 
\end{example}

We will frequently use the following lemma concerning the curves in Theorem \ref{theo:goodboundsforp2}.

\begin{lemma} \label{lemm:curvesaregood}
Let $X$ be a surface and let $W$ be a variety equipped with a generically finite morphism $W \to \mathcal{M}_{g,0}(X)$ parametrizing a dominant family of curves on $X$.  Suppose that:
\begin{itemize}
\item the general map $s: C \to X$ parametrized by $W$ is a birational immersion, and
\item $\dim(W) \geq kg$ for some $k \geq 2$.
\end{itemize}
Then the normal sheaf $N_{s}$ of the general map $s: C \to X$ is a line bundle of degree $\geq kg$.  In particular $s$ is contained in the smooth locus of $\mathcal{M}_{g,0}(X)$.
\end{lemma}

\begin{proof}
Since $s$ is an immersion $N_{s}$ is a line bundle.  We have $H^{0}(C,N_{s}) \geq \dim(W) \geq kg$.  By Riemann-Roch and Serre vanishing we conclude that $\deg(N_{s}) \geq (k+1)g-1$ if $g\geq 1$ and $\deg(N_{s}) \geq 0$ if $g = 0$.
\end{proof}

\subsection{Dominant families with connected fibers} \label{sect:p2domconn}
To control the codimension of $W$ in this case we must understand the birational properties of stability.  An advantage of working with surfaces is that one can use convexity techniques instead of the log canonical threshold -- this will greatly improve the bounds we obtain.  The following theorem is an analogue of Theorem \ref{theo:kxbound} for stable rank $2$ bundles on $\mathbb{P}^{2}$ with explicit constants.  Recall that for a rank $2$ vector bundle $\mathcal{E}$ on $\mathbb{P}^{2}$ the discriminant is $\Delta(\mathcal{E}) = 4c_{2}(\mathcal{E}) - c_{1}(\mathcal{E})^{2}$.

\begin{theorem} \label{theo:p2bound}
Let $\mathcal{E}$ be a stable rank $2$ vector bundle on $\mathbb{P}^{2}$.  Let $\phi: X' \to \mathbb{P}^{2}$ be a birational map and let $\alpha' \in \Nef_{1}(X')$ be such that $\phi_{*}\alpha'$ is the hyperplane class $[H]$.  Let $\mathcal{E} \to \mathcal{Q}$ be a torsion-free quotient of rank $1$ and let $\mathcal{Q}'$ denote the birational transform of $\mathcal{Q}$ on $X'$.  Suppose that $\mu_{\alpha'}(\mathcal{Q}') \leq \mu_{\alpha'}(\phi^{*}\mathcal{E})$.  Then
\begin{equation*}
K_{X'/\mathbb{P}^{2}} \cdot \alpha' \geq  \sqrt{ 1 - \frac{  \Delta(\mathcal{E}) }{4\left( c_{1}(\mathcal{Q}) \cdot H - \mu(\mathcal{E}) \right)^{2} + \Delta(\mathcal{E}) } }.
\end{equation*}
\end{theorem}

Note that $\Delta(\mathcal{E}) \geq 0$ by Bogomolov's inequality.  Furthermore $c_{1}(\mathcal{Q}) \cdot H > \mu(\mathcal{E})$ by stability.  Thus the right hand side is between $0$ and $1$ and is minimized when $c_{1}(\mathcal{Q})$ is as small as possible.

\begin{proof}
For convenience we set $c_{1}(\mathcal{E}) = e[H]$ and $c_{2}(\mathcal{E}) = f[H^{2}]$.  We also write $c_{1}(\mathcal{Q}) = d[H]$; note that $d$ is an integer satisfying $d > \mu(\mathcal{E})$.

We may assume that $\mathcal{Q}'$ is locally free.  Indeed, if this fails to be the case, then we can precompose $\phi$ with a birational morphism $\psi: \widetilde{X} \to X'$ that resolves $\mathcal{Q}'$.  Set $\widetilde{\alpha} = \psi^{*}\alpha'$.  Then the hypotheses of the theorem still hold for $(\widetilde{X},\widetilde{\alpha})$ and the desired statement for $(X',\alpha')$ follows from the analogous statement for $(\widetilde{X},\widetilde{\alpha})$.

Let $\mathcal{F}$ denote the kernel of $\mathcal{E} \to \mathcal{Q}$ and let $\mathcal{F}'$ denote the kernel of $\phi^{*}\mathcal{E} \to \mathcal{Q}'$.  Since $\mathcal{F},\mathcal{F}'$ are saturated subsheaves of a reflexive sheaf they are also reflexive, hence locally free.

We can realize $\phi$ as a sequence of point blow-ups.  We let $E_{i}$ denote the pullback of the $i$th exceptional divisor to $X'$.  (N.B.:~$E_{i}$ is usually not the strict transform of the $i$th exceptional divisor on $X'$.)  We choose (necessarily non-negative) integers $a_{i},b_{i}$ so that $c_{1}(\mathcal{Q}') = \phi^{*}c_{1}(\mathcal{Q}) - \sum b_{i}E_{i}$ and $\alpha' = \phi^{*}[H] - \sum a_{i}E_{i}$. Note that $K_{X'/\mathbb{P}^{2}} = \sum E_{i}$.  By the Cauchy-Schwarz inequality we have:
\begin{align*}
K_{X'/\mathbb{P}^{2}} \cdot \alpha' & = \sum a_{i} \\
& \geq \left( \sum a_{i}^{2} \right)^{1/2} \\
& \geq \frac{\sum a_{i}b_{i}}{ \left( \sum b_{i}^{2} \right)^{1/2} } 
\end{align*}
We will bound $\sum a_{i}b_{i}$ from below and compute $\sum b_{i}^{2}$ precisely.  The first bound follows from our condition on slopes:
\begin{align*}
\frac{c_{1}(\mathcal{E}) \cdot [H]}{2} & \geq c_{1}(\mathcal{Q}') \cdot \alpha' \\
& = c_{1}(\mathcal{Q}) \cdot  [H] - \sum a_{i}b_{i}
\end{align*}
or equivalently
\begin{equation*}
\sum a_{i}b_{i} \geq c_{1}(\mathcal{Q}) \cdot  [H] - \frac{c_{1}(\mathcal{E}) \cdot  [H]}{2} = d - \frac{e}{2}.
\end{equation*}
The second computation follows from a Chern class argument.  Note that
\begin{align*}
c_{1}(\mathcal{F})c_{1}(\mathcal{Q}) + c_{2}(\mathcal{Q}) = c_{2}(\mathcal{E}) \\
c_{1}(\mathcal{F}')c_{1}(\mathcal{Q}') = \phi^{*}c_{2}(\mathcal{E})
\end{align*}
Also we know that $c_{1}(\mathcal{Q}) + c_{1}(\mathcal{F}) = c_{1}(\mathcal{E})$ and similarly for $\mathcal{F}',\mathcal{Q}'$.  Thus if we
subtract the pullback of the first line from the second, we obtain
\begin{equation*}
\sum b_{i}^{2} = \phi^{*}c_{2}(\mathcal{Q}) = \phi^{*}c_{2}(\mathcal{E}) - \phi^{*}(c_{1}(\mathcal{E})- c_{1}(\mathcal{Q})) \cdot \phi^{*}c_{1}(\mathcal{Q}).
\end{equation*}
Combining with earlier equations, we see that
\begin{equation*}
K_{X'/\mathbb{P}^{2}} \cdot \alpha' \geq  \frac{ d - \frac{e}{2} }{\left( d^{2} - de + f \right)^{1/2} }.
\end{equation*}
The desired inequality is obtained by rearranging the right hand side.
\end{proof}

\begin{example} \label{exam:p2birstab}
Suppose $\mathcal{E} = T_{\mathbb{P}^{2}}$.  Keeping the notation used in Theorem \ref{theo:p2bound}, the result shows that
\begin{align*}
K_{X'/\mathbb{P}^{2}} \cdot \alpha' & \geq \sqrt{ 1 - \frac{ 3}{4 \left( d^{2} - 3d + 3 \right)}}
\end{align*}
Since $T_{\mathbb{P}^{2}}$ is stable we must have $d \geq 2$.  Thus the expression on the right is minimized when $d=2$, yielding
\begin{equation*}
K_{X'/\mathbb{P}^{2}} \cdot \alpha \geq \frac{1}{2}.
\end{equation*}
Furthermore, the argument shows that equality is achieved precisely when the following conditions are met:
\begin{itemize}
\item The equality $\sum a_{i} = \sqrt{ \sum a_{i}^{2}}$ shows that there is a single non-zero $a_{i}$.
\item The equality $b_{i} = (\sum b_{i}^{2})^{1/2}$ shows that there is a single non-zero $b_{i}$.
\item We must have $d=2$.
\end{itemize}
Thus (after ignoring unnecessary blow-ups) $\phi$ is the blow-up of a point and the quotient has the form
\begin{equation*}
0 \to \mathcal{O}(1) \to T_{\mathbb{P}^{2}} \to \mathcal{I}(2) \to 0
\end{equation*}
where $\mathcal{I}$ is an ideal sheaf concentrated at the point.  By calculating $c_{2}$ on the right, we see that in fact $\mathcal{I}$ must be the ideal sheaf of the point.

Conversely, suppose $\phi: X' \to \mathbb{P}^{2}$ is the blow-up of a point and $\alpha' = \phi^{*}H - tE$ is a nef class for some $0 \leq t \leq 1$.  The computation above shows that $\phi^{*}T_{\mathbb{P}^{2}}$ is semistable for $0 \leq t \leq \frac{1}{2}$.  For $\frac{1}{2} < t \leq 1$ we have the destabilizing sequence
\begin{equation*}
0 \to \mathcal{O}_{X'}(H+E) \to \phi^{*}T_{\mathbb{P}^{2}} \to \mathcal{O}_{X'}(2H - E) \to 0.
\end{equation*}
where the leftmost term can be interpreted as the relative tangent bundle $T_{X'/\mathbb{P}^{1}}$ for the $\mathbb{P}^{1}$-bundle structure on $X'$ and the leftmost map is the inclusion $T_{X'/\mathbb{P}^{1}} \to T_{X'} \to \phi^{*}T_{\mathbb{P}^{2}}$.  When we pushforward this exact sequence to $\mathbb{P}^{2}$ we obtain the sequence identified in the previous paragraph.  Note that we find destabilizing classes $\alpha'$ such that $K_{X'/\mathbb{P}^{2}} \cdot \alpha'$ is arbitrarily close to $1/2$.
\end{example}

Combining with dimension estimates, we get a codimension bound.  We will focus on the version with more assumptions and nicer bounds; the general version is stated in Remark \ref{rema:generalversion}. 

\begin{theorem} \label{theo:cleanp2connectedfibers}
Let $\mathcal{E}$ be a stable rank $2$ vector bundle on $\mathbb{P}^{2}$.
Suppose that $W \to {\mathcal{M}}_{g,0}(\mathbb{P}^{2})$ is a generically finite morphism and that the normalization of the universal family defines a dominant family of curves on $\mathbb{P}^{2}$ of degree $d$ such that the evaluation map has connected fibers.   Furthermore suppose that:
\begin{itemize}
\item the general curve $s: C \to \mathbb{P}^{2}$ parametrized by $W$ is a birational immersion, and
\item $\dim(W) \geq kg$ for some $k \geq 2$.
\end{itemize}
Define
\begin{equation*}
\mu = \Vert \SP_{\mathbb{P}^{2},s_{*}[C]}(\mathcal{E}) - \SP_{C}(s^{*}\mathcal{E}) \Vert_{sup}.
\end{equation*}
where $s: C \to \mathbb{P}^{2}$ is a general morphism parametrized by $W$.  Then
\begin{equation*}
\codim_{M}(W) \geq \min \left\{ \frac{2(k-1)\mu}{k} - 1, d\zeta \right\}
\end{equation*}
where $\zeta$ is defined as an infimum over all quotients $\mathcal{E} \twoheadrightarrow \mathcal{Q}$ of the constant from Theorem \ref{theo:p2bound}:
\begin{equation*}
\zeta = \inf_{\mathcal{E} \twoheadrightarrow \mathcal{Q}}   \sqrt{ 1 - \frac{  \Delta(\mathcal{E}) }{4\left( c_{1}(\mathcal{Q}) \cdot H - \mu(\mathcal{E}) \right)^{2} + \Delta(\mathcal{E}) } }.
\end{equation*}
\end{theorem}

\begin{proof}
Let $\phi: X' \to \mathbb{P}^{2}$ be a birational map flattening the family of curves.  Let $W'$ denote the parameter space of the strict transforms $C'$ on $X'$.  Note that the general $s': C' \to X'$ parametrized by $W'$ is still a birational immersion that moves in dimension $\geq kg$. By Lemma \ref{lemm:curvesaregood} $W'$ is contained in the smooth locus of $\mathcal{M}_{g,0}(X')$.  Thus the irreducible component $M'$ of ${\mathcal{M}}_{g,0}(X')$ containing $W'$ must have the expected dimension.

We split into two cases.  In the first case, we have that $\phi^{*}\mathcal{E}$ is $[C']$-semistable.   By Lemma \ref{lemm:curvesaregood} $\mu^{min}(N_{s'}) \geq kg$ and so by Corollary \ref{coro:cleantangentgaps} we have
\begin{equation*}
\codim_{M}(W) \geq \codim_{M'}(W') \geq \frac{2(k-1)}{k} \cdot \mu - 1.
\end{equation*}

In the second case, we have that $\phi^{*}\mathcal{E}$ is no longer $[C']$-semistable.  Thus there is some quotient $\mathcal{E} \to \mathcal{Q}$ whose birational transform destabilizes $\phi^{*}\mathcal{E}$.  If we let $M$ denote the irreducible component of ${\mathcal{M}}_{g,0}(\mathbb{P}^{2})$ containing the image of $W$,
then
\begin{align*}
\dim(W') & \leq -K_{X'} \cdot C' + g-1 \\
& = -K_{X'/\mathbb{P}^{2}} \cdot C' + (-K_{\mathbb{P}^{2}} \cdot C) + g-1 \\
& = \dim(M) - K_{X'/\mathbb{P}^{2}} \cdot C'.
\end{align*} 
so that $\codim_{M}(W) \geq (K_{X'/\mathbb{P}^{2}} \cdot C')$.  Applying Theorem \ref{theo:p2bound}, we conclude
\begin{equation*}
\codim_{M}(W) \geq d \zeta.
\end{equation*}
\end{proof}

\begin{remark} \label{rema:generalversion}
By a similar argument (which we omit), one can prove a more general statement: we can remove the two itemized restrictions on $W$ in Theorem \ref{theo:cleanp2connectedfibers} at the cost of weakening the inequality to 
\begin{equation*}
\codim_{M}(W) \geq \min \left\{ \frac{2\mu}{2g+3} - (g+1), d\zeta - g \right\}.
\end{equation*}
\end{remark}

\begin{remark} \label{rema:chernonly}
We can formulate a version of Theorem \ref{theo:cleanp2connectedfibers} that only relies on the Chern classes of $\mathcal{E}$ by defining
\begin{equation*}
\zeta' = \inf_{d \in \mathbb{Z}, d > \mu(\mathcal{E})}   \sqrt{ 1 - \frac{  \Delta(\mathcal{E}) }{4\left( d - \mu(\mathcal{E}) \right)^{2} + \Delta(\mathcal{E}) } }.
\end{equation*}
It is clear that $\zeta'$ only depends on the Chern classes of $\mathcal{E}$ and that $0 < \zeta' \leq \zeta$.
\end{remark}

\begin{example}
Consider again the tangent bundle of $\mathbb{P}^{2}$ and suppose we have a family $W$ of curves $s: C \to \mathbb{P}^{2}$ whose evaluation map is dominant with connected fibers. Suppose furthermore that $T_{\mathbb{P}^{2}}$ fails to be semistable with respect to the strict transform family on the flattening birational model $\phi: X' \to \mathbb{P}^{2}$.  Combining Remark \ref{rema:generalversion} with Example \ref{exam:p2birstab}, we see that $\frac{\dim(W)}{\dim(M)}$ is bounded above by a number that approaches $5/6$ as the degree gets large. 
\end{example}

\subsection{Dominant families with disconnected fibers} \label{sect:p2domnonconn}
We next turn to families of curves which factor through a generically finite morphism.  According to Theorem \ref{theo:nonconnectedfiberscase}, we can bound the codimension of $W$ using the Fujita invariant of covers $f: Y \to X$.  Using some basic birational geometry, we have:

\begin{lemma} \label{lemm:fujinvsurf}
Let $f: Y \to \mathbb{P}^{2}$ be a dominant generically finite morphism of degree $\geq 2$ from a normal surface $Y$.  Then $a(Y,-f^{*}K_{\mathbb{P}^{2}}) \in [0,\frac{1}{3}] \cup \{\frac{1}{2}, \frac{2}{3} \}$.  Furthermore:
\begin{enumerate}
\item If $a(Y,-f^{*}K_{\mathbb{P}^{2}}) = \frac{2}{3}$ then $f: Y \to \mathbb{P}^{2}$ is birationally equivalent to either:
\begin{enumerate}
\item the degree $2$ map $f': Y' \to \mathbb{P}^{2}$ obtained by projecting a (possibly singular) quadric surface, or
\item a map $f': Y' \to \mathbb{P}^{2}$ from a ruled surface $Y'$ which takes the ruling to lines in $\mathbb{P}^{2}$.
\end{enumerate}
\item If $a(Y,-f^{*}K_{\mathbb{P}^{2}})= \frac{1}{2}$ then $f$ is birationally equivalent to a degree $2$ morphism $f': \mathbb{P}^{2} \to \mathbb{P}^{2}$.
\end{enumerate}
\end{lemma}

\begin{proof}
Note that $a(Y,-f^{*}K_{\mathbb{P}^{2}}) = \frac{1}{3}a(Y,f^{*}H)$ where $H$ is the hyperplane class on $\mathbb{P}^{2}$.  By \cite[Proposition 1.3]{Horing10} the possible values of $a(Y,f^{*}H)$ in the range $[1,\infty)$ are $1,\frac{3}{2}$, $2$, and $3$ and the final case could only occur if $f$ were birational.  The explicit description of $Y$ and $f$ is implied by \cite[Lemma 5.3]{LTJAG} and \cite[Proposition 1.3]{Horing10}.
\end{proof}

The following lemma clarifies the connection between Fujita invariants and dimension in this case.

\begin{lemma} \label{lemm:disconnectedfibersp2}
Suppose that $W \to {\mathcal{M}}_{g,0}(\mathbb{P}^{2},d)$ is a generically finite morphism such that the normalization of the evaluation map is dominant and factors rationally through a non-birational generically finite morphism $f: Y \to \mathbb{P}^{2}$ from a smooth projective variety $Y$.  Furthermore assume that the general map parametrized by $W$ is a birational immersion and that $\dim(W) \geq kg$ for some $k \geq 2$. 

Let $M$ be an irreducible component of ${\mathcal{M}}_{g,0}(\mathbb{P}^{2},d)$ containing $W$.  Then
\begin{equation*}
\dim(W) \leq a(Y,-f^{*}K_{\mathbb{P}^{2}}) \dim(M) + \sup\{ g-1,0\}.
\end{equation*}
\end{lemma}

\begin{proof}
The hypotheses of Lemma \ref{lemm:curvesaregood} hold both for $W$ and for the induced morphisms $s_{Y}: C \to Y$ corresponding to the curves parametrized by $W$.  We conclude that $\dim(W) \leq -K_{Y} \cdot s_{Y*}C +g-1$ and $\dim(M) = -K_{\mathbb{P}^{2}} \cdot s_{*}C + g-1$.   Furthermore since $K_{Y} - a(Y,-f^{*}K_{\mathbb{P}^{2}})f^{*}K_{\mathbb{P}^{2}}$ is pseudo-effective it has non-negative intersection against $s_{Y*}C$ so that $-K_{Y}\cdot s_{Y*}C \leq -a(Y,-f^{*}K_{\mathbb{P}^{2}})K_{\mathbb{P}^{2}} \cdot s_{*}C$.  Thus
\begin{align*}
\dim(W)  & \leq -a(Y,-f^{*}K_{\mathbb{P}^{2}})K_{\mathbb{P}^{2}} \cdot s_{*}C + g - 1 \\
& = a(Y,-K_{\mathbb{P}^{2}}) \dim(M) + (1-a(Y,-f^{*}K_{\mathbb{P}^{2}}))(g-1) \\
& \leq  a(Y,-K_{\mathbb{P}^{2}}) \dim(M) + \sup\{ (g-1), 0 \}.
\end{align*}
\end{proof}

\begin{remark} \label{rema:worseboundgenfinite}
If we remove the assumption $\dim(W) \geq kg$, then we obtain the weaker inequality
\begin{equation*}
\dim(W) \leq a(Y,-f^{*}K_{\mathbb{P}^{2}}) \dim(M) + \sup\{ 2g-1, 0\}.
\end{equation*}
\end{remark}

Combining Lemma \ref{lemm:fujinvsurf} and Lemma \ref{lemm:disconnectedfibersp2} we obtain:

\begin{lemma} \label{lemm:genfinitecodim}
Suppose that $W \to {\mathcal{M}}_{g,0}(\mathbb{P}^{2},d)$ is a generically finite morphism such that the normalization of the evaluation map is dominant and factors rationally through a non-birational generically finite morphism $f: Y \to \mathbb{P}^{2}$.  Furthermore assume that the general map parametrized by $W$ is a birational immersion and that $\dim(W) \geq kg$ for some $k \geq 2$.  Let $M$ denote the irreducible component of $\mathcal{M}_{g,0}(\mathbb{P}^{2})$ containing the image of $W$.  Then
\begin{equation*}
\dim(W) \leq \frac{2}{3} \dim(M) + g.
\end{equation*}
\end{lemma}

\begin{proof}[Proof of Theorem \ref{theo:goodboundsforp2}]
Follows from Theorem \ref{theo:cleanp2connectedfibers} and Lemma \ref{lemm:genfinitecodim}.
\end{proof}

\subsection{Non-stable bundles}
Finally we consider the behavior of restrictions when $\mathcal{E}$ is a non-stable rank $2$ bundle on $\mathbb{P}^{2}$.  It turns out that this case can be handled directly.  Consider the exact sequence
\begin{equation*}
0 \to \mathcal{F} \to \mathcal{E} \to \mathcal{Q} \to 0
\end{equation*}
where $\mathcal{F}$ is a maximal destabilizing subsheaf (if $\mathcal{E}$ is unstable) or a rank $1$ subsheaf with the same slope as $\mathcal{E}$ (if $\mathcal{E}$ is strictly semistable).
Since $\mathcal{Q}$ is torsion-free, $\mathcal{F}$ is a saturated subsheaf of the reflexive sheaf $\mathcal{E}$.  This implies that $\mathcal{F}$ is reflexive, hence a line bundle.

\begin{lemma} \label{lemm:unstablehnfilt}
In the situation above, for a general curve $s: C \to \mathbb{P}^{2}$ in a dominant family then either $s^{*}\mathcal{E}$ is strictly semistable or the quotient $s^{*}\mathcal{E} \to (s^{*}\mathcal{Q})_{tf}$ defines the Harder-Narasimhan filtration of $s^{*}\mathcal{E}$.
\end{lemma}

\begin{proof}
We have an exact sequence
\begin{equation*}
s^{*}\mathcal{F} \to s^{*}\mathcal{E} \to s^{*}\mathcal{Q} \to 0
\end{equation*}
Since the non-free locus of $\mathcal{Q}$ has codimension $\geq 2$, the curve $C$ meets the locus where the map $\mathcal{F} \to \mathcal{E}$ is non-zero.
Thus $s^{*}\mathcal{F} \to s^{*}\mathcal{E}$ is generically injective, hence injective.

Note that $c_{1}(s^{*}\mathcal{F}) \geq c_{1}(s^{*}\mathcal{Q})$.  Thus the kernel $\mathcal{F}'$ of the map $s^{*}\mathcal{E} \to (s^{*}\mathcal{Q})_{tf}$ also satisfies $c_{1}(\mathcal{F}') \geq c_{1}((s^{*}\mathcal{Q})_{tf})$.  If we have a strict inequality, then this quotient must define the Harder-Narasimhan filtration of $s^{*}\mathcal{E}$.  If we have equality, then $s^{*}\mathcal{E}$ is strictly semistable.  
\end{proof}

\begin{proposition}
Suppose $\mathcal{E}$ is a non-stable vector bundle on $\mathbb{P}^{2}$.  There is an affine linear function $S$ with positive leading coefficient such that the following holds.

Let $W$ be a variety equipped with a generically finite morphism $W \to {\mathcal M}_{g,0}(X)$ and let $p: U^{\nu} \to W$ denote the normalization of the universal family over $W$ with evaluation map $ev^{\nu}: U^{\nu} \to X$.  Assume that $ev^{\nu}$ is dominant.  For a general $s: C \to \mathbb{P}^{2}$ parametrized by $W$ set
\begin{equation*}
\mu = \Vert \SP_{\mathbb{P}^{2},s_{*}[C]}(\mathcal{E}) - \SP_{C}(s^{*}\mathcal{E}) \Vert_{sup}.
\end{equation*}
Then we have $\codim(W) \geq S(\mu)$.
\end{proposition}

\begin{proof}
Let $\mathcal{Q}$ be the minimal destabilizing quotient (if $\mathcal{E}$ is unstable) or a torsion-free rank $1$ quotient of the same slope as $\mathcal{E}$ (if $\mathcal{E}$ is strictly semistable).  Let $\phi: X' \to \mathbb{P}^{2}$ be a birational map that resolves the top Fitting ideal of $\mathcal{Q}$ and write $E = \phi^{*}c_{1}(\mathcal{Q}) - c_{1}((\phi^{*}\mathcal{Q})_{tf})$.  Let $\alpha'$ be the numerical class of the strict transform of the general curve parametrized by $W$.  By Lemma \ref{lemm:unstablehnfilt} we have $\mu = E \cdot \alpha'$.

We can choose a positive constant $\rho$ such that $\rho K_{X'/\mathbb{P}^{2}} \geq E$.  This readily yields the desired statement by arguing as in Theorem \ref{theo:maincodim}.
\end{proof}

\bibliographystyle{alpha}
\bibliography{tangentgaps}

\def\cprime{$'$}
\begin{thebibliography}{MnOSC12}

\bibitem[AC81]{AC81}
E.~Arbarello and M.~Cornalba.
\newblock On a conjecture of {P}etri.
\newblock {\em Comment. Math. Helv.}, 56(1):1--38, 1981.

\bibitem[AR15]{AR15}
A.~Alzati and R.~Re.
\newblock {$PGL(2)$} actions on {G}rassmannians and projective construction of
  rational curves with given restricted tangent bundle.
\newblock {\em J. Pure Appl. Algebra}, 219(5):1320--1335, 2015.

\bibitem[Asc88]{Ascenzi88}
M.-G. Ascenzi.
\newblock The restricted tangent bundle of a rational curve in {${\bf P}^2$}.
\newblock {\em Comm. Algebra}, 16(11):2193--2208, 1988.

\bibitem[Asc22]{Ascenzi22}
M.-G. Ascenzi.
\newblock The tangent bundle restricted to a rational curve spanning {$\Bbb
  P^3$}.
\newblock {\em J. Algebra}, 610:703--727, 2022.

\bibitem[Bar77]{Barth77}
W.~Barth.
\newblock Moduli of vector bundles on the projective plane.
\newblock {\em Invent. Math.}, 42:63--91, 1977.

\bibitem[BHJ17]{BHJ17}
S.~Boucksom, T.~Hisamoto, and M.~Jonsson.
\newblock Uniform {K}-stability, {D}uistermaat-{H}eckman measures and
  singularities of pairs.
\newblock {\em Ann. Inst. Fourier (Grenoble)}, 67(2):743--841, 2017.

\bibitem[Bog94]{Bogomolov95}
F.~A. Bogomolov.
\newblock Stable vector bundles on projective surfaces.
\newblock {\em Mat. Sb.}, 185(4):3--26, 1994.

\bibitem[BR97]{BR97}
E.~Ballico and B.~Russo.
\newblock On the stability of the restriction of {$T{\bf P}^n$} to projective
  curves.
\newblock In {\em Complex analysis and geometry ({T}rento, 1995)}, volume 366
  of {\em Pitman Res. Notes Math. Ser.}, pages 7--18. Longman, Harlow, 1997.

\bibitem[BR00]{BR00}
E.~Ballico and L.~Ramella.
\newblock The restricted tangent bundle of smooth curves in {G}rassmannians and
  curves in flag varieties.
\newblock {\em Rocky Mountain J. Math.}, 30(4):1207--1227, 2000.

\bibitem[BS23]{BrowningSawin}
T.~Browning and W.~Sawin.
\newblock Free rational curves on low degree hypersurfaces and the circle
  method.
\newblock {\em Algebra Number Theory}, 17(3):719--748, 2023.

\bibitem[But94]{Butler94}
D.~C. Butler.
\newblock Normal generation of vector bundles over a curve.
\newblock {\em J. Differential Geom.}, 39(1):1--34, 1994.

\bibitem[CP11]{CP11}
F.~Campana and T.~Peternell.
\newblock Geometric stability of the cotangent bundle and the universal cover
  of a projective manifold.
\newblock {\em Bull. Soc. Math. France}, 139(1):41--74, 2011.
\newblock With an appendix by Matei Toma.

\bibitem[CR18]{CR18}
I.~Coskun and E.~Riedl.
\newblock Normal bundles of rational curves in projective space.
\newblock {\em Math. Z.}, 288(3-4):803--827, 2018.

\bibitem[dJS17]{dJS17}
A.~J. de~Jong and J.~Starr.
\newblock Divisor classes and the virtual canonical bundle for genus 0 maps.
\newblock In {\em Geometry over nonclosed fields}, Simons Symp., pages 97--126.
  Springer, Cham, 2017.

\bibitem[FHS80]{FHS80}
O.~Forster, A.~Hirschowitz, and M.~Schneider.
\newblock Type de scindage g\'{e}n\'{e}ralis\'{e} pour les fibr\'{e}s stables.
\newblock In {\em Vector bundles and differential equations ({P}roc. {C}onf.,
  {N}ice, 1979)}, volume~7 of {\em Progr. Math.}, pages 65--81. Birkh\"{a}user,
  Boston, MA, 1980.

\bibitem[GHI13]{GHI13}
A.~Gimigliano, B.~Harbourne, and M.~Id\`a.
\newblock On plane rational curves and the splitting of the tangent bundle.
\newblock {\em Ann. Sc. Norm. Super. Pisa Cl. Sci. (5)}, 12(3):587--621, 2013.

\bibitem[GKP14]{GKP14}
D.~Greb, S.~Kebekus, and T.~Peternell.
\newblock Reflexive differential forms on singular spaces. {G}eometry and
  cohomology.
\newblock {\em J. Reine Angew. Math.}, 697:57--89, 2014.

\bibitem[GKP16]{GKP16}
D.~Greb, S.~Kebekus, and T.~Peternell.
\newblock Movable curves and semistable sheaves.
\newblock {\em Int. Math. Res. Not. IMRN}, (2):536--570, 2016.

\bibitem[GM75]{GM75}
H.~Grauert and G.~M\"{u}lich.
\newblock Vektorb\"{u}ndel vom {R}ang {$2$} \"{u}ber dem {$n$}-dimensionalen
  komplex-projektiven {R}aum.
\newblock {\em Manuscripta Math.}, 16(1):75--100, 1975.

\bibitem[Hei00]{Hein00}
G.~Hein.
\newblock Curves in {$\bold P^3$} with good restriction of the tangent bundle.
\newblock {\em Rocky Mountain J. Math.}, 30(1):217--235, 2000.

\bibitem[HK96]{HK93}
G.~Hein and H.~Kurke.
\newblock Restricted tangent bundle on space curves.
\newblock In {\em Proceedings of the {H}irzebruch 65 {C}onference on
  {A}lgebraic {G}eometry ({R}amat {G}an, 1993)}, volume~9 of {\em Israel Math.
  Conf. Proc.}, pages 283--294. Bar-Ilan Univ., Ramat Gan, 1996.

\bibitem[HL97]{HL97}
D.~Huybrechts and M.~Lehn.
\newblock {\em The geometry of moduli spaces of sheaves}.
\newblock Aspects of Mathematics, E31. Friedr. Vieweg \& Sohn, Braunschweig,
  1997.

\bibitem[HL20]{HL20}
Jingjun Han and Zhan Li.
\newblock On {F}ujita's conjecture for pseudo-effective thresholds.
\newblock {\em Math. Res. Lett.}, 27(2):377--396, 2020.

\bibitem[H{\"o}r10]{Horing10}
A.~H{\"o}ring.
\newblock The sectional genus of quasi-polarised varieties.
\newblock {\em Arch. Math. (Basel)}, 95(2):125--133, 2010.

\bibitem[HTT15]{HTT15}
B.~Hassett, S.~Tanimoto, and Y.~Tschinkel.
\newblock Balanced line bundles and equivariant compactifications of
  homogeneous spaces.
\newblock {\em Int. Math. Res. Not. IMRN}, (15):6375--6410, 2015.

\bibitem[Hul79]{Hulek79}
K.~Hulek.
\newblock Stable rank-{$2$} vector bundles on {${\bf P}\sb{2}$} with
  {$c\sb{1}$} odd.
\newblock {\em Math. Ann.}, 242(3):241--266, 1979.

\bibitem[Kol97]{Kollar97}
J.~Koll\'{a}r.
\newblock Singularities of pairs.
\newblock In {\em Algebraic geometry---{S}anta {C}ruz 1995}, volume~62 of {\em
  Proc. Sympos. Pure Math.}, pages 221--287. Amer. Math. Soc., Providence, RI,
  1997.

\bibitem[Kop20]{Kopper20}
John Kopper.
\newblock Stability conditions for restrictions of vector bundles on projective
  surfaces.
\newblock {\em Michigan Math. J.}, 69(4):711--732, 2020.

\bibitem[Lar16]{Larson16}
E.~Larson.
\newblock Interpolation for restricted tangent bundles of general curves.
\newblock {\em Algebra Number Theory}, 10(4):931--938, 2016.

\bibitem[Laz04a]{Lazarsfeld04a}
R.~Lazarsfeld.
\newblock {\em Positivity in algebraic geometry. {I}}, volume~48 of {\em
  Ergebnisse der Mathematik und ihrer Grenzgebiete. 3. Folge. A Series of
  Modern Surveys in Mathematics [Results in Mathematics and Related Areas. 3rd
  Series. A Series of Modern Surveys in Mathematics]}.
\newblock Springer-Verlag, Berlin, 2004.
\newblock Classical setting: line bundles and linear series.

\bibitem[Laz04b]{Lazarsfeld04b}
R.~Lazarsfeld.
\newblock {\em Positivity in algebraic geometry. {II}}, volume~49 of {\em
  Ergebnisse der Mathematik und ihrer Grenzgebiete. 3. Folge. A Series of
  Modern Surveys in Mathematics}.
\newblock Springer-Verlag, Berlin, 2004.
\newblock Positivity for vector bundles, and multiplier ideals.

\bibitem[LRT23a]{LRT24}
B.~Lehmann, E.~Riedl, and S.~Tanimoto.
\newblock Non-free curves on {F}ano varieties.
\newblock to appear in Osaka J.~Math, 2023.

\bibitem[LRT23b]{LRT23}
B.~Lehmann, E.~Riedl, and S.~Tanimoto.
\newblock Non-free sections of {F}ano fibrations.
\newblock submitted, arXiv:2301.01695, 2023.

\bibitem[LT21]{LTJAG}
B.~Lehmann and S.~Tanimoto.
\newblock Rational curves on prime {F}ano threefolds of index 1.
\newblock {\em J. Algebraic Geom.}, 30(1):151--188, 2021.

\bibitem[Man90]{Manaresi90}
M.~Manaresi.
\newblock On the jumping conics of a semistable rank two vector bundle on
  {${\bf P}^2$}.
\newblock {\em Manuscripta Math.}, 69(2):133--151, 1990.

\bibitem[Man19]{Mandal}
S.~Mandal.
\newblock On the loci of morphisms from p1 to g(r,n) with fixed splitting type
  of the restricted universal sub-bundle or quotient bundle.
\newblock {\em arXiv:1908.09978 [math.AG]}, 2019.

\bibitem[Mar81]{Maruyama81}
M.~Maruyama.
\newblock The theorem of {G}rauert-{M}\"{u}lich-{S}pindler.
\newblock {\em Math. Ann.}, 255(3):317--333, 1981.

\bibitem[Mar24]{Marangone24}
E.~Marangone.
\newblock The non-lefschetz locus of conics, 2024.

\bibitem[MnOSC12]{MOS12}
R.~Mu\~{n}oz, G.~Occhetta, and L.~E. Sol\'{a}~Conde.
\newblock Uniform vector bundles on {F}ano manifolds and applications.
\newblock {\em J. Reine Angew. Math.}, 664:141--162, 2012.

\bibitem[MR84]{MR84}
V.~B. Mehta and A.~Ramanathan.
\newblock Restriction of stable sheaves and representations of the fundamental
  group.
\newblock {\em Invent. Math.}, 77(1):163--172, 1984.

\bibitem[MR82]{MR82}
V.~B. Mehta and A.~Ramanathan.
\newblock Semistable sheaves on projective varieties and their restriction to
  curves.
\newblock {\em Math. Ann.}, 258(3):213--224, 1981/82.

\bibitem[Nak04]{Nakayama04}
N.~Nakayama.
\newblock {\em Zariski-decomposition and abundance}, volume~14 of {\em MSJ
  Memoirs}.
\newblock Mathematical Society of Japan, Tokyo, 2004.

\bibitem[Neu09]{Neumann10}
S.~Neumann.
\newblock A decomposition of the {M}oving cone of a projective manifold
  according to the {H}arder-{N}arasimhan filtration of the tangent bundle.
\newblock Thesis, {U}niversit{\"a}t Freiburg,
  {https://freidok.uni-freiburg.de/fedora/objects/freidok:7287/datastreams/FILE1/content},
  2009.

\bibitem[Ou23]{Ou23}
W.~Ou.
\newblock On generic nefness of tangent sheaves.
\newblock {\em Math. Z.}, 304(4):Paper No. 58, 23, 2023.

\bibitem[Pey17]{Peyre17}
E.~Peyre.
\newblock Libert\'{e} et accumulation.
\newblock {\em Doc. Math.}, 22:1615--1659, 2017.

\bibitem[PRT20]{PatelRiedlTseng}
A.~Patel, E.~Riedl, and D.~Tseng.
\newblock Moduli of linear slices of high degree hypersurfaces.
\newblock {arXiv:2005.03689 [math.AG]}, 2020.

\bibitem[Ram90]{Ramella}
Luciana Ramella.
\newblock La stratification du sch\'{e}ma de {H}ilbert des courbes rationnelles
  de {${\bf P}^n$} par le fibr\'{e} tangent restreint.
\newblock {\em C. R. Acad. Sci. Paris S\'{e}r. I Math.}, 311(3):181--184, 1990.

\bibitem[Ran01]{Ran01}
Z.~Ran.
\newblock The degree of the divisor of jumping rational curves.
\newblock {\em Q. J. Math.}, 52(3):367--383, 2001.

\bibitem[Ray72]{Raynaud72}
M.~Raynaud.
\newblock Flat modules in algebraic geometry.
\newblock {\em Compositio Math.}, 24:11--31, 1972.

\bibitem[Sch61]{Schwarzenberger61}
R.~L.~E. Schwarzenberger.
\newblock Vector bundles on the projective plane.
\newblock {\em Proc. London Math. Soc. (3)}, 11:623--640, 1961.

\bibitem[Spi79]{Spindler79}
H.~Spindler.
\newblock Der {S}atz von {G}rauert-{M}\"{u}lich f\"{u}r beliebige semistabile
  holomorphe {V}ektorb\"{u}ndel \"{u}ber dem {$n$}-dimensionalen
  komplex-projektiven {R}aum.
\newblock {\em Math. Ann.}, 243(2):131--141, 1979.

\bibitem[{Sta}]{stacks}
The {Stacks Project Authors}.
\newblock {\itshape Stacks Project}.

\bibitem[Sta03]{Starr03}
J.~Starr.
\newblock The {K}odaira dimension of spaces of rational curves on low degree
  hypersurfaces, 2003.

\bibitem[Str84]{Stromme84}
S.~A. Str\o{}mme.
\newblock Ample divisors on fine moduli spaces on the projective plane.
\newblock {\em Math. Z.}, 187(3):405--423, 1984.

\bibitem[Vit04]{Vitter04}
A.~Vitter.
\newblock Restricting semistable bundles on the projective plane to conics.
\newblock {\em Manuscripta Math.}, 114(3):361--383, 2004.

\end{thebibliography}

\end{document}